\documentstyle{amlts}
\begin{document}
\annalsline{151}{2000}
\received{December 29, 1995}
\revised{August 2, 1999}
\startingpage{459}
\def\bye{\end{document}}
 \font\tenrm=cmr10
\catcode`\@=11
\font\twelvemsb=msbm10 scaled 1100
\font\tenmsb=msbm10
\font\ninemsb=msbm10 scaled 800
\newfam\msbfam
\textfont\msbfam=\twelvemsb  \scriptfont\msbfam=\ninemsb
  \scriptscriptfont\msbfam=\ninemsb
\def\msb@{\hexnumber@\msbfam}
\def\Bbb{\relax\ifmmode\let\next\Bbb@\else
 \def\next{\errmessage{Use \string\Bbb\space only in math
mode}}\fi\next}
\def\Bbb@#1{{\Bbb@@{#1}}}
\def\Bbb@@#1{\fam\msbfam#1}
\catcode`\@=12

 \catcode`\@=11
\font\twelveeuf=eufm10 scaled 1100
\font\teneuf=eufm10
\font\nineeuf=eufm7 scaled 1100
\newfam\euffam
\textfont\euffam=\twelveeuf  \scriptfont\euffam=\teneuf
  \scriptscriptfont\euffam=\nineeuf
\def\euf@{\hexnumber@\euffam}
\def\frak{\relax\ifmmode\let\next\frak@\else
 \def\next{\errmessage{Use \string\frak\space only in math
mode}}\fi\next}
\def\frak@#1{{\frak@@{#1}}}
\def\frak@@#1{\fam\euffam#1}
\catcode`\@=12

\def\ritem#1{\item[{\rm #1}]}
\def\RR{{\Bbb R}}
\input boxedeps.tex 
\SetepsfEPSFSpecial 
\HideDisplacementBoxes
\def\figin#1#2{\medbreak
$$
 {\BoxedEPSF{#1 scaled
#2}%
}%
$$
\medbreak\noindent}

\title{Double bubbles minimize}

\acknowledgements{The first author was partially supported by the National Science Foundation.}
\twoauthors{Joel Hass}{Roger Schlafly}

\institutions{University of California, Davis, CA\\
{\eightpoint  http://math.ucdavis.edu/hass}\\
{\eightpoint {\it E-mail address\/}:  hass@math.ucdavis.edu}\vglue6pt
Dept.\ of Computer Science, University of California, Santa Cruz, CA\\
{\eightpoint  http://bbs.cruzio.com/schlafly}\\
{\eightpoint {\it E-mail address\/}:  real@ieee.org}}

\bigbreak \centerline{\bf Abstract}
\medbreak

 The classical isoperimetric inequality in $\RR^3$ states that the surface of smallest
area enclosing a given volume is a sphere. We show that the least area surface enclosing
two equal volumes is a double bubble, a surface made of two pieces of round
spheres separated by a flat disk, meeting along a single circle at an angle of $120^\circ$.
 
\section{Introduction} \label{sec.intro} 

In this paper we find the unique surface of smallest area enclosing two equal volumes.
The surface is called a {\it double bubble}, and is made of two pieces of round
spheres separated by a disk, meeting along a single circle at an angle of $120^\circ$. This is the
form assumed by two equally sized soap bubbles which are brought together until their
boundaries conglomerate to form a common wall. 
See Figure~1, due to J. Sullivan.

Isoperimetric problems, which study maximizing the size of an enclosed bounded region  whose boundary
is of fixed size, are among the oldest problems in mathematics. 
For a broad discussion of isoperimetric problems see Osserman \cite{O}.

The two volume isoperimetric problem in $\RR^3$ was considered by the Belgian
physicist J. Plateau \cite{Pl}, and appears in C.V. Boys' famous book on soap bubbles. As
Boys wrote,

\begin{quote}
``When however the bubble is not single, say two have been blown in real contact
with one another, again the bubbles must together take such a form that the total surface of the
two spherical segments and of the part common to both, which I shall call the interface, is the
smallest possible surface which will contain the two volumes of air and keep them separate."
\end{quote}

\figin{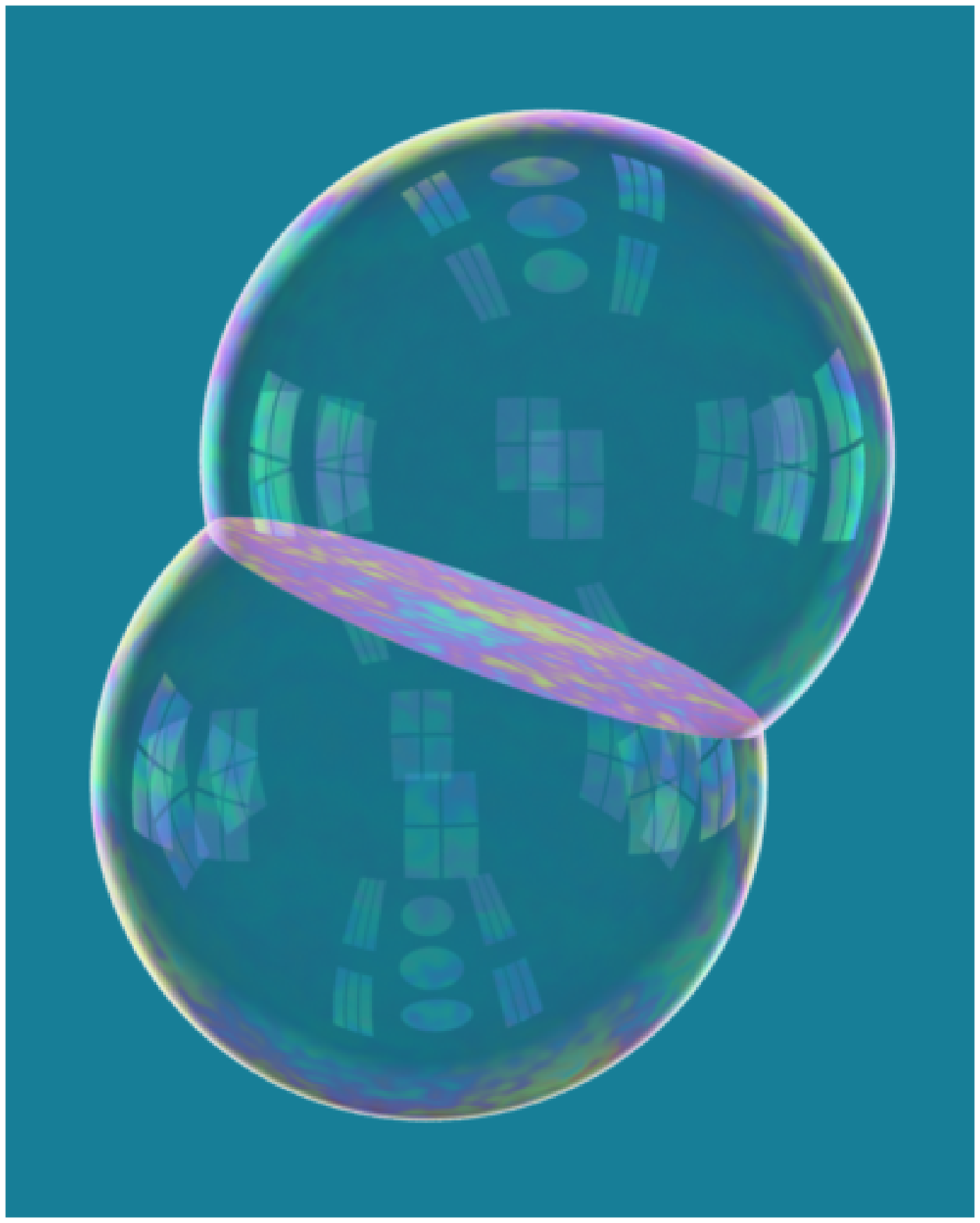}{300}
\centerline{{Figure 1}. A double bubble enclosing two equal volumes.}
\medbreak

We now outline briefly the  recent history of the problem, with more details given later.
F. Almgren proved the existence of bubble clusters enclosing a given collection of volumes in
three-space  and minimizing area among a 
general class of surfaces. J. Taylor established a regularity result for these bubble clusters,
showing that they satisfy regularity properties physically observed by Plateau.
B. White suggested an argument for showing that a two-region minimizer must be a surface of
revolution, and J. Foisy and M. Hutchings filled in the details of this argument.
At this point it was known that the minimizer consisted of pieces of Delaunay surfaces
meeting at $120^\circ$. Hutchings went on to eliminate many remaining configurations, leaving just
the standard double bubble and an additional family known as torus bubbles as the potential minimizers.
In this paper we closely examine the family of torus bubbles and show they cannot be minimizers for
the equal volume double bubble problem. Each torus bubble is determined by an ODE, and we can
rigorously eliminate the entire collection by an examination of the solutions of these
ODEs using a computational technique called interval arithmetic,
described in Section~\ref{Computation}.
At the end, only the standard double bubble remains as a possible minimizer. 

Extensive interest in the problem was generated in recent years by work of F. Morgan. The
planar case has been solved in \cite{ABFHZ} by methods special to two dimensions. 
A fundamental
paper of Almgren \cite{A} established the existence of solutions to a great variety of geometric
minimization problems, including multiple component isoperimetric problems in
$\RR^3$. Despite the fact that existence was established some time ago,
our result is the first explicit example of a
surface in $\RR^3$ solving a multiple region finite volume isoperimetric problem.
Lawlor and Morgan have shown that the cone over the regular tetrahedron is a minimizing
surface among all those with its boundary that separate the four faces of the tetrahedron,
but it is still not known whether this surface is minimizing without the separation condition.

Such multiple region isoperimetric problems arise in many fields, including the
growth and shape of biological cells \cite{M-S-W}, \cite{Th}. There were
extensive studies of such problems by physical and biological scientists in the 19$^{\rm th}$ century
\cite{Kel}, \cite{Pl}. Plateau established experimentally that a soap bubble cluster is a
piecewise-smooth surface having only two types of singularities. The first type of
singularity occurs when three smooth surfaces come together along a smooth triple curve
at an angle of $120^\circ$. The second type of singularity occurs when six smooth surfaces and
four triple curves converge at a point. The angles are then equal to
those of the cone over the 1-skeleton of a regular tetrahedron. A mathematical proof that
these types of singularities are the only ones possible in a minimizing bubble in $\RR^3$ was
given by Taylor \cite{T}. However there were no explicit minimizing bubbles known for any
collection of volumes exhibiting either of these singularities. Thus we have found the
first explicit example of a closed minimizing surface in $\RR^3$ known to exhibit some of the
singularities predicted by Plateau.

The arguments presented in this paper are a mixture of geometrical analysis and
of estimates of geometric quantities obtained by the use of numerical computation. We perform
these calculations with strict estimates on the accuracy of the computations.
Since it is still somewhat unusual in a mathematical proof to use digital computers 
to do calculations involving real numbers, we will say a few words about the nature of
this part of the argument. 

Computers are widely known in the mathematical community to have been used successfully 
for analysis of discrete and combinatorial problems, such as the 4-color theorem
\cite{A-H}, \cite{Seymour}. Problems such as isoperimetric inequalities are qualitatively
different, since they inherently involve real numbers. Real numbers are represented in digital
computers in a floating point format which allows exact description of only a finite number of
rationals. Simple calculations, like division, or even addition, lead to unrepresentable numbers,
and so the computer must round off. For purposes of mathematical proof, the size of the round-off
must be accurately tracked throughout the calculation.

Methods for strictly estimating solutions of differential equations exist, but are
not yet widely used in the mathematical community. The vast majority of numerical
work is for approximation and simulation that does not meet the standards of a mathematical
proof. However there have been some important results achieved through rigorous
use of floating point numerical methods to achieve
traditional mathematical proofs, notably Lanford's work on the Feigenbaum
Conjectures \cite{Lanford} and Fefferman and de la Llave's work on the stability of matter
\cite{delaLlave}, \cite{Fefferman} and the work of McKay and Percival on nonexistence of invariant
tori \cite{MP}. 
See also recent work of Hales on the Kepler Conjecture \cite{Hales} and Gabai, Meyerhoff and Thurston \cite{GMT}
on hyperbolic 3-manifolds.
Computer calculations are essential to our proof that equal volume double bubbles
minimize area. There are too many calculations to be done by hand. 

The proof parameterizes the space of possible solutions by a two-dimensional
rectangle, one dimension corresponding to an angle and the other to a mean curvature. This
rectangle is subdivided into 15,016 smaller rectangles which are investigated by calculations
involving a total of 51,256 numerical integrals.
Every calculation is done with strict
error bounds, and all results are precise mathematical statements.  
All operations conform to the IEEE 754 standard for computer arithmetic,
a widely adopted standardized method for
implementing real number (floating point) computations on computers \cite{ANSI}.
Our methods indicate that numerical techniques are likely to play an important role in future
geometrical arguments.

The main result is the following, announced in \cite{ERA}: 

\specialnumber{1}
\proclaim{Theorem} \label{main.thm}
The unique surface of least area enclosing two equal volumes in $\RR^3$ is a double bubble{\rm .}
\endproclaim

\specialnumber{2}
\proclaim{Corollary}
For any surface in $\RR^3$ enclosing two regions{\rm ,} each of volume $v${\rm ,} the area $a$ satisfies
$$
a^3 \ge 243 \pi v^2 , 
$$
with equality if and only if the surface is a double
bubble enclosing two regions of volume $v${\rm .}
\endproclaim

The first issue in establishing Theorem~\ref{main.thm} is to find an appropriate category of
surface in which to minimize area. 

It suffices in this paper to consider piecewise smooth two-dimensional surfaces.
It is sometimes useful to consider a much more general notion of surface, such as the
$(F, \varepsilon, \delta)$ sets described in \cite{A}.
Our arguments actually imply that the double bubble
minimizes in this larger class; see Theorem~\ref{thm.bubble}.
Consideration of such a larger class of surfaces is needed
primarily in the establishment of the existence and regularity of a minimizer, carried out in
\cite{A} and \cite{T}, and we will not need to be overly concerned with it in this paper.

Define a {\it piecewise-smooth curve} to be an embedded finite union of smooth curves,
with any two either disjoint or having intersection contained in their endpoints.
Define a {\it piecewise-smooth surface} to be an embedded finite union of smooth surfaces with
piecewise-smooth boundary curves, with any two surfaces either disjoint or intersecting along
piecewise-smooth curves contained in their boundaries. The set of points which are not
in the interior of a smooth subsurface of a piecewise-smooth surface is called the
{\it singular set}.

Define a {\it bubble} to be a piecewise-smooth surface satisfying:
\begin{itemize}
\item[1.] Each two dimensional surface has constant mean curvature.
\item[2.] The singular set is of the type described by Plateau. It consists of smooth triple
curves along which three smooth surfaces come together at an angle of $120^\circ$ and
isolated vertices where six smooth surfaces and four triple curves converge at a point.
The angles at the point are equal to those of the cone over the 1-skeleton of a
regular tetrahedron.
\item[3.] The mean curvatures around an edge in the singular set where three surfaces have
common boundary sum up to zero.
\end{itemize}

The above conditions are necessary for no local perturbation of a piecewise-smooth surface to
decrease the area while preserving the volume in each of its complementary
regions. 
 
A bubble enclosing regions of prescribed volumes is called a {\it minimizing bubble}
if it minimizes area among all bubbles enclosing the same volumes. The regions are not
necessarily connected.

Given positive constants $v_1$ and $v_2$, let $a(v_1,v_2)$ denote the infimum of the
area of piecewise-smooth surfaces enclosing two regions $R_1$ and $R_2$ in $\RR^3$ which are
closed bounded sets with disjoint interiors such that ${\rm volume}(R_1)=v_1$ and
${\rm volume}(R_2)=v_2$. We will refer to the interior of the
complement of $R_1 \cup R_2$ as
the exterior region $R_0$. Note that $R_0$ may not be connected, in which case
the bubble encloses some compact ``empty regions''.

A {\it double bubble} enclosing volumes $v_1,v_2$ is a surface made of three pieces of
round spheres,
meeting along a single circle at an angle of $120^\circ$, and enclosing two
connected regions having volumes $v_1$ and $v_2$.
We consider the plane to be a sphere of
infinite radius in this setting,
allowing the interface of a double bubble to be a flat disk, as occurs when
$v_1$ and $v_2$ are equal.
See Figure~2.
In this and all other figures, we assume that the axis of rotational 
symmetry is the $x$-axis.
Furthermore, we will always take the generating curves for a surface of 
revolution to be in the upper half of the $xy$-plane.

It has been conjectured since the work of Plateau that double bubbles give the most
efficient shape for enclosing two given volumes.

\specialnumber{3} \proclaim{{C}onjecture} \label{conj-double.bubble}
The double bubble enclosing volumes $v_1,v_2$ is the unique surface having
area equal to $a(v_1,v_2)${\rm .}
\endproclaim

Theorem~\ref{main.thm} solves this conjecture in the case that $v_1 = v_2$.

Deep results of Almgren and Taylor, summarized in Section~\ref{sec.exist}, imply that for any
two positive numbers $v_1$, $v_2$ there exists a minimizing bubble $S(v_1,v_2)$ in $\RR^3$
which encloses volumes $v_1$ and $v_2$. Arguments based on ideas of B. White and F. Morgan,
and developed in \cite{Fo}, \cite{Hu}, \cite{Mo2}, show that any solution must be a surface of
revolution. We will need to refer to this proof, so we present a simple version
for the case of two regions in $\RR^3$ in Theorem~\ref{thm.revolution}. The lack of such an
argument for isoperimetric problems involving three or more regions in $\RR^3$ makes those
problems more formidable.

While the reduction to a surface of revolution gives an enormous simplification in the scale of
the problem, F. Morgan has pointed out several major topological obstacles to solving the double
bubble conjecture. 
The first is that the regions $R_1$ and $R_2$ bounded by the minimizing bubble $S(v_1,v_2)$
may not be connected, as illustrated in Figure~3.

A second problem is that $S(v_1,v_2)$ may enclose bounded regions in $\RR^3$ which do not
form part of either $R_1$ or $R_2$, as in Figure~4. These regions are
called {\it empty regions} by Morgan.

 \begin{center}
\BoxedEPSF{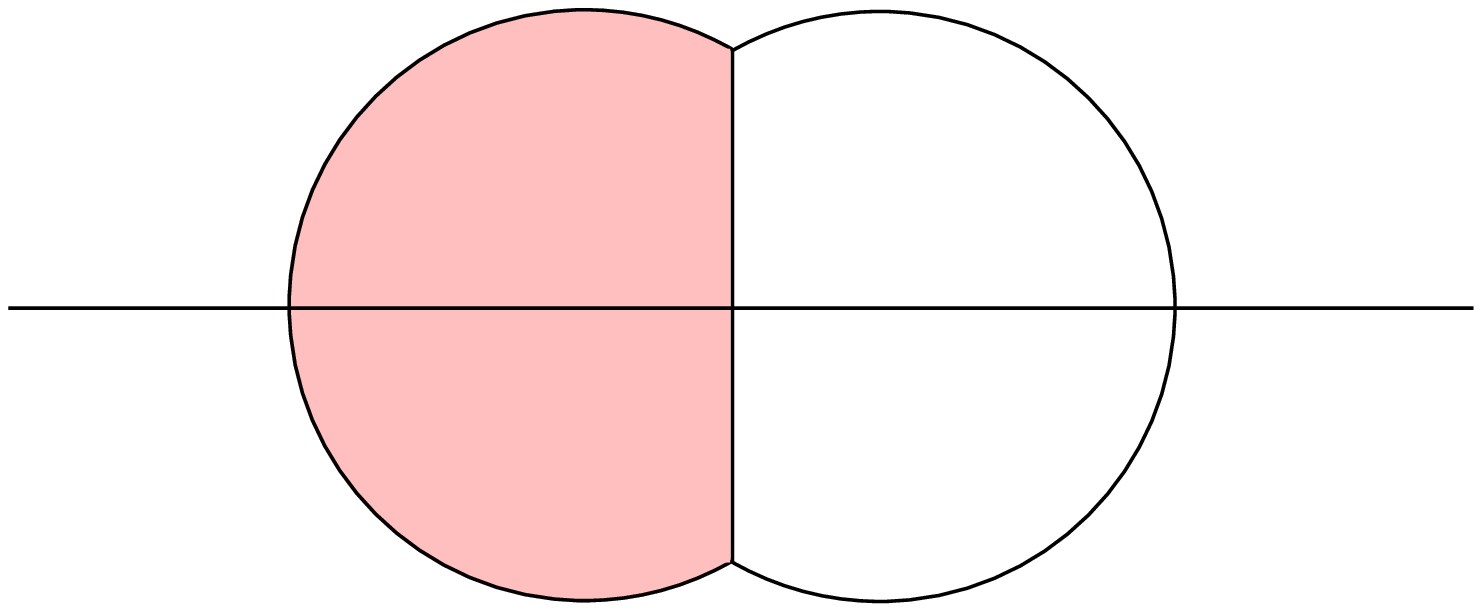 scaled 300}
\quad
\BoxedEPSF{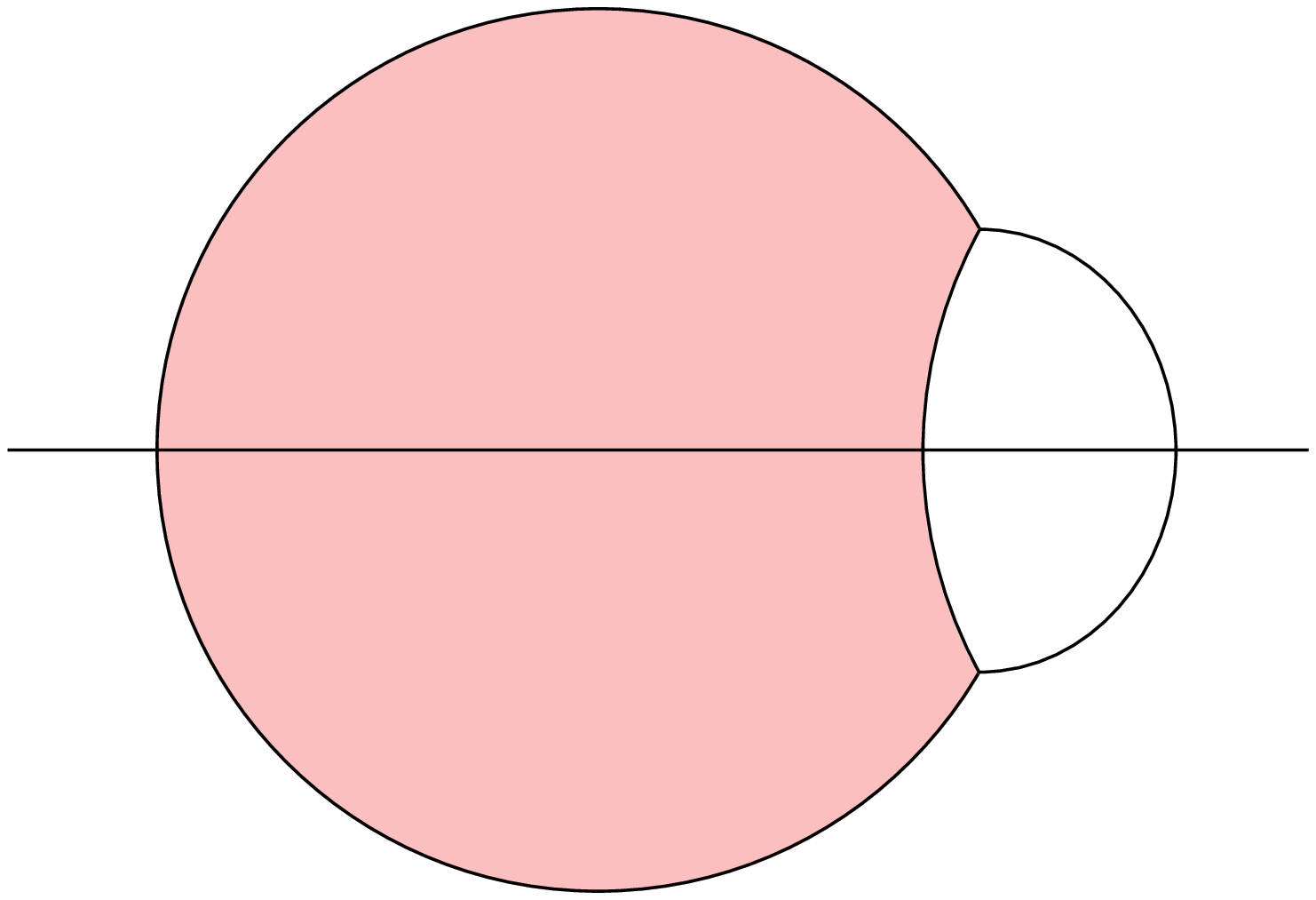 scaled 300}
\end{center}

 {{Figure 2}. Equal and unequal volume double bubbles in cross-section.

The axis
of revolution is the  $x$-axis.}

\figin{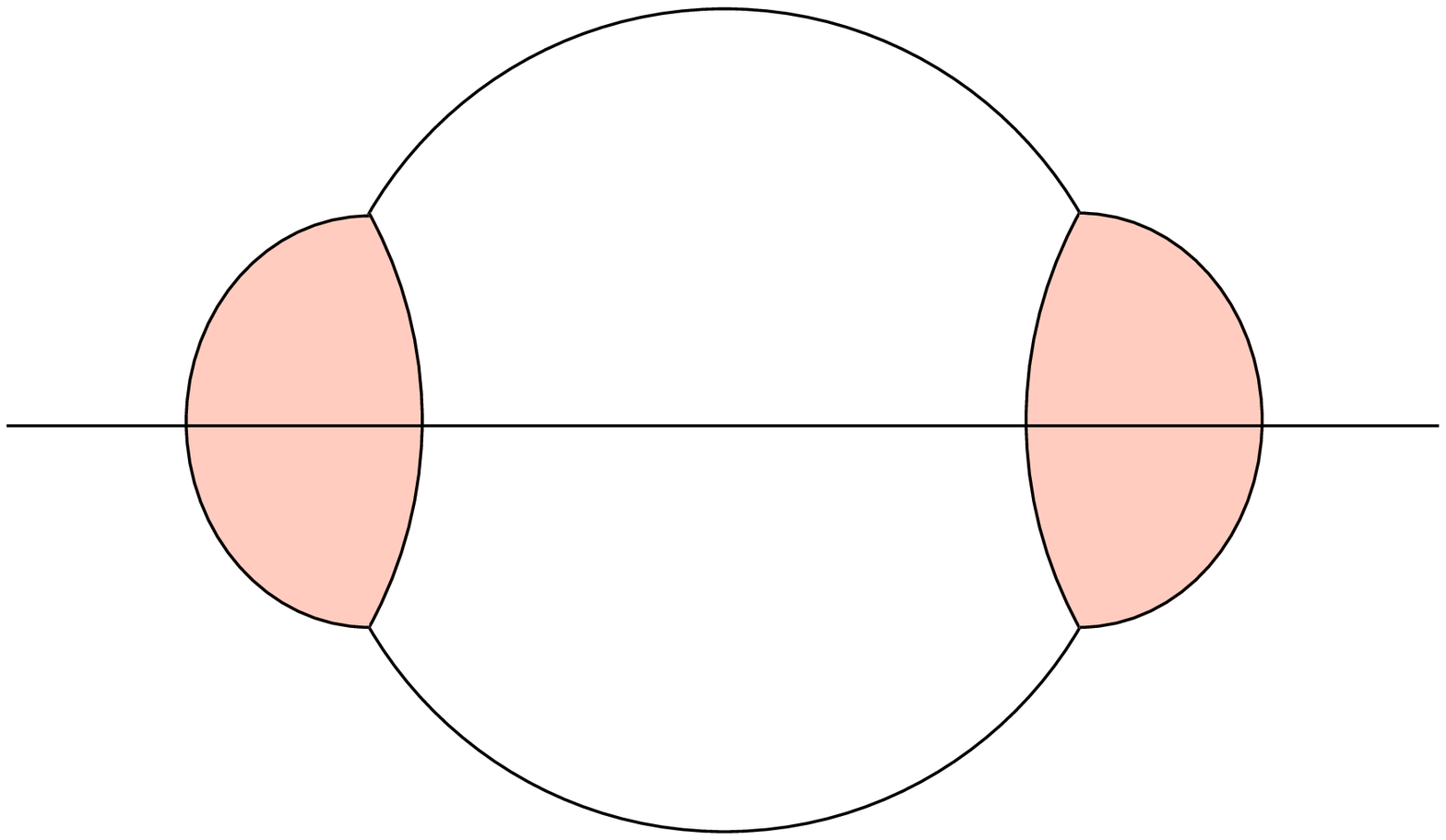}{250}

{{Figure 3}. Cross-section of a bubble with a nonconnected region.  

The
nonconnected region is shaded.}

\figin{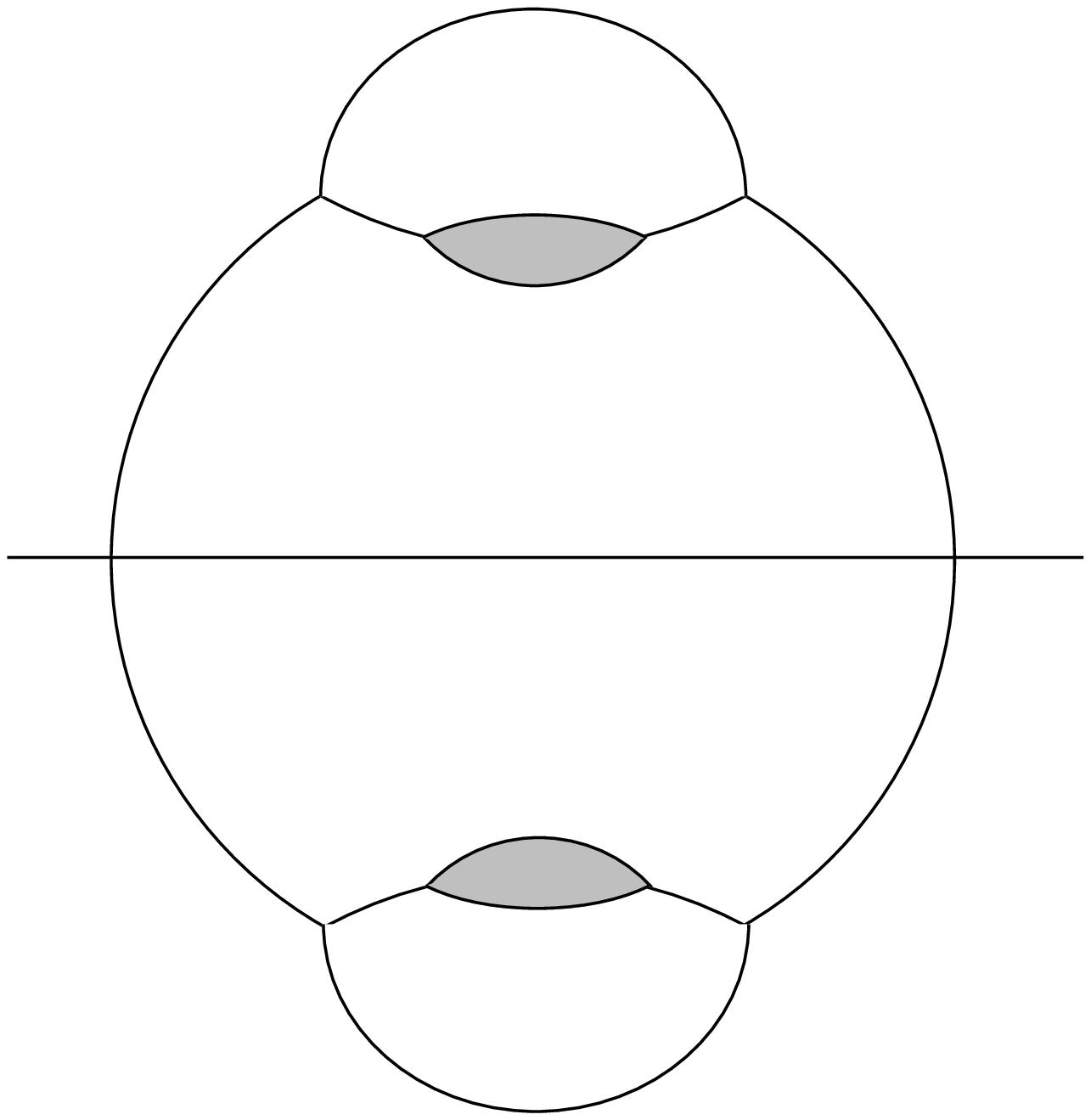}{250}

{{Figure 4}. \enspace Cross-section of a bubble with an empty torus region. 
 The

 shaded 
region generates
a solid torus when revolved around the $x$-axis.}
\medbreak

A third problem is that $S(v_1,v_2)$ may enclose nonsimply connected regions. Since
the surface is a surface of revolution, these regions are homeomorphic to solid tori, and we call
them {\it torus components}.

There are therefore numerous possible configurations which a minimizing bubble might take. A
recent breakthrough due to Hutchings \cite{Hu} has given restrictions on the type of surfaces
that can arise in the double bubble problem, eliminating many possibilities. Some remaining
possibilities are depicted in Figure~5. For the case of equal volumes, Hutchings
showed that there were further constraints. Each region must be connected, leaving only
two possible configurations, the double bubble and an additional class of possibilities,
called {\it torus bubbles}, whose properties will be discussed in Section~\ref{sec.torus}.
\begin{eqnarray*}
\BoxedEPSF{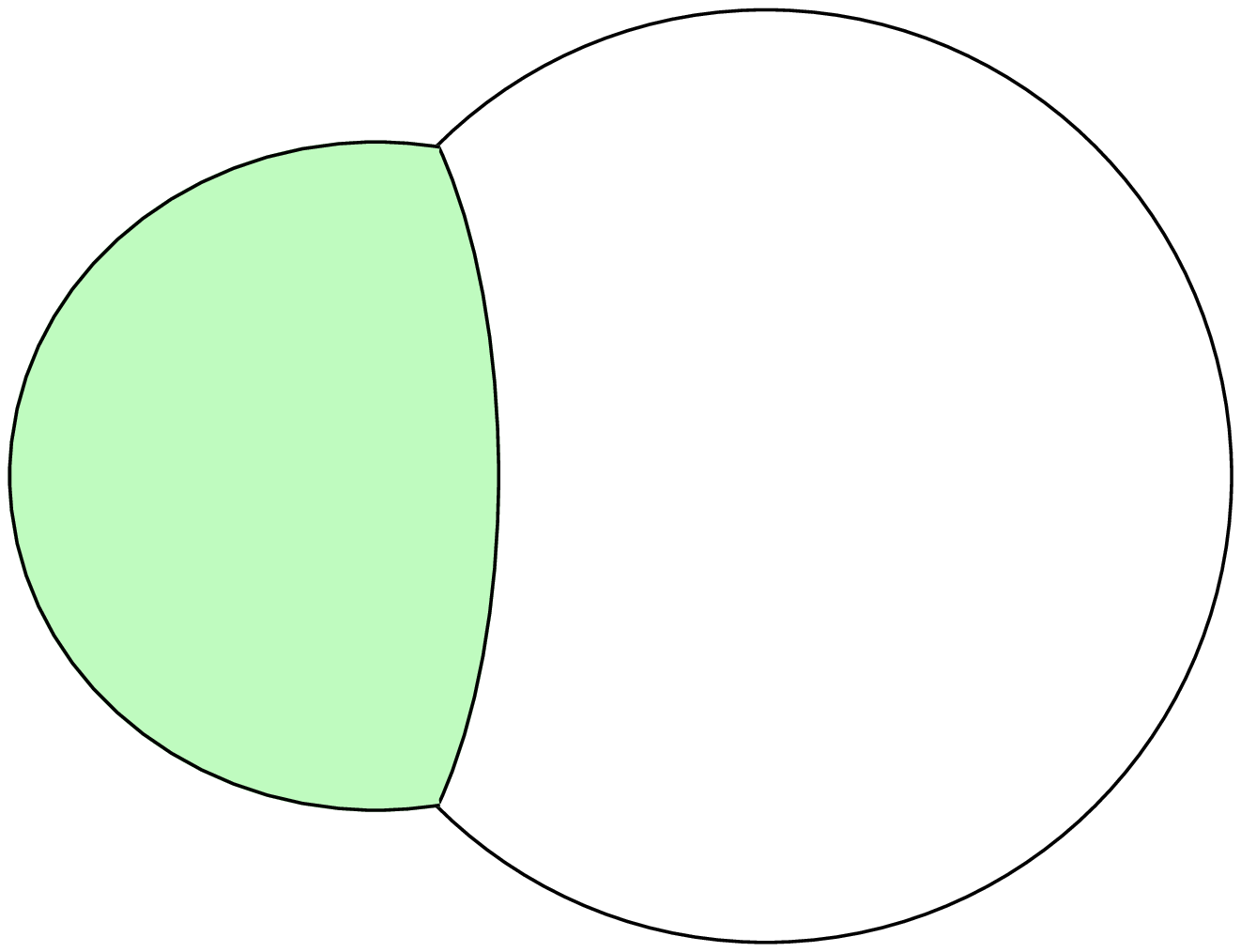 scaled 250} && \BoxedEPSF{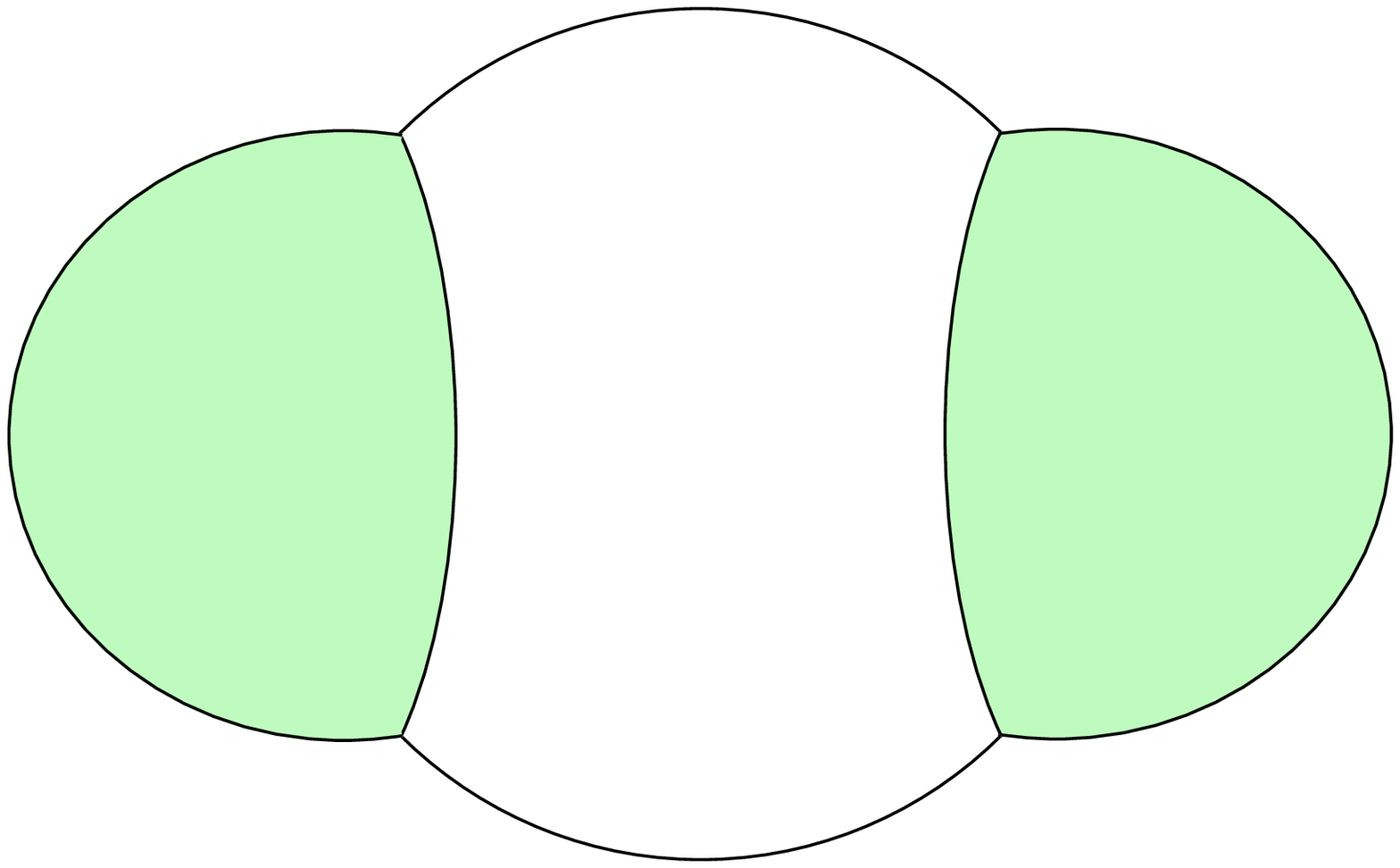 scaled 250} \\
\BoxedEPSF{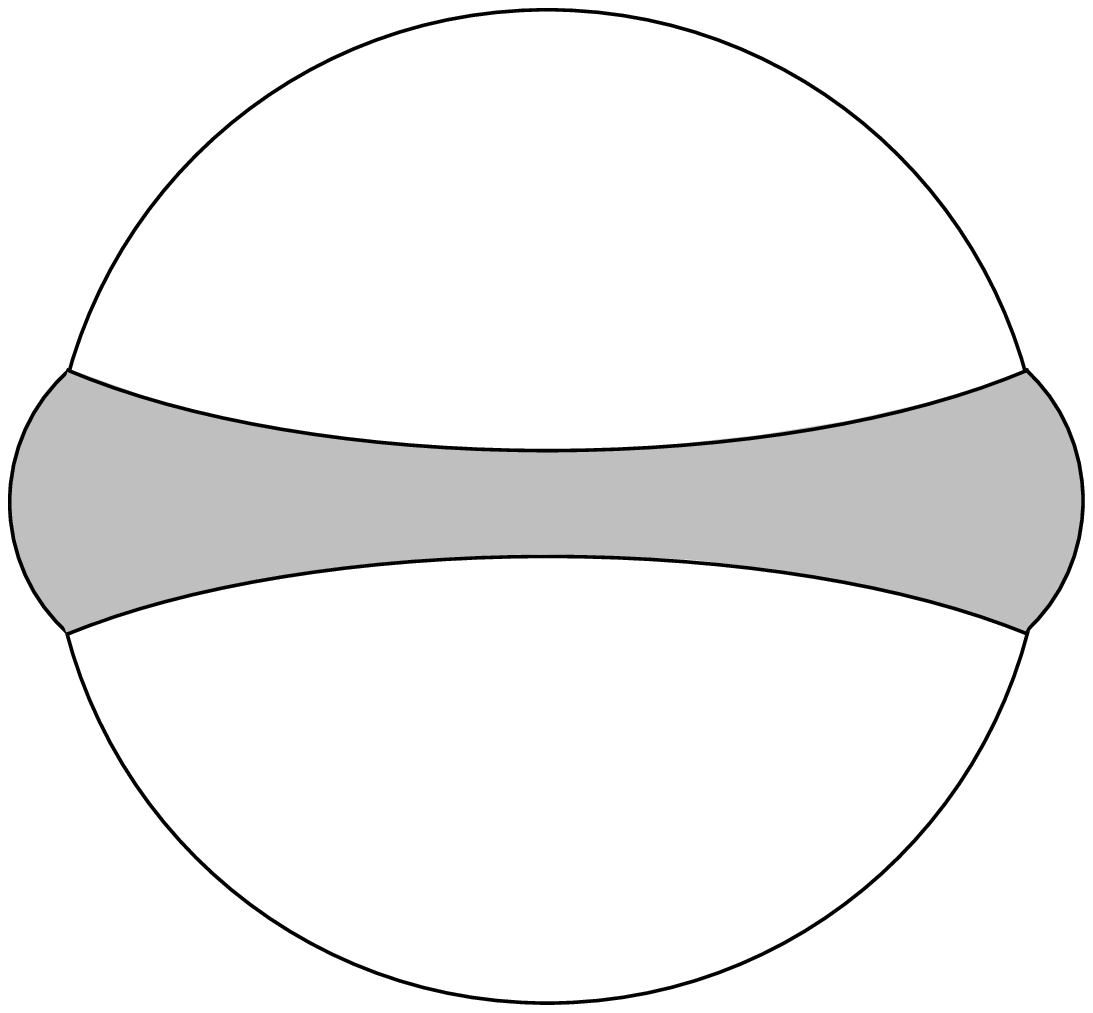 scaled 250} && \BoxedEPSF{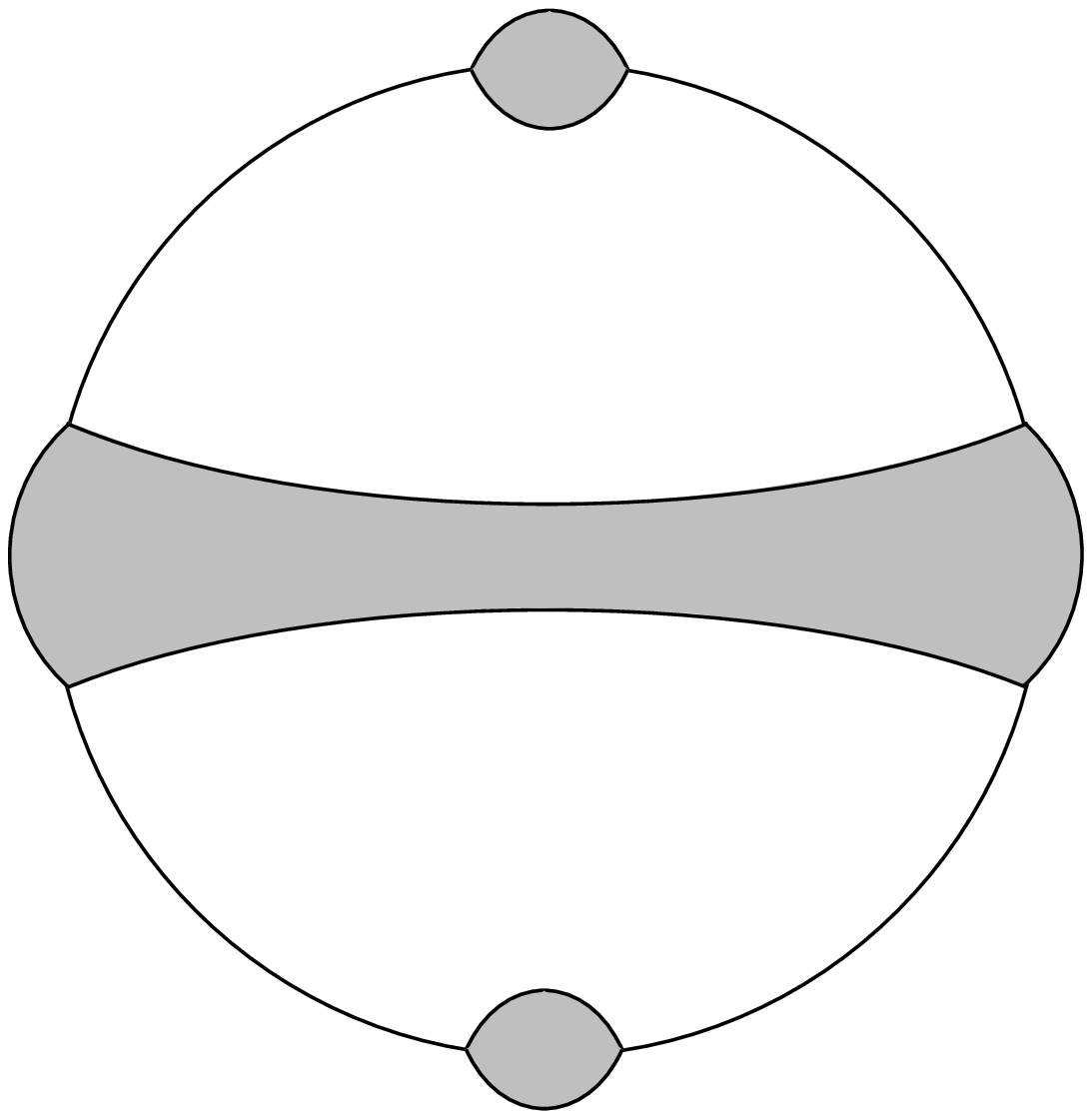 scaled 250} 
\end{eqnarray*} 
\centerline {{Figure 5}. Cross-sections of some possible bubble configurations.}
\medbreak

In Section~\ref{Computation} we describe the algorithm used in a series of computations which
show that torus bubbles are not minimizers for the two equal volume isoperimetric
problem in $\RR^3$.

Techniques extending those developed in this paper also prove
Conjecture~\ref{conj-double.bubble} for other volume ratios, but at this time do
not seem to suffice to cover all cases.

\section{Existence and regularity} \label{sec.exist} 

Almgren showed in \cite{A} that there exists an area minimizing surface $S(v_1,v_2)$ in $\RR^3$ among the
set of surfaces enclosing volumes $v_1,v_2$. Here {\it surface} refers to a generalized
notion defined using the methods of geometric measure theory, specifically what Almgren
calls $(F, \varepsilon, \delta)$ sets. For our purposes it suffices that
this class includes the piecewise-smooth surfaces.

\proclaimtitle{\cite{A}} 
\proclaim{Theorem}\label{thm.exist} 
An area minimizing surface $S(v_1,v_2)$ exists and is a smooth surface almost everywhere{\rm . }
\endproclaim

The nature of the singularities of $S(v_1,v_2)$ was established by Taylor. 

\proclaimtitle{\cite{T}}
\proclaim{Theorem}  \label{thm.bubble} 
$S(v_1,v_2)$ is a piecewise\/{\rm -}\/smooth surface{\rm .} Its singularities consist of smooth triple curves
along which three smooth surfaces come together at an angle of $120^\circ$ and isolated points
where four triple curves and six pieces of surface converge{\rm .} At these isolated points the
asymptotic cone is the cone over the {\rm 1-}\/skeleton of a regular tetrahedron{\rm . }
\endproclaim

Thus Taylor's work established that $S(v_1,v_2)$, which Almgren's theorem established
as a minimizer among a generalized class of surfaces called 
$(F, \varepsilon, \delta)$ surfaces, is a bubble, in our terminology.
It follows that if we can establish the minimality of the double bubble among our class of bubbles then
we will also show it minimizes among the more general class of surfaces considered by
Almgren. We now establish some properties of minimizing bubbles. Lemma~\ref{curv.sum},
Corollary~\ref{curv.sum.cor} and Lemma~\ref{lem.connected},
are standard results in variational geometry. We include brief proofs for completeness.
Lemma~\ref{curv.sum} states that the sum of the (oriented) mean curvatures of all points
crossed by a path starting and ending in the same region is zero.
It implies that the sum of the (oriented) mean curvatures of all points
crossed by any path depend only on its starting and ending points.
These results hold even if the path crosses through many surfaces, and if
it starts and ends in different components of the same region.

\proclaim{Lemma} \label{curv.sum}
Let $\gamma$ be an oriented curve in $\RR^3$
intersecting a minimizing bubble $B$ transversely
at regular points such that the initial and final points of
$\gamma$ lie in the interior of the same region{\rm .}
Then the sum of the mean curvatures of all the points of $\gamma \cap B${\rm ,} oriented by
$\gamma${\rm ,} is zero{\rm .} 
\endproclaim

\demo{Proof} Perturb the curve slightly so that each of its intersections with $B$
becomes perpendicular. Consider a deformation of
$B$ which pushes points in $\RR^3$ near $\gamma$ a uniform distance along the curve. To first
order this preserves the volume of each region. The derivative of the area, to first order, is
given by the sum of the mean curvatures over the points of $\gamma \cap B$.
If this sum is nonzero, a deformation can be defined which decreases
area while preserving the volume of each region. \enddemo

A special case of the above lemma occurs when $\gamma$ is a simple closed curve encircling a
triple curve of the bubble. The lemma then implies that the sum of the mean
curvatures around
the triple curve adds up to zero. This local minimization condition is built into our definition of a
bubble. More generally, the lemma implies that any two surfaces separating the same pair of
regions have the same mean curvature.

\proclaim{{C}orollary} \label{curv.sum.cor}
If $S_1$ and $S_2$ are surfaces in a bubble{\rm ,} each of which
separate regions $R_1$ and $R_2${\rm ,} then their mean curvatures are equal{\rm . }
\endproclaim
\demo{Proof} Consider an arc $\alpha$ which starts in $R_1$ passes through $S_1$ into $R_2$,
then transversely through the bubble and
on through $S_2$ back into $R_1$. A subarc $\beta$ of $\alpha$ starts and ends in
$R_2$. Applying Lemma~\ref{curv.sum} to both $\alpha$ and $\beta$ implies the corollary.
\enddemo

\proclaim{Lemma} \label{lem.connected}
A minimizing bubble $B$ is connected{\rm .}
\endproclaim

{\it Proof}. If not, we can move a component by an isometry of $\RR^3$ until it touches a distinct
component. The resulting singularity violates those allowed in Theorem~\ref{thm.bubble}. 
\hfill\qed \medbreak

The following result was first observed by B. White and F. Morgan, and first written down by Foisy
in \cite{Fo}.  Since we will
need to refer to the proof, we present an argument for the case of two regions in $\RR^3$.
More detailed and general arguments can be found in the work of Hutchings \cite{Hu} and Morgan \cite{Mo2}.

\proclaim{Theorem} \label{thm.revolution}
$S(v_1,v_2)$ is a surface of revolution{\rm . }
\endproclaim

{\it Proof}. Since $S(v_1,v_2)$ is a minimizing bubble,
each of its faces separates a pair of distinct
regions - either separating $R_1$ from $R_2$ or one of them from $R_0$ (which is
not connected if there are empty regions.)

Given any unit vector $Z$ in $\RR^3$, there is a plane $P_Z$
in $\RR^3$ perpendicular to $Z$ which bisects the total volume $ v_1 + v_2$ enclosed by $S(v_1,v_2)$.
If $Q_Z$ is a parallel plane which also bisects the total volume, then $P_Z$ and
$Q_Z$ split the regions $R_1$ and $R_2$ in the same proportions.
(In fact $P_Z$ must coincide with~$Q_Z$, but we do not need to show this.
It is also irrelevant to this argument
whether $S(v_1,v_2)$ contains empty regions.)
Let $f(Z)$ denote the proportion of
the volume of $R_1$ on the side of $P_Z$ to which $Z$ points.
The function $f(Z)$ is well defined, continuous on the unit 2-sphere, $0 \leq f(Z) \leq 1$ and
$f(-Z) = 1 - f(Z)$. 
Along any great circle of directions on the unit 2-sphere,
the Intermediate Value Theorem implies that
there are at least two points where $f(Z)=1/2$ and the plane $P_Z$ bisects the volume of
both regions. We fix $Z$ to be such a vector. 

Consider the intersection of $S(v_1,v_2)$ with each of the two half-spaces determined by $P_Z$.
The intersection of a face of $S(v_1,v_2)$ with the plane $P_Z$ is given by a graph in $P_Z$ with
smooth edges and isolated vertices, or else the face is completely contained in $P_Z$,
since the face is a smooth constant mean curvature surface \cite{Gulliver}.
Reflection of the smaller area piece of $S(v_1,v_2)$ lying in one of these half-spaces gives a
new surface $S_1(v_1,v_2)$ which encloses regions of the same volumes.
Since $S_1(v_1,v_2)$
cannot have less area, it must have the same area as $S(v_1,v_2)$. $S_1(v_1,v_2)$ has the
property that reflection through $P_Z$ preserves the new regions, which we continue to call
$R_0$, $R_1$ and $R_2$.
If the plane $P_Z$ contained a face of $S(v_1,v_2)$ then this face would now separate
two components of the same region in $S_1(v_1,v_2)$, and area could be reduced by
removing that face from $S_1(v_1,v_2)$, preserving
the volume of each region. So $P_Z$ cannot contain a face of $S(v_1,v_2)$.

We repeat the above argument 
along a great circle of the unit 2-sphere consisting of the unit vectors which are
perpendicular to $Z$, 
to find a plane $P_W$ perpendicular to $P_Z$ which bisects the volume of
the regions $R_1$ and $R_2$ in $S_1(v_1,v_2)$. This process constructs
a surface $S_2(v_1,v_2)$, still having the same area as
$S(v_1,v_2)$, enclosing regions of the same volumes,
and for which reflection through each of $P_Z$ and $P_W$ preserves $R_1$ and $R_2$.

Since composing reflections through two perpendicular planes gives a rotation of angle
$180^\circ$, it follows that rotation of angle $180^\circ$ about $L = P_Z \cap P_W$ preserves the
regions $R_1$ and $R_2$.
Now consider any plane $Q$ containing $L$.
Note that half of the volume of $R_1$ and $R_2$ lies on each side of $Q$.
As before the plane $Q$ cannot contain a face of $S_2(v_1,v_2)$.
If the intersection of
$Q$ with $S_2(v_1,v_2)$ is not perpendicular to $Q$, then replacing half of $S_2(v_1,v_2)$
with its reflected image through $Q$ gives a new minimizing surface which is not smooth, and
has singularities
not allowed by Theorem~\ref{thm.bubble}. Thus $S_2(v_1,v_2)$ is perpendicular
to each plane through $L$ and is therefore a surface of revolution around $L$. 
$S_1(v_1,v_2)$ and $S_2(v_1,v_2)$ may not a priori be identical,
but they coincide on a half-space
of $\RR^3$. Moreover this half-space was chosen arbitrarily, so that $S_1(v_1,v_2)$ is cut
by $P_W$ into two pieces, each half of some surface of revolution. The axis of each of
these surfaces of revolution is $L = P_Z \cap P_W$, so they coincide and
$S_1(v_1,v_2)$ itself is a surface of revolution. Similar reasoning shows that
$S(v_1,v_2)$ is a surface of revolution, and coincides with $S_2(v_1,v_2)$
everywhere. If we make $S_1$ by reflecting the other half of $S$ across $P_Z$, then the same
argument shows that the other half of $S$ is also a surface of revolution about a line
$L'$ in $P_Z$. The Almgren and Taylor regularity results, and in particular the unique continuation
property of constant mean curvature surfaces, imply that $L'=L$, so $S$ itself is a surface
of revolution about $L$. \hfill\qed\medbreak

We next summarize some key results obtained by Hutchings in \cite{Hu}. Recall
that $S (v_1,v_2)$ is a minimizing bubble that separates $\RR^3$ into regions $R_1$ of volume
$v_1$ and $R_2$ of volume $v_2$ and that $a(v_1,v_2)$ is the area of $S (v_1,v_2)$.

\proclaim{Theorem} \label{Hutchings} The function $a(v_1,v_2)$ is strictly concave on
$[0,+\infty)\times [0,+\infty)${\rm .}
\endproclaim

\proclaim{{C}orollary} \label{area.increases}
If $v_1' > v_1$ then $a(v_1' , v_2) > a(v_1,v_2)${\rm .}
\endproclaim

{\it Proof}. As $v_1 \to +\infty$, $a(v_1, v_2) \ge a(v_1) \to +\infty$.
Concavity of the function $a$ implies the corollary. \hfill\qed

\proclaim{{C}orollary} \label{no.empty}
$S (v_1 ,v_2 )$ has no empty regions{\rm .}
\endproclaim

{\it Proof}. If there is an empty region, join it it to
 one of the other regions by removing a face on
its boundary  and apply the previous corollary.  \hfill\qed
\medbreak

Deeper connectedness results also follow from concavity.
In particular, Hutchings deduced that each of $R_1$
and $R_2$ is connected if the volumes are equal.

\proclaim{Theorem} \label{connected}
$S (v_1 , v_1 )$ encloses exactly two connected components{\rm .}
\endproclaim

{\it Proof}. See Hutchings \cite[Th.\ 4.2]{Hu}.  \hfill\qed

\section{Delaunay surfaces} \label{sec:Delaunay}

In this section, we summarize the classification theory of constant
mean curvature surfaces of
revolution, the Delaunay surfaces, and present some properties of these surfaces that will be
used in our study of bubbles. 

The {\it mean curvature} of a surface in $\RR^3$ is the trace of the
second fundamental form of the surface.
Computing this requires choosing a unit normal vector field to the surface.
The mean curvature vector field, formed by scaling a unit normal vector field
by the mean curvature, gives the direction
of a variation which decreases the area of a surface
as quickly as possible \cite{Spivak}.

For a surface of revolution we have simple formulas for the mean curvature.
Consider a surface of revolution about the $x$-axis in $\RR^3$,
with generating curve contained in the upper half-plane.
If the generating curve is a graph, $y = y(x)$, then
the mean curvature $h$,
equal to the sum of the principal curvatures, is given by the formulas:
$$
k_m = \frac {- \ddot {y}} {(1 + { \dot {y}}^2) ^{3 /2} } \ ,
$$
$$
k_p = \frac{1}{ y \sqrt{1 + \dot {y}^2 }} \ ,
$$
$$
h= k_m + k_p ,
$$ 
where $k_m$ is the curvature at $(x,y)$ of the generating curve in the $xy$-plane
(sometimes called the meridian, or profile curve)
and $k_p$ is the normal curvature of the parallel curve.
The value of $k_p$ is equal to the reciprocal of the distance to the $x$-axis
along the perpendicular to the generating curve \cite{Spivak}. 
Note that the mean curvature of a unit sphere equals two. 
For surfaces of revolution generated by curves
which are not graphs, we need to specify an orientation on the
generating curve to fix the sign of the mean curvature.
We use the convention that the
sign of the mean curvature is given by the formulas above if the curve is oriented left-to-right,
or equivalently if the graph's tangent vector has a positive $x$-component.
If its tangent vector has a negative $x$-component,
the signs in the formulas for $k_m$ and $k_p$
should be reversed in the above formulas.
For curves with vertical tangent vectors the above formulas do not apply directly, but in
our applications vertical tangents are isolated points, and the limits of the above
formulas will give  $k_m$ and $k_p$ at vertical tangencies.

Surfaces of revolution having constant mean curvature were first studied by Euler, and
classified by Delaunay \cite{D}.
The term {\it Delaunay surface} is used to refer to these surfaces.

\proclaim{Theorem} \label{thm.Delaunay}
The surfaces of revolution of constant mean curvature are the plane{\rm ,}
sphere{\rm ,} catenoid{\rm ,} cylinder{\rm ,}
unduloid{\rm ,} and nodoid{\rm .} The generating curves of a nodoid and an unduloid{\rm ,}
called the {\rm nodary} and
the {\rm undulary,} are periodic along the
$x$\/{\rm -}\/axis{\rm ,} and have exactly one local minimum and one local maximum in each period{\rm .} The
undulary is a graph over the $x$\/{\rm -}\/axis{\rm .} The nodary
has one local maximum{\rm ,} one local minimum and two vertical tangencies
in each period{\rm . }
\endproclaim

\demo{Proof} Expositions of the classification of Delaunay surfaces can be found in
\cite{Eells} and \cite{K}. \enddemo

The nodoid and unduloid are not as well known as the other surfaces.
An undulary can be obtained by rolling an ellipse along the $x$-axis
and tracing the path taken by one of the foci.
A nodary can be obtained by rolling a hyperbola \cite{Eells}.
See Figure~6.

$$
\BoxedEPSF{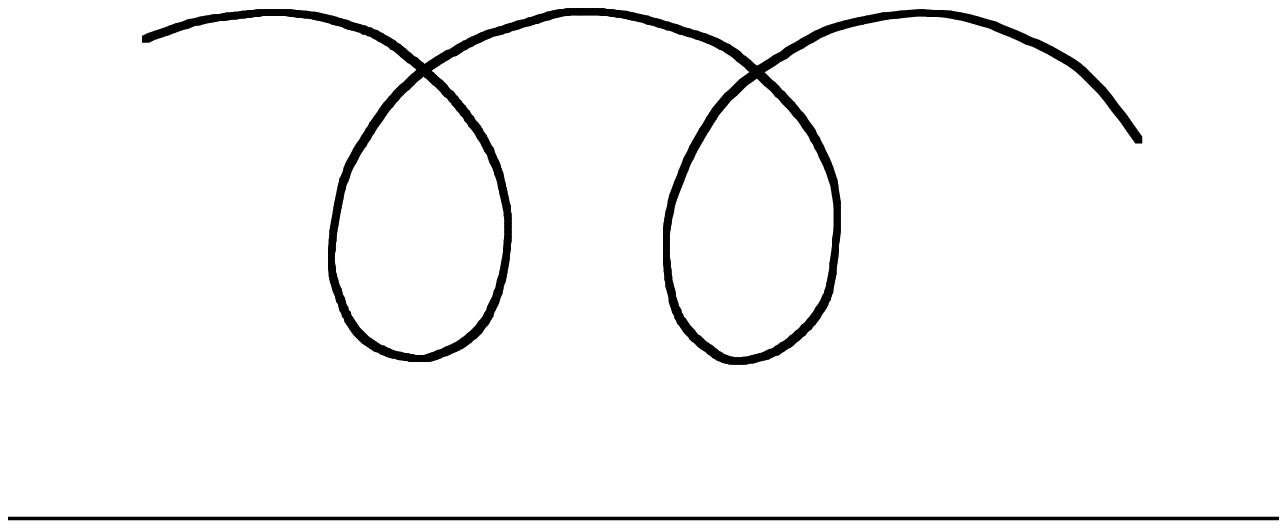 scaled 400} \qquad
\BoxedEPSF{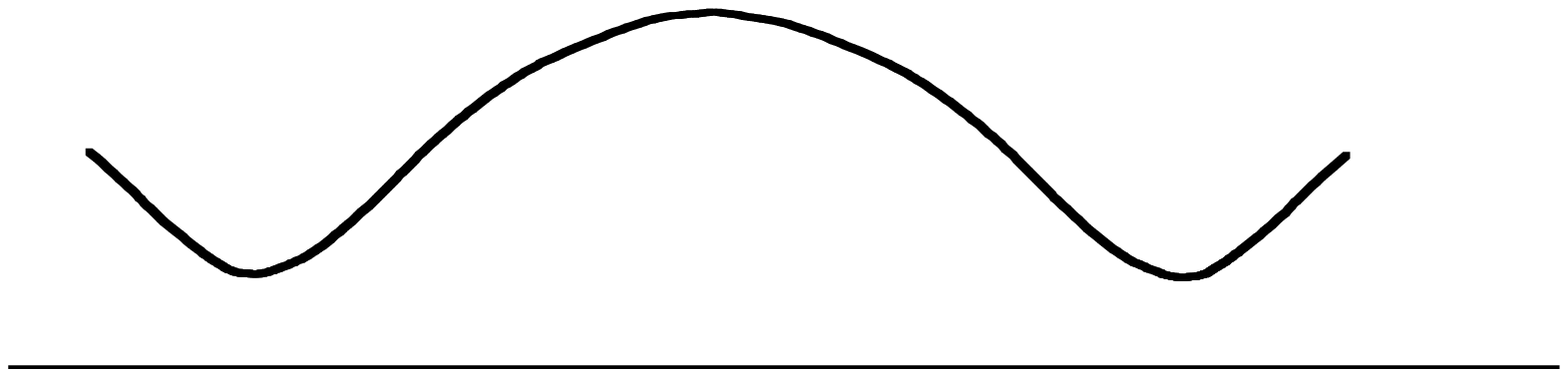 scaled 400}
$$

{{Figure 6}. Portions of a nodary and an undulary, generating curves of a 

nodoid and an
unduloid.}
\medbreak

With their generating curves oriented left to right,
the mean curvature of each of the Delaunay surfaces is nonnegative, with the
exception of the nodoid. The sign of the
mean curvature of a nodoid depends upon the choice
of orientation for its generating nodary, and can be positive or negative.

If a surface of revolution minimizes area among all surfaces of revolution surrounding a
given volume, then its generating curve satisfies an associated Euler-Lagrange equation. This
second order ODE implies that the mean curvature is constant \cite{Eells}. 

The mean curvature of the surface of revolution generated by a $C^2$ graph $y(x)$ is given by:
\begin{equation}
 h = k_m + k_p = \frac { - \ddot {y}} {(1 + \dot {y}^2)^{3/2}} + 
\frac{1}{ y \sqrt{1 + \dot {y}^2}} \ . \label{eq:cmc}
\end{equation}

We can integrate this ODE once, after multiplying through by $2 y \dot {y}$. 

\begin{eqnarray*}
0 & = & h + \frac { \ddot {y}} {(1 + \dot {y}^2)^{3/2}} - \frac{1}{ y \sqrt{1 + \dot {y}^2}} \\
& = & 2hy \dot {y} + \frac {2y\dot {y} \ddot {y}} {(1 + \dot {y}^2)^{3/2}} - \frac{2 \dot {y}}{
\sqrt{1 + \dot {y}^2}} \\
& = & \frac{df(x)}{dx}
\end{eqnarray*}
where $ f(x)$ is the function
\begin{equation}
f(x) = h y^2 - \frac {2y} { \sqrt{1+{\dot {y}}^2} } . \label{eq:force}
\end{equation}
Therefore $f$ is constant along the graph. The fact that the equation for constant
mean curvature surfaces of revolution has a first integral was known to Plateau \cite[pp.~138--139]{Pl},
 who in turn
references work of Beer. Korevaar, Kusner and Solomon called
$f$ the {\it force} of a constant mean curvature surface, and pointed out that it has a physical
interpretation as the net force exerted by a soap film on a plane cutting off an end of the
surface. They also showed that it could be defined for more general constant
mean curvature surfaces \cite{K-K-S}, where it serves as a useful analytical tool.

A convenient way of expressing the force of a constant mean curvature surface of revolution
which holds even when the generating curve is not a graph is given by the formula
\begin{equation}
f = h y^2 - 2y \cos \alpha \ .  \label{cos alpha}
\end{equation}
Here $\alpha$ denotes the angle between the positive $x$-axis and the generating curve of the
surface of revolution, which is an oriented curve.
This makes sense for all values of $\alpha$.
In our applications, we alway orient the
generating curve of a Delaunay surface, and thus fix an unambiguous sign for the mean curvature.

\proclaim{Lemma} \label{lem:trap}
Given $h${\rm ,} $x_0${\rm ,} $ y_0>0$ and $\dot {y_0}${\rm ,} there exists a unique Delaunay curve
$\delta$ through the point $(x_0,y_0)$ having slope $\dot {y_0}$ and generating a Delaunay
surface $D$ having mean curvature $h${\rm .} The slope $\dot {y_0}$ can be infinite{\rm .} The Delaunay
surface is a sphere $S$ if and only if
\begin{equation}
h = h_S = \frac{2}{ y_0 \sqrt{1 + {\dot {y}_0}^2}}\ . \label{eq:hs}
\end {equation} 

If $h > h_S${\rm ,} then $D$ is a nodoid and $\delta$ lies strictly below the circle
$\sigma$ generating $S$ in a neighborhood of $(x_0,y_0)${\rm .} An arc of $\delta$ which is
decreasing as it leaves $(x_0,y_0)$ remains beneath $\sigma$ as long as it remains a graph 
over the $x$\/{\rm -}\/axis{\rm .}

If $ h < h_S${\rm ,} $\delta$ lies strictly above $\sigma$ in a neighborhood of
$(x_0,y_0)${\rm .} If $0 < h < h_S$ then $D$ is an unduloid or cylinder{\rm ,} if
$h< 0$ then $D$ is a nodoid{\rm ,} and if $h=0$ then $D$ is a catenoid{\rm .}
\endproclaim

\demo{Proof} Solving equation \ref{eq:cmc} for $ \ddot {y}$ gives a second order ODE satisfied by
the generating curves of Delaunay surfaces at points where $\dot {y}$ is finite:
\begin{equation}
\ddot {y}= - h(1 + \dot {y}^2)^{3/2} + \frac {1 + \dot {y} ^2}{y}\ . \label{eq:y''}
\end{equation}
Given a constant $h$ and a choice of initial conditions $y$ and $\dot {y}$ this equation has a
unique solution from standard existence and uniqueness of solutions to ODEs.
These solutions
can also be obtained by rolling an appropriate ellipse, hyperbola or parabola \cite{Eells}.
The latter approach also shows that we can find a nodoid passing through $y$ with any mean
curvature $h \neq 0$ in the case where $\dot {y}$ is infinite.

The mean curvature of the unique sphere centered on the $x$-axis and having slope
$\dot {y_0}$ at $(x_0,y_0)$ is given by $h_S$. Assume now that $\sigma$ is centered at $(0,0)$
on the $xy$-plane. 

If $h>h_S$ then Equation~\ref{eq:y''} implies that $\delta$ lies strictly below
$\sigma$ in a neighborhood of $(x_0,y_0)$. If $D$ is an unduloid then $\delta$ is an
undulary whose graph
crosses $\sigma$ again in at least two more points $(x_1,y_1)$ and $(x_2,y_2)$.
We can assume, after performing a left/right reflection if necessary,
that $ x_0 \ge 0$ and $x_1 > x_0$.
Denote by $\gamma$ the subcurve of
$\delta$ starting at $(x_0,y_0)$ and running to $(x_1,y_1)$.
Let $\gamma_t = \gamma + (t,0)$
be a horizontal translate of $\gamma$ by a distance $t$ and let $T = \sup \{ t : \gamma_t \cap S
\ne \emptyset \}$. Then $\gamma_T$ is tangent to $S$ at a point $P$ and lies above $\sigma$
in a neighborhood of $P$. This
violates Equation~\ref{eq:y''},
so $D$ cannot be
an unduloid. $D$ cannot be a cylinder or a catenoid, so it can only be a nodoid.
If a subarc $\gamma$ of $\delta$ which is decreasing as it leaves $(x_0,y_0)$ intersects $\sigma$ 
while $\delta$ remains a graph
then an identical argument shows that we can horizontally
translate $\gamma$ until it is tangent to $\sigma$ while lying above it, which
is impossible by  equation~\ref{eq:y''}.
(In fact it will remain beneath $\sigma$ for somewhat longer than this, but this
will not be relevant for us.) 

If $0<h<h_S$ then equation~\ref{eq:y''} implies that $\delta$ lies
strictly above $\sigma$ in a neighborhood of $(x_0,y_0)$. If $D$ is a nodoid, consider the
subarc $\beta$ of $\delta$ decreasing from $(x_0,y_0)$ until it reaches a
vertical tangency $(x_v,y_v)$. The arc $\beta$ is a graph which may or may not cross $\sigma$ at
an additional point. In either case we can translate $\beta$ horizontally to the left until the
last time at which it intersects $\sigma$. The final intersection must be at an interior point of
$\beta$, so after translation $\beta$ becomes tangent to $\sigma$ at an interior point while
lying underneath $\sigma$, again
impossible by  equation~\ref{eq:y''}.
Thus the surface must be an unduloid or cylinder.

If $h < 0 $ then the surface must be a nodoid. \enddemo

Nodoids and unduloids can be distinguished by the sign of the product $fh$. 

\proclaim{{C}orollary} \label{cor:signhf}
For a nodoid{\rm ,} $fh > 0${\rm .} For an unduloid and a
cylinder $fh < 0${\rm .} For a sphere{\rm ,} a catenoid and a plane{\rm ,} $fh = 0${\rm .} 
\endproclaim

\demo{Proof} For a sphere of radius $r$, we can compute $f$ at a maximum point where
$$f = hr^2 - 2r.$$ 
Since $h = 2/r$, we have $f = 0$. We also compute $fh$ for an unduloid or cylinder at
a local maximum. By Lemma~\ref{lem:trap} we know that $0< h < h_S = 2/r$. Then 
$$f = h r^2 - 2r < (2/r) r^2 - 2r = 0 .$$  So $h>0$ and $f<0$ implying that $hf < 0$.

We compute $fh$ for a nodoid at a vertical tangency, where
$\dot {y} = +\infty$ and $f = hy^2$. Then $fh = h^2 y^2 >0$. 
\enddemo

The next proposition deduces some useful properties of Delaunay surfaces.

\proclaim{Proposition} \label{prop.Delaunay}
Let $D$ be a Delaunay surface with generating curve $\delta${\rm ,}
mean curvature $h${\rm ,} maximum
$y$\/{\rm -}\/value $y_M$ and minimum $y$\/{\rm -}\/value $y_m${\rm .} Then
\begin{itemize}
\ritem{1.} If $\delta$ is a nodary or undulary{\rm ,}
the arc\/{\rm -}\/length of one period of $\delta$ is $2
\pi / h${\rm . }

\ritem{2.} If $\delta$ is an undulary{\rm ,} the period $P$ of $\delta$ satisfies
$$ 2 (y_M + y_m) \le P \le \pi (y_M + y_m).$$

\ritem{3.} If $\delta$ is a nodary{\rm ,}
$$
h = 2/(y_M - y_m).
$$
If $\delta$ is an undulary{\rm ,}
$$
h = 2/(y_M + y_m).
$$

\ritem{4.} If $\delta$ is a nodary{\rm ,} $\delta$ has nonzero curvature{\rm .}
\end{itemize}

\endproclaim

\demo{Proof} 
\begin{itemize}
\item[1.] In \cite{K} the equation for the generating curve of a Delaunay surface is
given in terms of arclength, and it is shown that nodaries and undularies are
periodic, with the arc length of one period given by $ 2 \pi / h$.

\item[2.] The period of an undulary is equal to the perimeter of the ellipse
that is rolled along the $x$-axis to generate it \cite{Eells}.
The inequalities follow by comparing twice the length of the major axis,
$2(y_M + y_m)$, with the perimeter of an enclosing circle
whose diameter is the length of the major axis, $\pi (y_M + y_m).$

\item[3.] If $\delta$ is a nodary, we calculate the force at $y_M$ and $y_m$. At
$y_M$ we get
$$ f = h y_M^2 - 2 y_M .$$
At $y_m$ we get
$$ f = h y_m^2 + 2 y_m ,$$
where the sign of the second term changes since the orientation has reversed.
Setting these equal we get
$$
h = \frac{ 2}{y_M - y_m}.
$$
For an unduloid, the proof is identical except for one sign change.

\item[4.] At a point where $\delta$ has zero curvature,
$ k_m = 0 $ so
$$
h = k_p = \frac {1 }{y\sqrt{1 + { \dot {y}}^2}}.
$$
Equation~\ref{eq:hs} gives the mean curvature of a sphere with slope $ \dot {y}$ at
height $y$ as:
$$
h_S = \frac {2}{y \sqrt{1 + { \dot {y}}^2}} > 0.
$$
Where $\dot {y}$ is defined, this is exactly twice the value of $h$ calculated above.
For a nodary, $h>h_S$ or $h<0$ by Lemma~\ref{lem:trap}, so we cannot have
$h_S = 2h$ and there are no points of zero curvature on a nodary. At points
where the nodary is vertical, $h=k_m$ since $k_p = 0$.
Thus the curvature is nonzero at
vertical tangencies except possibly when $h=0$,
which would imply that the Delaunay surface is a
plane or a catenoid, rather than a nodoid as assumed.
\end{itemize}

\enddemo

Important stability formulas for constant mean curvature surfaces were
developed by J. L. Barbosa and M. do Carmo \cite{B-D}.
Stability refers to the behavior of a
surface when a compact subsurface is deformed by a variation,
while holding its boundary fixed.
A subsurface of a complete constant mean curvature surface
is called {\it stable} if there is no 
compactly supported normal variation
which decreases its area while preserving the volumes on each side,
and {\it unstable} otherwise.
Unstable constant mean curvature surfaces cannot form part of a minimizing bubble.

\proclaim{Proposition} \label{lem.stable}
The smooth subsurfaces of a minimizing bubble
$S (v_1 ,v_2 )$ are stable subsurfaces of Delaunay surfaces{\rm .}
\endproclaim

\demo{Proof} $S (v_1 ,v_2 )$ is a piecewise-smooth surface of revolution whose smooth
subsurfaces have constant mean curvature by Theorems~\ref{thm.revolution} and
\ref{thm.bubble}. Any constant mean curvature surface of revolution in
$\RR^3$ is a Delaunay surface.
If one of these subsurfaces is unstable, then there is a
variation which maintains the volume of each region while decreasing the area.
The surface then cannot be part of a minimizing bubble.
\enddemo

We next deduce a formula for the horizontal distance between two points on a Delaunay curve
whose $y$-coordinates are known.
We require that the second point be above the first, but it can be either to the left
or to the right.
We also get a formula for the volume under a Delaunay surface.

\proclaim{Proposition} \label{dx}
Let $(x_1,y_1)$ and $(x_2,y_2)$ be two points on a Delaunay curve{\rm .}
Suppose that $y_1 < y_2$ and that between the points{\rm ,} $\displaystyle \frac{dy}{dx} \neq 0${\rm .} Then
\begin{equation}
x_2 - x_1 = \int_{y_1}^{y_2} \frac{t} {\sqrt{ (2y + t)(2y - t)} } \ dy \label{eq:dx}
\end{equation}
where $t = t(y) = h y^2 - f${\rm .}

If $x_1 < x_2$ and the curve also satisfies $\displaystyle \frac{dx}{dy} > 0$
between the two points then the volume underneath the
surface of revolution generated by the curve between $x=x_1$ and $x=x_2$ is given by\/{\rm :}\/
\begin{equation}
V(x_1 , x_2) = \int_{y_1}^{y_2} \frac{\pi y^2 t} {\sqrt{ (2y + t)(2y - t)} } \ dy .
\label{eq:dv}
\end{equation}

If instead we have  $x_2 \le x_1$  in the previous case then the formula gives the negative
of the volume underneath the
surface of revolution generated by the curve between $x=x_1$ and $x=x_2${\rm .}
\endproclaim

\demo{Proof} We solve equation~\ref{eq:force} for
$dx/dy$, using that $ dy/dx \ne 0$, and then apply the change of variables formula.
\begin{eqnarray*}
f & = & h y^2 - \frac {2y} { \sqrt{ 1+{\dot{y}}^2} } \\
\sqrt{1+{\dot {y}}^2} & = & \frac {2y} { t } \\
\frac{dx}{dy} & = &\frac{t} {\sqrt{ (2y + t)(2y - t)} } \\
\end{eqnarray*}

Note that the sign of $t=t(y) = 2y \cos \alpha $
is positive if and only if the curve is oriented to the right, i.e.\ if 
$\cos \alpha > 0 $ in equation~\ref{cos alpha}. Hence the integrand gives the correctly signed measure of
displacement in the $x$ direction.

The volume of the region enclosed by a surface of revolution generated by a graph $y(x)$
between $x_1,x_2$ with $x_1 \le x_2$ is given by
$$v = \int_{x_1}^{x_2} \pi y^2 \ dx \ .$$
Since we have already deduced a formula for $dx$, the volume formulas follow immediately. \enddemo

Note that if $y_1 > y_2$ then the formulas in Proposition~\ref{dx} would
have the  sign reversed.
A choice in sign was made in taking the positive square root, and this
choice is correct if the curve is oriented so that $y$ is increasing.
However we need to assume that $x_2 > x_1$ only when calculating the volumes,
and not in equation~\ref{eq:dx}.

The integrals of Proposition~\ref{dx} are singular at a
local minimum or maximum, where $dy/dx = 0$, and we
need to apply a change of variables to obtain a formula which
holds near such points.
Note that the formulas below do not require that $y_2 > y_1$,
but do require that
the Delaunay curve is a graph over the $x$-axis, with $x_2 > x_1$.

\proclaim{Proposition} \label{prop:dxcrit}
Let $(x_1,y_1)$ and $(x_2,y_2)$ be two points on a Delaunay curve which is a graph over the 
$x$\/{\rm -}\/axis with $x_2 > x_1$ and suppose there is exactly one critical point between the two points{\rm ,}
with $y$\/{\rm -}\/value labeled $y_{{\rm min}}$ or $y_{{\rm max}}$ depending on 
whether it is a local minimum or maximum{\rm .}
Let $h,f$ be the mean
curvature and force associated to the Delaunay surface generated by the curve{\rm .}
Set $t = t(y(z)) = h (y(z))^2 - f ${\rm .}

If $y_{{\rm min}}$ is a local minimum{\rm ,} let $ y = y(z) = y_{{\rm min}} + z^2${\rm ,} 
$\displaystyle z_1 =-\sqrt{y_1-y_{{\rm min}}}$ and $ z_3 =\sqrt{y_2-y_{{\rm min}}}${\rm .}
Then
\begin{equation}
x_2 - x_1 = \int_{z_1}^{z_3} \frac{2t } {\sqrt{ (2y + t)(2 -
hy_{{\rm min}} - hy)} } \ dz
\label{Dxmin}
\end{equation}
and
\begin{equation}
V(x_1 , x_2) = \int_{z_1}^{z_3} \frac{2\pi y^2 t } {\sqrt{ (2y +
t)(2 - hy_{{\rm min}} - hy)} } \ dz \ .
\end{equation}

If $y_{{\rm max}}$ is a local maximum{\rm ,} let $y(z) = y_{{\rm max}} - z^2${\rm ,}
$\displaystyle z_2 =\sqrt{y_{{\rm max}}-y_2}$ and $z_4 =-\sqrt{y_{{\rm max}}-y_1}${\rm .} Then
\begin{equation}
x_2 - x_1 = \int_{z_4}^{z_2} \frac{2t } {\sqrt{ (2y + t)(hy - 2 + hy_{{\rm max}})} } \ dz
\end{equation}
and
\begin{equation}
V(x_1 , x_2) = \int_{z_4}^{z_2} \frac{2\pi y^2 t } {\sqrt{ (2y +
t)( hy - 2 + hy_{{\rm max}})} } \ dz \  .
\end{equation}
Finally{\rm ,} if the curve is an increasing graph over the
$x$-axis between $x_1$ and $x_2$
then equation~{\rm \ref{Dxmin}} continues to hold if the sign of $z_1$ is changed{\rm . }
\endproclaim

\demo{Proof} These formulas follow from applying a change of variables to $y$ in
equations~\ref{eq:dx} and~\ref{eq:dv}. Suppose $y_{{\rm min}}$ is a local minimum. 
Divide the Delaunay curve into the two pieces on
each side of the minimum, and apply Proposition~\ref{dx} and
the substitution $y = y_{{\rm min}} + z^2$ to each piece. 

One of the terms in the square root in the denominator in equation~\ref{eq:dx} factors as
$$2y - t = 2y -hy^2 +f = (y-y_{{\rm min}})(2-hy_{{\rm min}} - hy) = z^2 (2-hy_{{\rm min}} - hy)$$
since $f =hy_{{\rm min}}^2 - 2y_{{\rm min}}$.
A resulting factor of $z$ in the denominator, which causes the singularity,
cancels with a $z$  in the numerator from the change of variables.
The integrands are even functions of $z$, and $y=y_{{\rm min}}$
corresponds to $z=0$, so the curve portions on opposite sides of the minimum can be mapped to
positive and negative $z$ values, and combined into one integral. 

The maximum case is similar. Near a maximum the substitution $y = y_{{\rm max}} - z^2$
is used. Several minus signs appear in the derivation, resulting in the given expression.

The final statement in the proposition is a direct application of the change of variables formula
to equation~\ref{eq:dx}.
\enddemo
 
\section{Torus bubbles} \label{sec.torus} 

We established in Lemma~\ref{lem.stable} that a minimizing bubble must be obtained
by revolving a union of Delaunay curves contained in the upper half plane.
A key case of such a bubble is the {\it torus bubble}, constructed as follows.
Take two circular arcs of the same
radius, facing each other, each with one endpoint and center on the $x$-axis,
and connect the other endpoints with two different
Delaunay curves meeting at 120 degrees. Rotating around the $x$-axis,
we get a piecewise-smooth surface surrounding two components, one homeomorphic to a torus,
which we call the {\it torus component} $T$, and one homeomorphic to a ball, which we call
the {\it ball component} $B$. It is not immediately clear whether it is possible to make such a
construction so that the curves meet at $120^\circ$ angles and the mean curvatures sum to zero
around each triple curve. Such torus bubbles do indeed exist, and we will need
to show that none of them are minimizers.
Figure~7, due to J. Sullivan, shows a torus bubble.

The reason that torus bubbles play so central a role in our argument is due to another key
result of Hutchings:

\proclaim{Theorem} \label{two.types} A minimizing bubble enclosing two equal volumes must
be either a double bubble or a torus bubble{\rm .}
\endproclaim

\figin{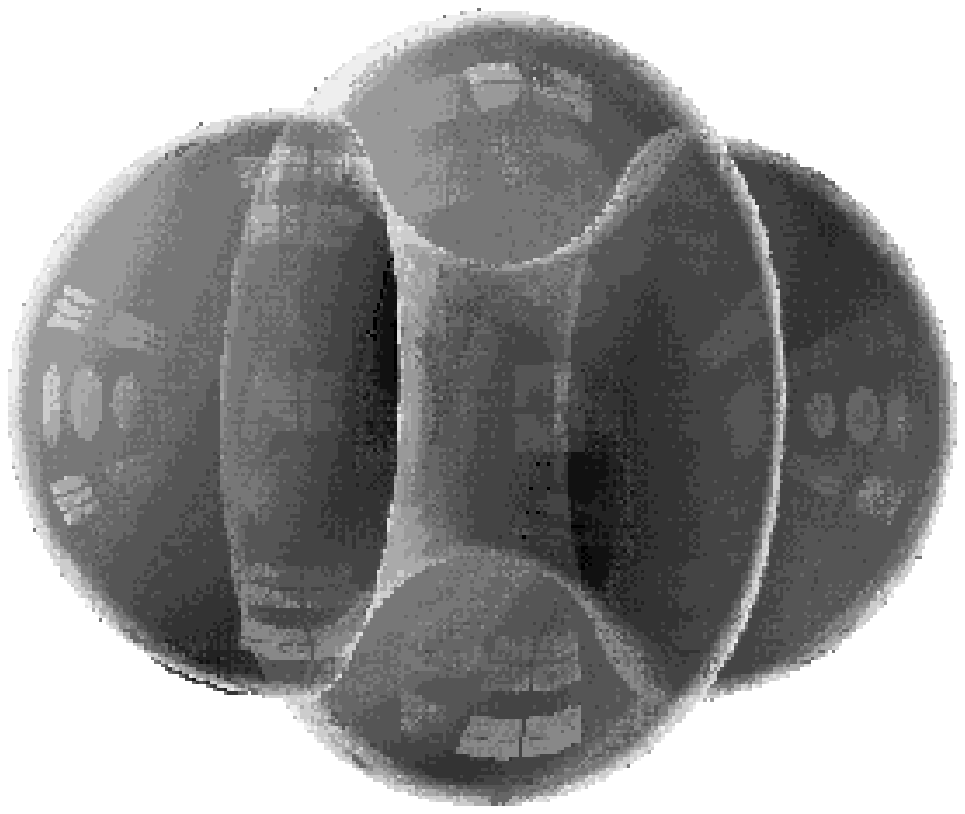}{400}
\centerline{{Figure 7}. A torus bubble in ${\Bbb R}^3$.}
\medbreak

\demo{Proof} Theorem~\ref{connected} states that $S(v,v)$ contains two connected
components. The two possible configurations are then a consequence of Hutchings' structure
theorem for minimizing bubbles \cite[Th.~5.1]{Hu}. \enddemo

We label the surfaces of a minimizing torus bubble as indicated in
Figure~8.
The component homeomorphic to a solid torus is denoted $T$
and the component homeomorphic to a ball is denoted $B$.
The volumes of the  inner ball component and outer torus component
are denoted by $v_i$ and $v_o$ respectively.
The inner Delaunay surface on the boundary of
$T$ is denoted $T_i$ and the outer Delaunay surface is denoted $T_o$.
The generating curve
of $T_o$ is denoted by $\tau_o$ and that of $T_i$ by $\tau_i$.
The angles subtended by
the generating curves for $S_1$ and $S_2$ are denoted by $\theta_1$ and
$\theta_2$. The mean curvatures of $T_i$ and $T_o$ are $h_i$ and $h_o$,
with signs chosen so
that so that $h_i$ is positive if the mean
curvature vector points into $B$ and negative if it
points into $T$, and $h_o$ is positive if the mean curvature vector points into
$T$ and negative if it points towards outside the bubble.

Denote by $(x_1,y_1)$ the initial point where the two Delaunay surfaces start, and
$(x_2,y_2)$ the point where they rejoin. 
To specify a sign for the mean curvature,
we orient the generating curves $\tau_i$ and $\tau_o$ of $T_i$ and $T_o$
so that they run from $(x_1,y_1)$ to $(x_2,y_2)$.

\figin{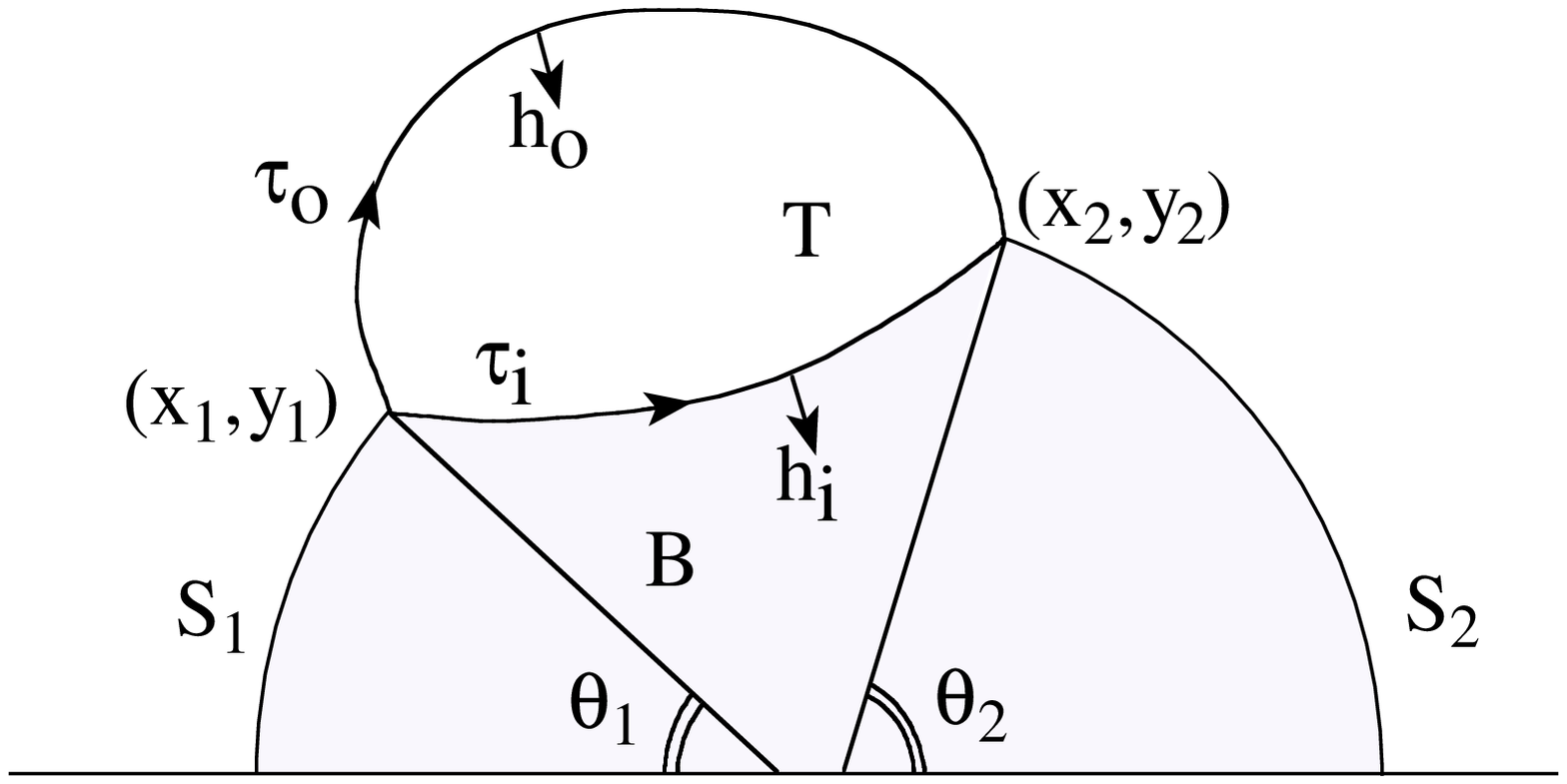}{350}
\centerline{{Figure 8.} Parameters of a torus bubble.}
\medbreak
Torus bubbles for various ranges of $\theta_1$ are indicated in
Figures~18,~19 and~20.
We will show that $\theta_1$ and $h_o$ parameterize the space of torus bubbles and that no
choice of $\theta_1$ and $h_o$ gives a minimizing bubble. To do so, we first present some
properties of minimizing torus bubbles, mainly derived from the work of Hutchings.

\proclaim{Lemma} \label{radius} In a minimizing torus bubble{\rm ,} $S_1$ and $S_2$ have the same mean
curvature{\rm .}
\endproclaim

\proclaim{Lemma} \label{param.bounds} In a minimizing torus bubble{\rm ,} $h_o \ge
0$ and $h_i = 2- h_o${\rm .}
\endproclaim

{\it Proof}. The three mean curvatures around a triple curve sum to zero by
Lemma~\ref{curv.sum} implying that $h_o = 2- h_i$. If $h_o < 0$ then there is a deformation
which pushes $T_o$ outward, increasing the volume of the torus component and decreasing its
area. Corollary~\ref{area.increases} then gives a contradiction. 
\hfill\qed

\proclaim{Proposition} \label{extra.reflect} The regions bounded by a minimizing torus bubble are not
all preserved by reflection through some plane $P$ perpendicular to the $x$\/{\rm -}\/axis{\rm .}
\endproclaim

{\it Proof}. If there is such a symmetry, then the bubble is invariant by reflection through the
$xy$-plane and by reflection through a perpendicular plane~$P$. The argument given in
Theorem~\ref{thm.revolution} shows that the bubble is a surface of revolution around the line of
intersection of these two planes. But it is also a surface of revolution around the
$x$-axis, and therefore a surface of revolution around two perpendicular axis. It must then be a
union of spheres, and not a torus bubble. \hfill\qed

\proclaim{{C}orollary} \label{theta_1<theta_2}
Given a minimizing equal\/{\rm -}\/volume torus bubble{\rm ,} we can rescale and reflect so
that $\theta_1 < \theta_2$ and  each of the two circular arcs generating $S_1$ and $S_2$
has curvature one{\rm .}
\endproclaim

{\it Proof}.
If $\theta_1 = \theta_2$, Lemma~\ref{extra.reflect} gives a contradiction.
If $\theta_1 > \theta_2$, a left/right reflection interchanges the two angles.
Lemma~\ref{radius} implies that the two circular arcs generating $S_1$ and $S_2$ have the same
curvature, and thus the same radius. We can rescale and assume without loss of generality that each
has radius one (and thus mean curvature two).
\hfill\qed \bigbreak

From now on we will restrict attention to torus bubbles with $\theta_1 < \theta_2$ and with
spherical caps of radius one. 

\proclaim{Proposition} \label{unstable.nodoid}
If a torus bubble contains a surface $T_1$ which is a
subsurface of a nodoid and the
nodary generating $T_1$ contains two interior vertical tangencies{\rm ,}
then $T_1$ is unstable and not a part of a minimizing bubble{\rm .}
\endproclaim

{\it Proof}. We will construct a volume preserving Jacobi field supported on a proper subset of
$T_1$, establishing instability as in \cite[Prop.~2.24]{B-D}. 

We start with the constant unit vector field $v$  
pointing in the direction of the positive $y$-axis.
The vector field $v$ is tangent to $T_1$ along two circles which are
generated by rotating the interior vertical tangencies of the nodary.
Let $T_1'$ denote the compact portion of $T_1$ between these two circles.
Define the vector field $W = \langle v, N\rangle N $. where $N$ is the normal vector field along $T_1'$.
and extend $W$ to be the zero vector field elsewhere on $T_1$. 
The integral of the function $f = \langle v,N\rangle $ along $T_1$ is zero, since the involution
$(x,y) \to (-x,-y)$ of $T_1$ takes $f$ to $-f$.
So the normal variation defined by $f$ preserves volume in the sense
of \cite[Prop.~2.24]{B-D}.
Moreover \cite[Prop.~2.24]{B-D} shows that $fN$ is a Jacobi field, since
it is constructed by taking the normal part of a constant vector field, and
since the nodary extends beyond the two vertical tangencies, the
corresponding Jacobi field is supported on a proper subset of $T_1$.
The Morse Index Theorem implies that the surface $T_1$ is unstable
and that there is a smaller area surface separating the same volumes.
\hfill\qed \bigbreak

The next proposition establishes that minimizing torus bubbles do not contain Delaunay curves
which are longer than a full period.

\proclaim{Proposition} \label{period.size}
If $\tau_i$ or $\tau_o$ contains a full period{\rm ,} then the torus bubble is not a
minimizer{\rm . }
\endproclaim

\demo{Proof} Suppose to the contrary that a minimizing torus bubble
contains a subsurface generated by a curve which contains a full period.
If a boundary surface of the torus component is a subsurface of a nodoid with
width more than a period, then the associated nodary contains two vertical tangencies,
which implies instability by Proposition~\ref{unstable.nodoid}.
Similarly, if a boundary surface of the torus component is a subsurface of a
cylinder of radius $r$ and its width is greater than $2\pi r$, the period of
a cylinder, than the surface is shown to be unstable in \cite{B-D}.

If a boundary surface of the torus component is a subsurface $U$ of an unduloid then the lemma
is less obvious. Note that we cannot necessarily construct a volume preserving Jacobi field from
the vector field $\partial / \partial x$, since the Delaunay curve may not contain two maxima
or two minima. We establish the lemma with a reflection argument.  Assume first that the
unduloid containing a full period is the inner surface $T_i$.

An arc $u$ of an undulary which is longer than a period contains three distinct points
with the same $y$ value, $(x_1,y_1),(x_2,y_1),(x_3,y_1)$, with $x_1 < x_2 < x_3$, and with
the slopes at $(x_1,y_1), (x_3,y_1)$ each being nonzero and
equal to the negative of the slope at $(x_2,y_1)$.
Reflect the part of $u$ between the planes $x=x_1$ and $x=x_3$ through the plane $x= (x_1 +
x_3) / 2$ to get a new surface of revolution
$U'$ generated by $u'$, the reflection of $u$.
Reflection preserves both the area and the volume
under $U$ between the planes $x=x_1$ and $x=x_3$. 
Suppose first that $u'$ does not intersect $\tau_o$. Then
the volume of the torus bubble is also preserved. The reflection creates a surface which
is not smooth along the planes $x=x_1$ and $x=x_3$ since reflection changes the sign
of the slope of the reflected surface along these planes, and this slope was not zero.
Thus we have found a solution to the minimization problem which contains
a surface which contradicts Theorem~\ref{thm.bubble}.

Possibly $u'$ does intersect the outer curve $\tau_o$. See Figure~9.
In that case the total volume enclosed by
the union of $U'$ and $T_o$ between $x=x_1$ and $x=x_3$ is increased, or at least
not decreased. We know that $U'$ does not intersect
the spherical caps of the torus bubble, since they do not meet the region between the planes
$x=x_1$ and $x=x_3$. Since the volume under $U'$ is equal to the volume of the ball region
of the torus bubble $B$, the volume in
the remaining (disconnected)
region must be larger than that in the original torus region $T$. There is no
increase in area, contradicting Corollary~\ref{area.increases}.

\figin{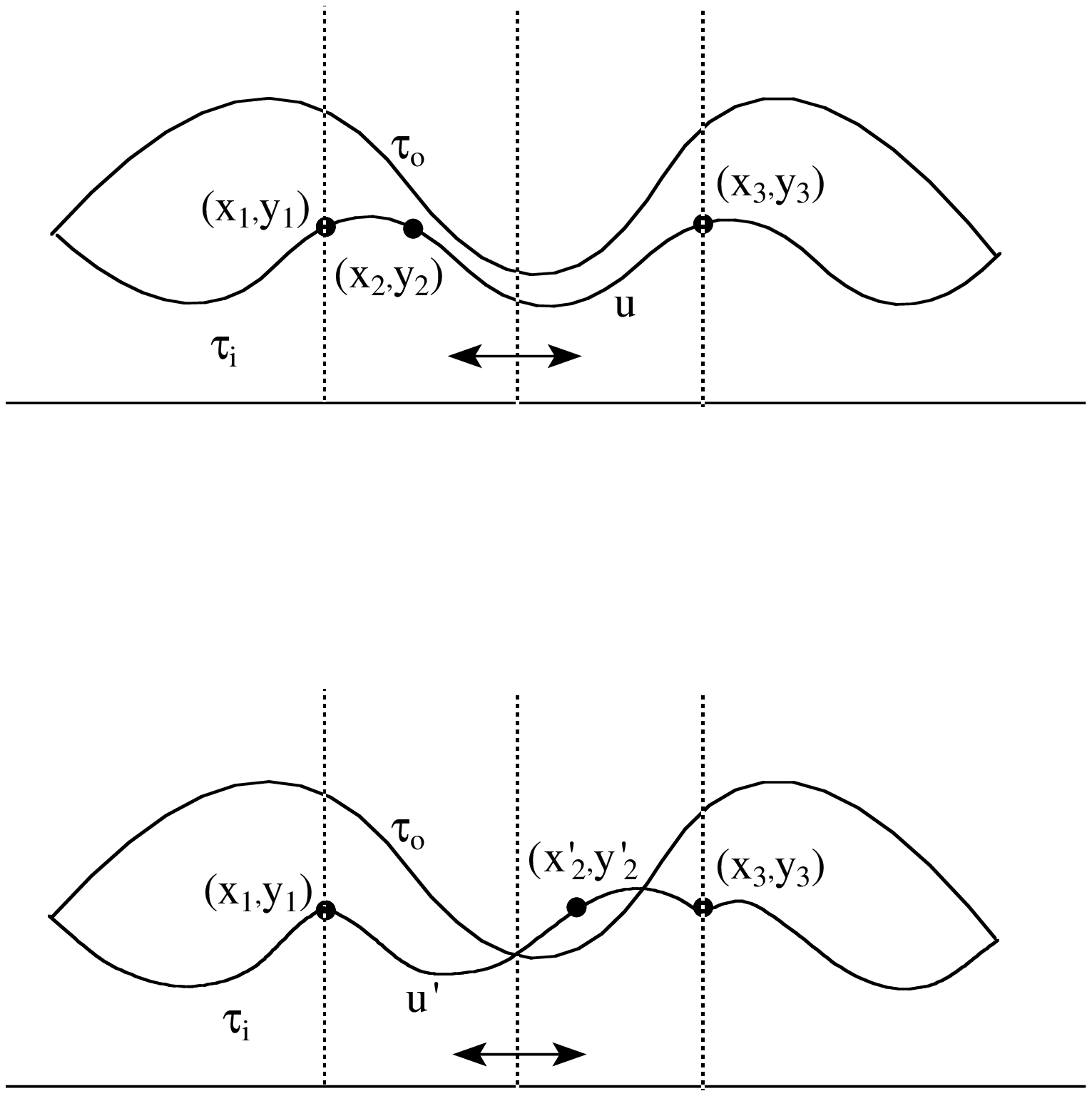}{350}
\centerline{{Figure 9}. Reflection of a piece of $T_i$ could lead to extra intersections.}
\medbreak

Note that have used here the fact that two distinct 
Delaunay curves intersect locally at finitely many points.

Assume now that the unduloid is the outer surface $T_o$. We again reflect it to get a new
boundary surface $U'$. The total volume underneath $U'$ is the same as that underneath $U$. 
Possibly $T_i$ intersects the reflected surface $U'$. In this case the
total volume underneath the union of the two surfaces has increased. Divide the region
underneath $U'$ so that the region under $T_i$ is allocated entirely to $B$, as
before, and all other regions are allocated to $T$. Then the total volume of $B$ is
preserved, that of $T$ is increased, and the area is  unchanged, a contradiction as before. \enddemo

We note that a somewhat similar argument was used by Athanassenas
\cite{Athanassenas} in studying the stability of minimizing capillary surfaces. See also
\cite{Vogel}. Our situation has some additional complications due to the possible intersections
resulting from the reflection. 

\proclaim{Lemma} \label{lem.60} 
If a torus bubble has $\theta_1 < 60^\circ$ and $\theta_2 < 60^\circ$ then it is not minimizing. In
particular{\rm ,} $\theta_2 \ge 60^\circ$ in a minimizing torus bubble{\rm .}
\endproclaim

\demo{Proof} If both $\theta_1$ and $\theta_2$ are less than $60^\circ$, the outer bubble $T_o$ is
generated by a nodary that becomes vertical twice. Proposition~\ref{unstable.nodoid} implies it is unstable. 
Since $\theta_2$ is the larger angle, $\theta_2 \ge 60^\circ$.
\enddemo

\proclaim{Lemma} \label{120}
Suppose that $T$ is a torus bubble with $\theta_2 \ge 120^\circ${\rm .}
The the torus component is completely contained inside the smaller region of
the double bubble whose larger outer sphere subtends an angle equal to $\theta_2${\rm .}
Similarly the ball component contains in its interior the larger region of this double bubble{\rm .}
\endproclaim

\demo{Proof} Consider a double bubble which is symmetric around
the $x$-axis and whose larger spherical cap $S_2$
is on the right and subtends an angle of $\theta_2$ with the  $x$-axis.
Denote the smaller (leftmost) spherical cap by $S_o$ and the spherical cap forming
the interface by $S_i$.
Then $S_o$ is tangent to $\tau_o$ at $(x_1,y_1)$.
We will show that a
torus region must be contained inside the smaller component
$B_1$ of this double bubble, the component bounded by  $S_2$ and $S_i$. 
Otherwise, we will show that  $\tau_i$ and $\tau_o$ cannot
rejoin after leaving $(x_1,y_1)$, and thus cannot generate a torus region.
See Figures~10 and  8.

Given a vector in the upper half-plane, there is a unique circle
perpendicular to the $x$-axis which is tangent to that vector.
A Delaunay curve tangent to the same vector and generating a Delaunay surface
with greater mean curvature than the corresponding sphere
is a nodary, by Lemma~\ref{lem:trap}.
A Delaunay curve tangent to the same vector and generating
a Delaunay surface with smaller,
but still positive, mean curvature is an undulary or horizontal line.

If the mean curvature of $T_o$ is equal to that of $S_o$, then
the mean curvature of $T_i$ is equal to that of $S_i$ by Lemma~\ref{param.bounds}.
It follows from Lemma~\ref{lem:trap} that
$T_o, T_i$ and  $S_o, S_i$ coincide, and that 
there is no torus component.

Suppose now that the mean curvature of $T_o$ is smaller than that of $S_o$.
Then the mean curvature of $T_i$ is smaller than that of
$S_i$ by Lemma~\ref{param.bounds}.
Since $h_o \ge 0 $, also by Lemma~\ref{param.bounds},
it follows that $T_o$ cannot be a nodoid, and therefore
$\tau_o$ is a graph over the $x$-axis. In particular, $\tau_o$
does not recross the spherical cap $S_o$.
The generating curve $\tau_i$ for $T_i$ leaves the point $(x_1,y_1)$
on the positive (right) side of $S_i$, and lies above $S_i$.
Since $S_2$ forms part of the torus bubble, $\tau_i$
must remeet $\tau_o$ without crossing $S_2$ or itself if
it is to generate a torus component with  $\tau_o$.
If it is a graph over the $x$-axis then it cannot remeet $\tau_o$ at all.
If it is not a graph,
then it is a nodary which generates a nodoid with negative mean curvature.
Such a curve crosses $S_2$ before it can remeet $\tau_o$.

Finally, suppose that
the mean curvature of $T_o$ is greater than that of $S_o$.
Then Lemma~\ref{param.bounds} implies that
the mean curvature of $T_i$ is greater than that of $S_i$,
and Lemma~\ref{lem:trap} implies that both $T_o$ and $T_i$ are subsurfaces of nodoids 
and that near $(x_1,y_1)$, $T_o$ lies below $S_o$ and $T_i$ lies below $S_i$.
See Figure~11.
The nodary $\tau_o$ generating $T_o$ either crosses itself before leaving
$B_1$ or leaves $B_1$ through $S_i$.
Similarly, the nodary $\tau_i$ generating $T_i$
either crosses itself before leaving
$B_1$ or leaves $B_1$ through $S_o$. 
The curves  $\tau_o$ and  $\tau_i$
cannot self-intersect before remeeting, since they generate a torus component.
Therefore either $\tau_o$ and $\tau_i$ are contained in $B_1$ or they remeet
for the first time after 
$\tau_o$ leaves $B_1$ through $S_i$ and $\tau_i$ leaves $B_1$ through $S_o$. 
The latter case is impossible, since $\tau_i$ must cross $\tau_o$
before it can leave $B_1$ through $S_i$.
It follows that the torus component is contained inside $B_1$.

The spherical caps $S_i$ in the double bubble and $S_1$ in the
torus bubble are disjoint, since $S_1$
has boundary circle lying to the left of $S_i$, and each is convex in opposite
directions.  So no part of the torus bubble
intersects the interior of the larger, rightmost component of the corresponding
double bubble, which must then be contained in the ball component of the torus bubble. 
This proves the last assertion of the lemma. \enddemo 

\proclaim{{C}orollary} \label{<120} 
In a minimizing equal\/{\rm -}\/volume torus bubble{\rm ,} $\theta_2 < 120^\circ${\rm .}
\endproclaim

\begin{center}
\BoxedEPSF{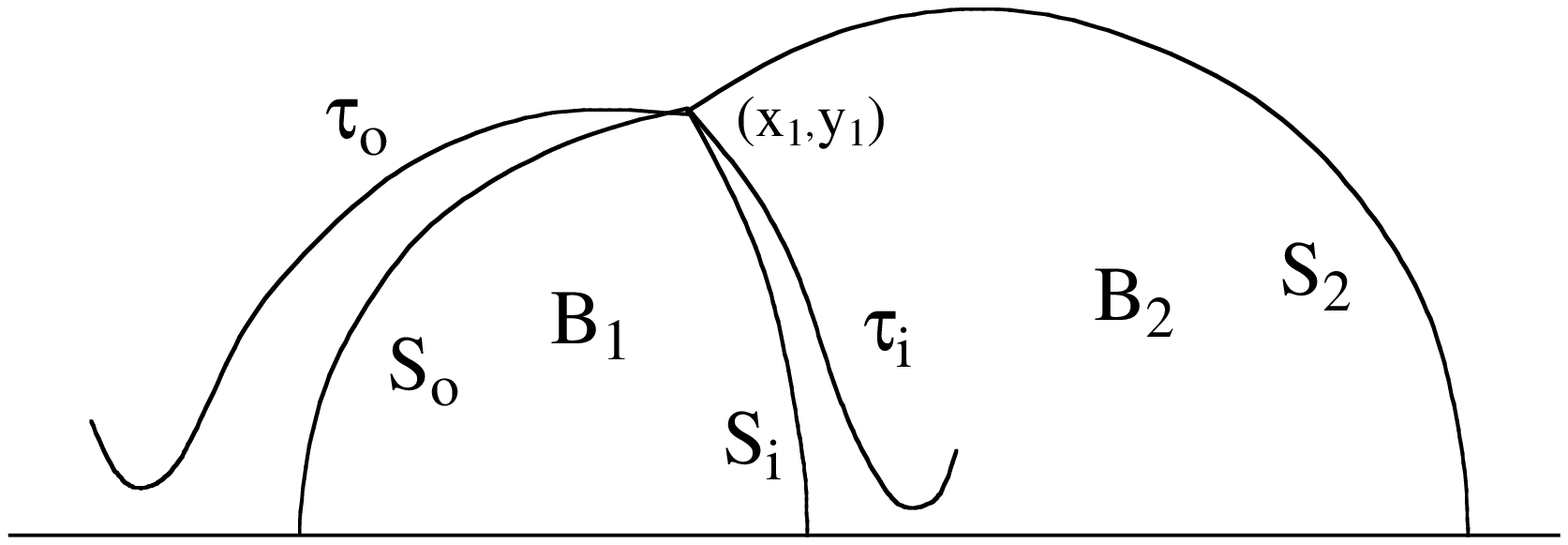 scaled 400}
\end{center}
 \centerline{{Figure 10}. If $h_o$ is too small, a torus component cannot exist.}

\figin{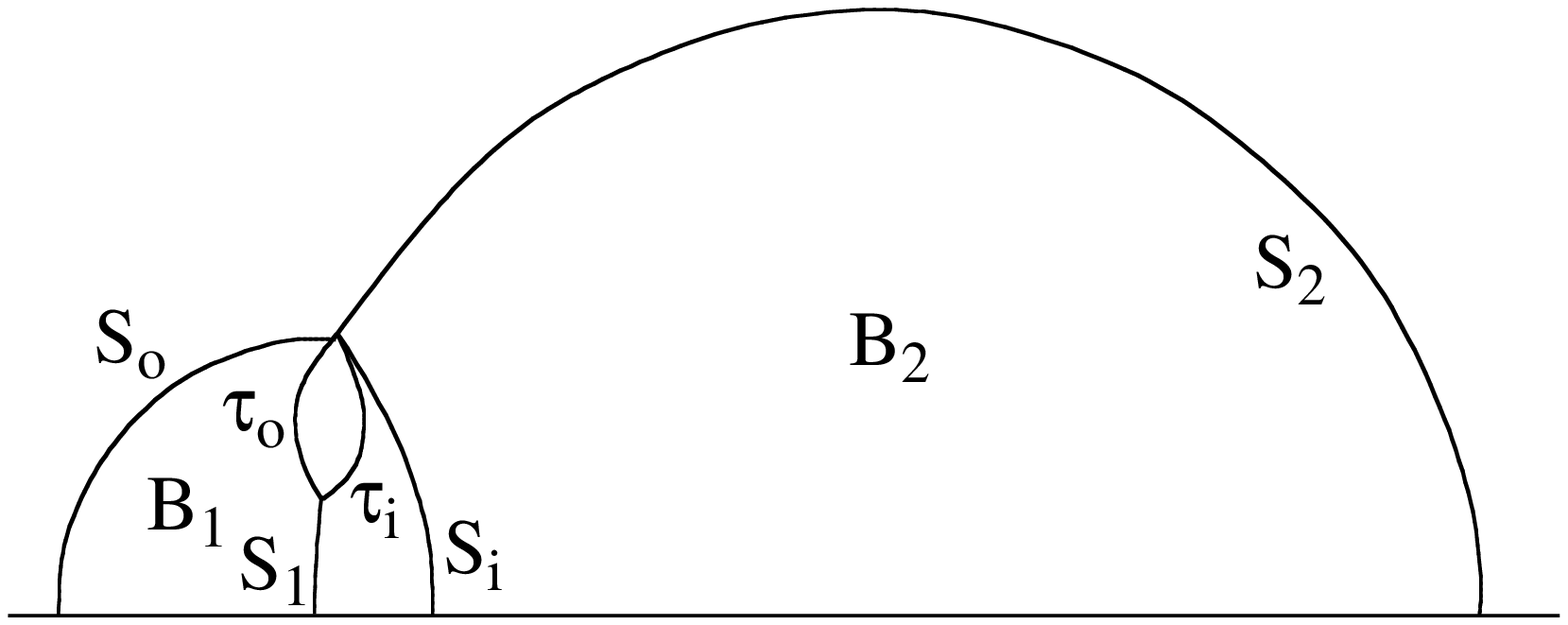}{350}%

 {{Figure 11}. The torus component is trapped inside a component of an

 associated double
bubble if $\theta_2 > 120^\circ$.}

\demo{Proof}
If $\theta_2 > 120^\circ$ then Lemma~\ref{120} implies that the volume of the
torus component is strictly less than the volume of the smaller component of the double
bubble whose spherical caps subtend angles $\theta_2$ and $(240 - \theta_2)^\circ$.
The other component of the torus bubble contains in its
interior the larger component of this double bubble, which in turn contains more than half the total volume
of the double bubble. It follows that the torus component has strictly less volume than the ball
component of the torus bubble, contradicting the assumption of equal volume.  \enddemo

Recall that $\theta_2$ is the larger of the two angles subtended by the
spherical caps, so the Corollary implies that both of $\theta_1$ and $\theta_2$ are
less than $120^\circ$.

\proclaim{Lemma} \label{h_o > 0} 
In a minimizing equal\/{\rm -}\/volume torus bubble{\rm ,} $h_o > 0${\rm .}
\endproclaim

\demo{Proof} We know from Lemma~\ref{param.bounds} that $h_o \ge 0$.  If $h_o = 0$ then
$\tau_o$ is a catenary. Since $y_2 \ge y_1$ the catenary is given
by a strictly increasing graph for $\{x \ge x_2 \}$,
implying in particular that $\theta_2 > 120^\circ$ which violates
Lemma~\ref{120}. \enddemo

Knowledge of $ h_o, \theta_1, \theta_2$ can be used to give a fairly
accurate qualitative picture of the torus bubble.

\proclaim{Proposition} \label{torus.properties}
A minimizing equal\/{\rm -}\/volume torus bubble has the following properties\/{\rm :}\/
\begin{itemize}
\ritem{1.} $\tau_i$ has negative slope at $(x_1,y_1)$ if and only if $\theta_1 > 30^\circ${\rm .}  Its
angle with the positive $x$\/{\rm -}\/axis is $30^\circ - \theta_1 ${\rm .}

\ritem{2.} $\tau_i$ has positive slope at $(x_2,y_2)$ if and only if $\theta_2 >
30^\circ${\rm .} Its angle with the positive $x$\/{\rm -}\/axis is $\theta_2 - 30^\circ ${\rm .}

\ritem{3.} $\tau_o$ has positive slope at $(x_1,y_1)$ if and only if $\theta_1 >
60^\circ${\rm .} Its angle with the positive $x$\/{\rm -}\/axis is $150^\circ - \theta_1 ${\rm .}

\ritem{4.} $\tau_o$ has angle at $(x_2,y_2)$ with the positive
$x$\/{\rm -}\/axis equal to $\theta_2 - 150^\circ${\rm .}

\ritem{5.} $\tau_i$ is a graph{\rm .} It has a unique local minimum if and only if
$\theta_1 \ge 30^\circ${\rm .}  It never has a local maximum{\rm .}  If $T_i$ is a nodoid then $h_i < 0${\rm .}

\ritem{6.} $\tau_o$ always has a unique local maximum{\rm .}  It has a vertical
tangent on the left if and only if $\theta_1 \le 60^\circ${\rm .}  It never has a vertical tangent on the
right{\rm .} It never has a local minimum{\rm .}

\ritem{7.} If $\theta_1 \ge 30^\circ$ the local minimum $y_{{\rm min}}$ of
$\tau_i$ has value 
$$
y_{{\rm min}} = \frac{- f_i }{1 + \sqrt{1 + f_i h_i} },
$$ 
where $f_i$ is the force associated to $T_i${\rm .}

\ritem{8.} The local maximum $y_{{\rm max}}$ of $\tau_o$ has value 
$$
y_{{\rm max}} = \frac {1 + \sqrt{1 + f_o h_o} }{h_o },
$$ 
where $f_o$ is the force associated to $T_o${\rm .}
\end{itemize}

\endproclaim

\demo{Proof} The first four assertions follow from the fact that the surfaces
forming a torus bubble meet at $120^\circ$ angles along a triple curve. 

If $\tau_i$ is not a graph then $\theta_2 > 120^\circ$, violating
Corollary~\ref{<120}.  If $\theta_1 \ge 30^\circ$ then $\tau_i$ has
nonpositive slope at $(x_1,y_1)$. Since  $\theta_2 >  \theta_1$ the curve must have a
local minimum in this case. 

Suppose $\tau_i$ has a local maximum. Since $\theta_2 > 60^\circ$, $\tau_i$ has
positive slope at $(x_2,y_2)$ and must have a local minimum as well. However
$\theta_2 > \theta_1$ so $\tau_i$ must be longer than a period, contradicting 
Proposition~\ref{period.size}. So it has no local maximum.  

If $h_i \ge 0$ and $T_i$ is a nodoid, then we can compare $T_i$ and $T_o$ to a double
bubble whose spherical interface $S_i$
coincides with $S_1$, and apply Lemma~\ref{lem:trap} as in
Lemma~\ref{120}.  $T_i$ is
trapped inside the component $B_2$
to the right of the interface and $T_o$ is trapped inside the
component $B_1$ to the left of the interface, as in Figure~12.
In particular, $h_o < 0$ in this case, a contradiction.
So if $T_i$ is a nodoid then $h_i < 0$.

\figin{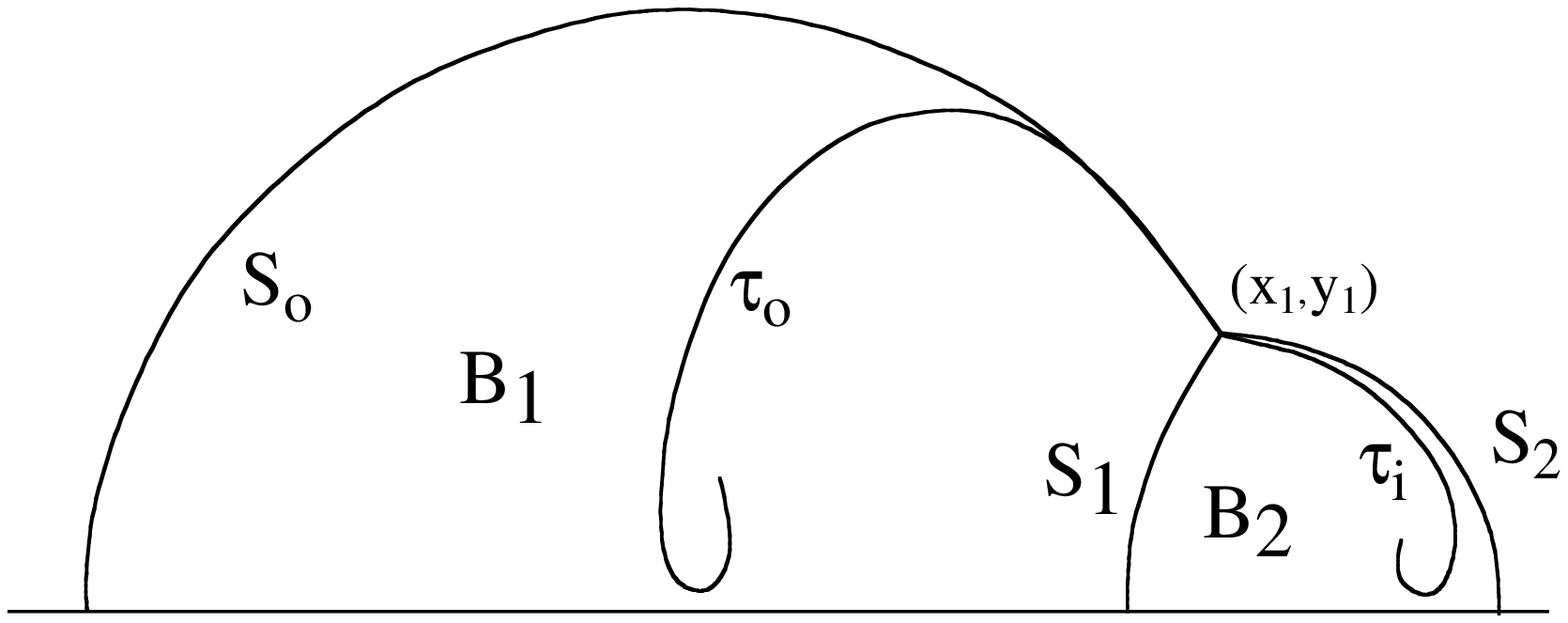}{350}
\medbreak\centerline{{Figure 12}. A double bubble traps $T_i$ and $T_o$ if $T_i$ is a nodoid with $h_i >0$.} 
\medbreak

If $\theta_1 < 30^\circ$ then $\tau_i$ has positive slope at $(x_1,y_1)$. If it
has a minimum then it must pass through a maximum first, so it has no minimum. 

$\tau_o$ is always oriented upwards at $(x_1,y_1)$ and downwards at
$(x_2,y_2)$, so it always has a local maximum. If it has a local minimum then it would
contain two local maxima, violating Proposition~\ref{period.size}. If $\theta_1 < 60^\circ$ then the
curve starts out to the left, therefore it has a vertical tangent at which it
changes direction from left to right. Proposition~\ref{unstable.nodoid} implies that there cannot be
a second vertical tangent, so it has no vertical tangent on the right, as claimed. If $\theta_1 \ge
60^\circ$ then $\tau_o$ initially is going right.  If it has a vertical tangent on the right
then the angle $\theta_2$ must be less than $60^\circ$, contradicting the assumption that $\theta_1
\le \theta_2$. This proves the assertion concerning
   vertical tangents of  $\tau_o$.   

To calculate the local minimum $y_{{\rm min}}$ of $\tau_i$, we use
equation~\ref{eq:force} of Section~\ref{sec:Delaunay} for the force of $T_i$,
$$
f_i = h_i y^2 - 2 y \cos \alpha \ .
$$
At a minimum, $\alpha = 0$ and $\cos \alpha \ = 1$, so
$$
f_i = h_i y^2 - 2y .
$$
If $h \neq 0$ this quadratic expression gives two roots for $y$:
$$
y = \frac {1 \pm \sqrt{1 + h_i f_i}}{h_i} \ .
$$
The value of the minimum is given by  
\advance\eqcount by 11
\begin{equation}
\frac {1 - \sqrt{1 + h_i f_i}}{h_i},  \label{eq:ymin1}
\end{equation}
whatever the sign of $h_i$. If $h_i f_i > 0$, as Corollary~\ref{cor:signhf} implies happens at a
local minimum for a nodoid, then the other root is negative and meaningless. If $h_i f_i < 0$, as
Corollary~\ref{cor:signhf} implies happens at a local minimum for an unduloid, then the other
root corresponds to the local maximum of the unduloid. A double root occurs only in the case of
a cylinder. equation~\ref{eq:ymin1} becomes linear when $h_i=0$.  An equivalent expression for
the minimum, which holds also when $h_i = 0$, is
\begin{equation}
y_{{\rm min}} = \frac{- f_i }{1 + \sqrt{1 + f_i h_i} }.  \label{eq:ymin}
\end{equation}
For a local maximum, a similar analysis gives that $y_{{\rm max}}$ is one of the two roots
\begin{equation}
y_{{\rm max}} = \frac {1 \pm \sqrt{1 + f_o h_o} }{h_o }.  
\end{equation}

Since by Corollary~\ref{T_o.is.nodoid} $T_o$ is always a nodoid, $h_o > 0$ and $f_o > 0$, there is
a unique positive root, and 
\medbreak
\hfill ${\displaystyle
y_{{\rm max}} = \frac {1 + \sqrt{1 + f_o h_o} }{h_o }.  }$\hfill
\enddemo
 
\proclaim{Proposition} \label{determined}  Given $\theta_1$ and $h_o${\rm ,} there is at most one
corresponding minimizing equal\/{\rm -}\/volume torus bubble{\rm .}
\endproclaim

\demo{Proof} Given $\theta_1$ and $h_o$, $h_i = 2-h_i$ is also determined. Thus $S_1,T_o$ and
$T_i$ are uniquely determined.
The spherical cap $S_2$ is uniquely determined by the second
intersection point $(x_2,y_2)$ of $\tau_o$ and $\tau_i$,
assuming that this second point exists.
If $\tau_o$ and $\tau_i$ do not meet in a second
point then there is no torus bubble corresponding to $\theta_1$ and $h_o$. In general, 
$\tau_o$ and  $\tau_i$ may intersect in many points, and it is necessary to show that only one of
these can possibly lead to a torus bubble.

If $\tau_o$ and $\tau_i$ are graphs, then it is clear that their first
intersection point with $x$ coordinate larger than $ x_1$ is the unique intersection point
which can give a torus bubble. In general,  $\tau_i$ is a graph, but not $\tau_o$.  However
Proposition~\ref{torus.properties} shows that $\tau_o$ has no vertical tangents for $ x > x_1$, so
that both curves are graphs in this region and the same argument applies.
\enddemo

\begin{eqnarray*}
\BoxedEPSF{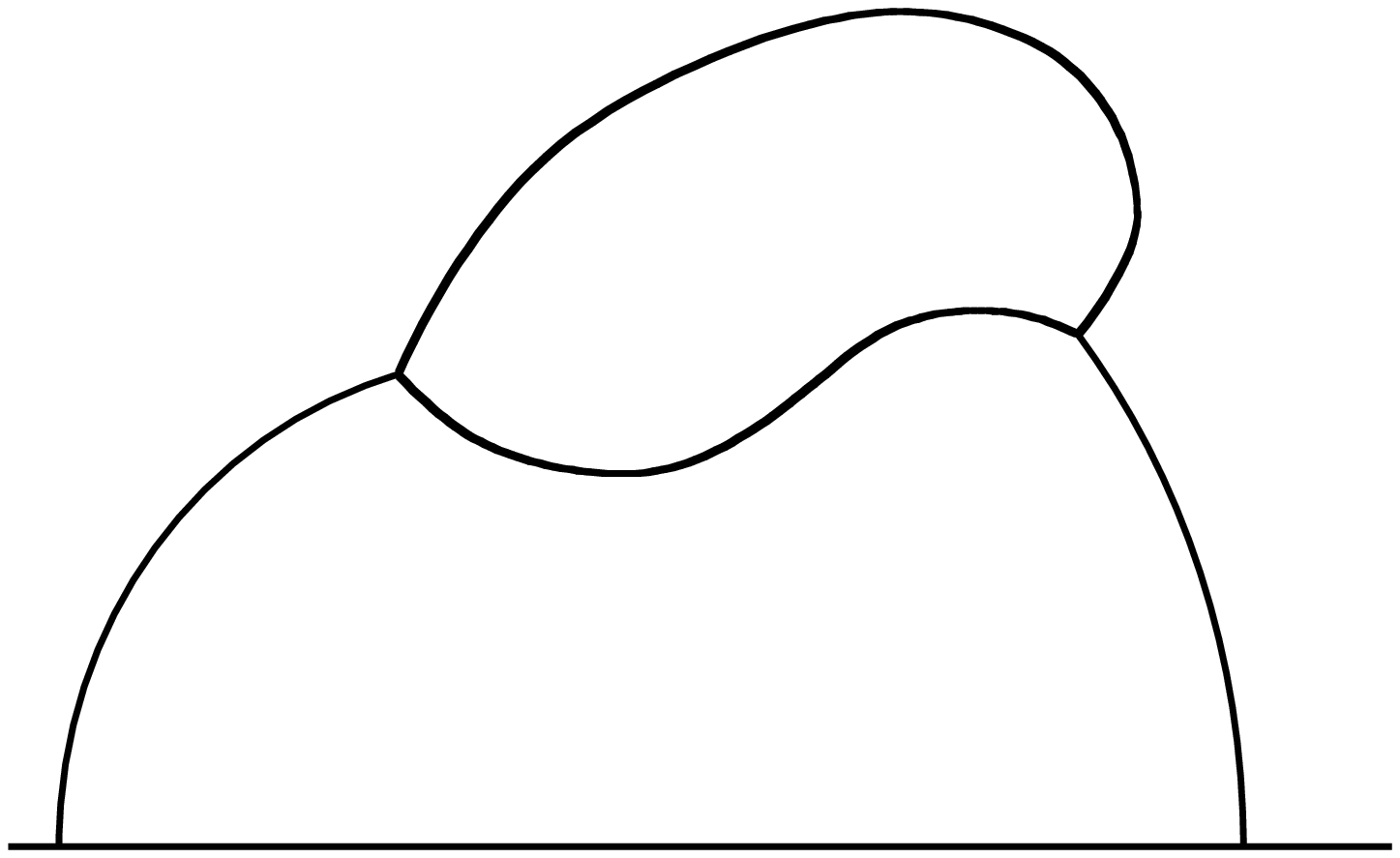 scaled 275}
&&\qquad \BoxedEPSF{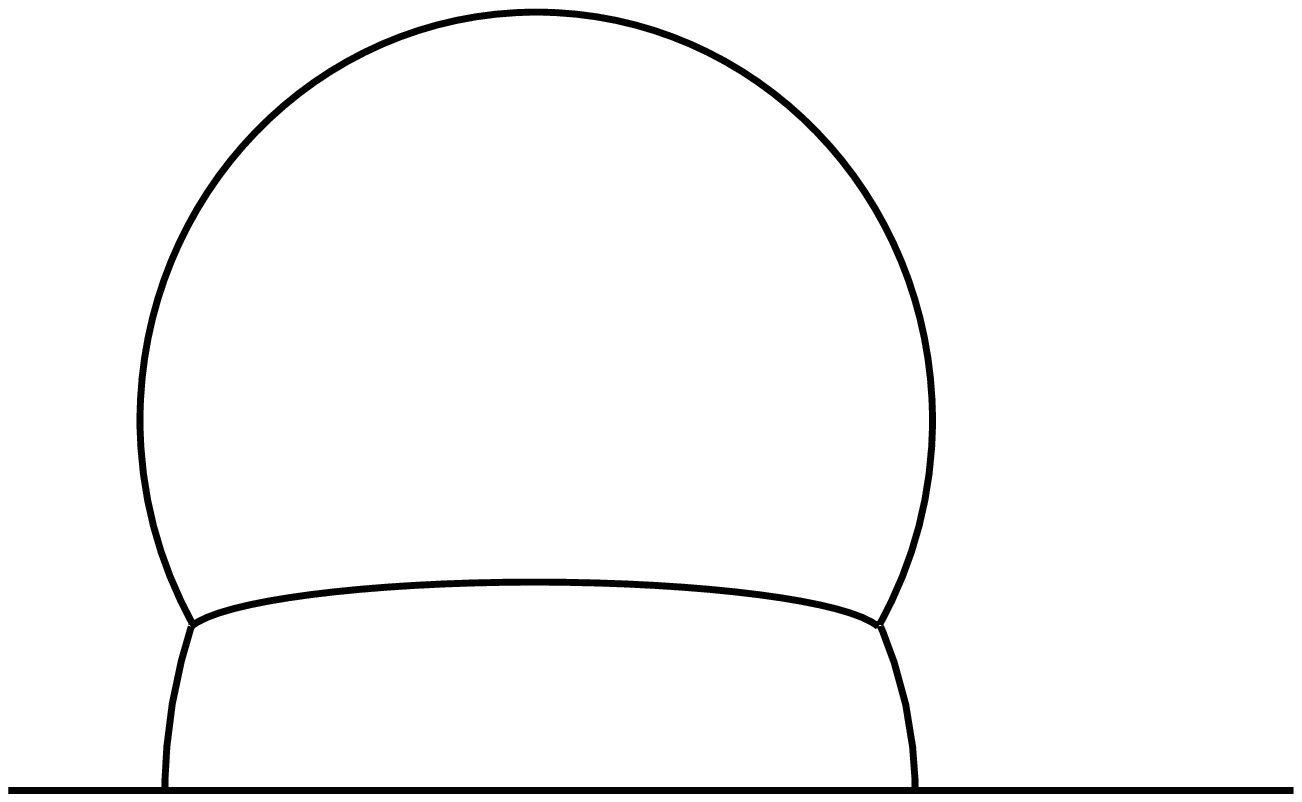 scaled 275}\\
\BoxedEPSF{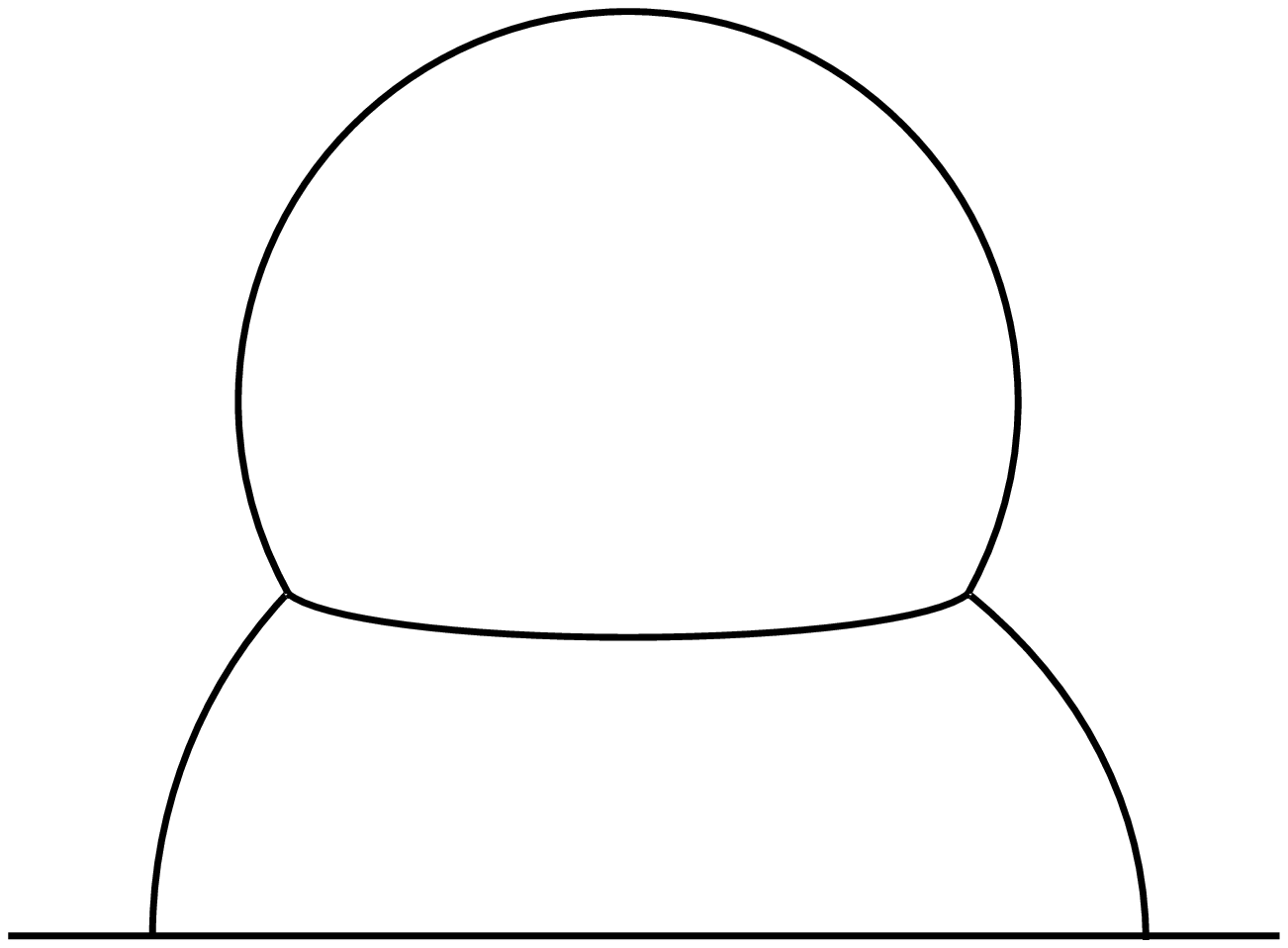  scaled 275}&&\qquad
\BoxedEPSF{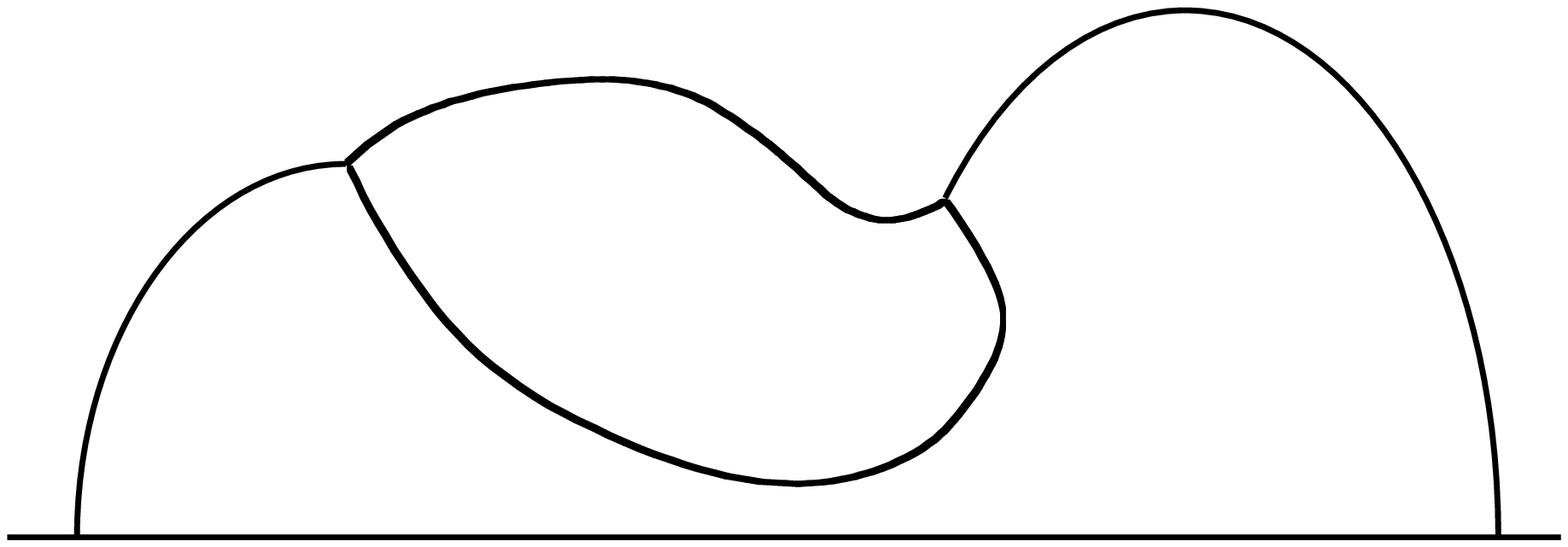 scaled 275}\\
&&\hskip-.75in\BoxedEPSF{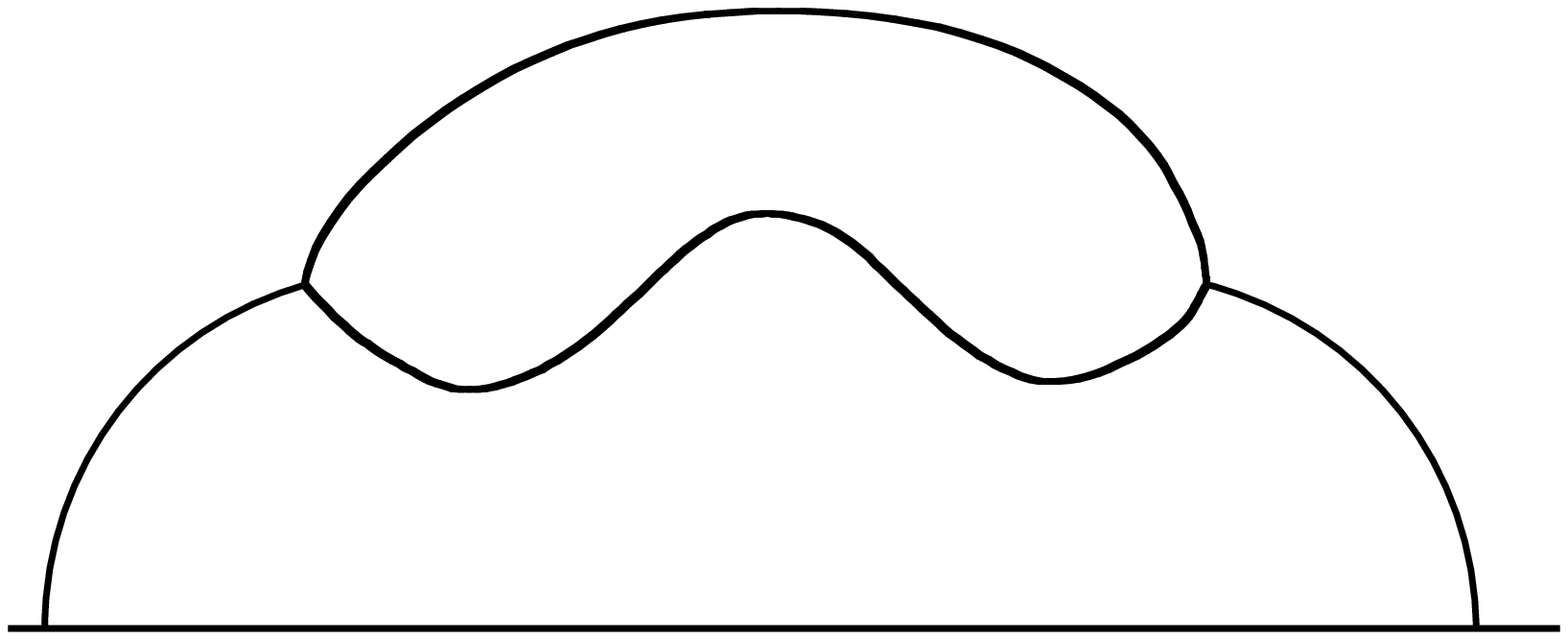 scaled 275}
\end{eqnarray*}
\centerline{{Figure 13}.
 Some eliminated possibilities for cross-sections of torus bubbles.}

\proclaim{Proposition} \label{sum180} If a torus bubble has $\theta_1 + \theta_2 > 180^\circ${\rm ,}
then it is not a minimizer.
\endproclaim 

\demo{Proof} Under this assumption there is a value of $y$ realized by three distinct points
on the circular arcs generating the spherical caps. Since $\theta_2 \ge \theta_1$ we have
in particular that $\theta_2 > 90^\circ$. We know from  Lemma~\ref{120} that $\theta_2 < 120^\circ$.
So $90^\circ < \theta_2 < 120^\circ, 60^\circ < \theta_1  \le 120^\circ$.  Then Lemma~\ref{torus.properties}
implies that the curves $\tau_i$ and $\tau_o$ are both graphs. If $\theta_1 + \theta_2 >180^\circ$ then
there are three points on $S_1 \cup S_2$ with $y$-value equal to $y_2$, 
$(a_1,y_2),(a_2,y_2),(a_3,y_2)$, with $a_1 < a_2 < a_3$ and $a_2 = x_2$. Reflect the part
of the bubble between the planes $x=a_1$ and $x=a_2$ through the plane $x= (x_1 + x_2) / 2$
to get a new piecewise-smooth surface with the same area enclosing the same volumes, as in
the proof of
Lemma~\ref{period.size}, but without needing to worry here about introducing self-intersections. The new
piecewise-smooth surface has three surfaces meeting at an angle not equal to $120^\circ$ and so is
not minimizing, a contradiction. \enddemo

Propositions~\ref{unstable.nodoid}, \ref{period.size},
\ref{sum180}, \ref{torus.properties}, Lemma~\ref{lem.60} and Corollary~ \ref{<120}
combine to rule out many possibilities for the shapes of minimizing torus bubbles.  See Figure~13
for cross-sections of some of these eliminated possibilities.

\proclaim{Lemma} \label{2-1}
For any constant $b${\rm ,} the function $f(c) = b(1-c^2) + c \sqrt {3 (1-c^2)} $ is at most
two\/{\rm -}\/to\/{\rm -}\/one on $(-1,1)${\rm .}
\endproclaim

\demo{Proof} $f'(c) = -2bc - \sqrt 3  g(c)$, where 
$$g(c) = (2c^2 - 1)/ \sqrt {1-c^2} .$$ 
Since $g''(c) = 3(1-c^2)^{-5/2}> 0$, $g$ is convex on $(-1,1)$.

A line intersects a convex curve at most twice, so $f'$ has at
most two roots on $(-1,1)$. Thus $f$ has at most one local maximum and one local minimum
on $(-1,1)$. Combined with $f(0) = f(1) = 0$, this completes the proof. \enddemo

We calculate in the next proposition the value of the forces $f_o$ of $T_o$ and $f_i$ of $T_i$. 

\proclaim{Proposition} \label{force.values}
For the torus bubble determined by $\theta_1$ and $h_o${\rm ,} the
forces $f_o$ of $T_o$ and  $f_i$ of $T_i$ are given by\/{\rm :}\/
$$
f_o = (h_o-1) (1- {c_1}^2) + {c_1} \sqrt {3 (1- {c_1}^2) } 
$$
and
$$
f_i = (h_i-1) (1- {c_1}^2) - {c_1} \sqrt {3 (1- {c_1}^2) } = - f_o \  
$$
where ${c_1} = \cos \theta_1${\rm .}
\endproclaim

\demo{Proof} We first calculate the force $f_o$ of $T_o$.
At $(x_1,y_1)$, $y_1 = \sin \theta_1$ and $\alpha = 5\pi/6 - \theta_1 $ is the angle between
the positive $x$-axis and the graph of $\tau_o$, therefore
\begin{eqnarray*}
f_o & =& h_o {y_1}^2 - 2y_1 \cos \alpha \\
& =& h_o \sin^2 \theta_1 - 2\sin \theta_1 \cos ( 5\pi/6 - \theta_1 ) \\
& =& h_o (1- \cos^2 \theta_1) - 2\sin \theta_1 ((-\sqrt 3 /2) \cos \theta_1 -
(1/2) \sin {(-\theta_1 )}) \\
& =& h_o (1- {c_1}^2) + {c_1} \sqrt {3 (1- {c_1}^2) } - ( 1 - {c_1}^2)\\
& =& (h_o-1) (1- {c_1}^2) + {c_1} \sqrt {3 (1- {c_1}^2) } \ .
\end{eqnarray*}
Similarly, 
\begin{eqnarray*}
f_i & =& h_i {y_1}^2 - 2y_1 \cos \alpha \\
& =& (h_i-1) (1- {c_1}^2) - {c_1} \sqrt {3 (1- {c_1}^2) } \ .
\end{eqnarray*}
Since $h_o-1 = 1-h_i = -(h_i-1)$ by Lemma~\ref{param.bounds}, $ f_o = -f_i $.  \enddemo

\demo{{N}ote} The fact that $ f_o = - f_i \ $ is an example of a general ``balancing principle'' for
constant mean curvature surfaces; see \cite{K-K-S}. 
\enddemo

\proclaim{Lemma} \label{second.angle}
For the torus bubble determined by $\theta_1$ and $h_o${\rm ,} the
angle $\theta_2$ is one of the {\rm (}\/at most\/{\rm )} two solutions of the equation
$$
f(\cos \theta_1) = f(\cos \theta_2),
$$
where
$$
f(c) = ( h_o - 1)(1-c^2) + c \sqrt {3(1-c^2)}.
$$
\endproclaim

\demo{Proof} We calculate the force $f_o$ of $T_o$ at both $(x_1,y_1)$ and $(x_2,y_2)$.

At $(x_1,y_1)$,
$$
f_o = (h_o-1) (1- {c_1}^2) + {c_1} \sqrt {3 (1- {c_1}^2) } ,
$$
where $y_1 = \sin \theta_1$ and $ c_1 = \cos \theta_1$.
A similar calculation at $\theta_2$ shows
$$
f_o = (h_o-1) (1- {c_2}^2) + {c_1} \sqrt {3 (1- {c_2}^2) } .
$$
Setting the two forces equal and applying Lemma~\ref{2-1} gives that there
are at most two possible values for $ c_2  = \cos \theta_2$.  Since $0 <  \theta_2 < 180^\circ$, where the
cosine function is monotonically decreasing, this implies that there are at most two possible
values for $ \theta_2$. \enddemo

One solution occurs when the angles are equal. An example of such a symmetric torus bubble,
depicted in Figure~7, occurs with $\theta_1 = \theta_2 = 90^\circ$, and $h_o$
approximately $1.9848$. However these symmetric solutions
never lead to a minimizing torus bubble by Corollary~\ref{theta_1<theta_2}.

Recall that we remarked after Corollary~\ref{theta_1<theta_2}
that we can without loss of generality restrict attention to nonsymmetric torus bubbles.
This means that we are working in the case where $\theta_1 < \theta_2$. 
For each $\theta_1$, $h_o$, there is at most one value of $\theta_2$ which can potentially
lead to a stable torus bubble.

\proclaim{Proposition} \label{cotan} If a minimizing torus bubble has $ 60^\circ \le \theta_1 \le 120^\circ$
then $h_o > 1-\sqrt 3 \cot \theta_1$ and $ h_i < 1 + \sqrt 3 \cot \theta_1${\rm .}
The same statement holds with $\theta_2$ replacing $\theta_1${\rm .}  
\endproclaim

\demo{Proof} We again apply Lemma~\ref{lem:trap} to compare the torus bubble with a double
bubble whose smaller spherical cap $S_o$ subtends an angle of $\theta_1$, as in the proof
1 Lemma~\ref{120}.
See Figure~14.
If $h_i$ is larger than the mean curvature of
$S_i$, then $T_i$ is a nodoid
contained inside the smaller component of the double bubble, as shown in Lemma~\ref{lem:trap}. 
See Figure~14. It follows that $\tau_o$ cannot remeet $\tau_i$
and the two curves cannot generate a torus component. 
\enddemo

\figin{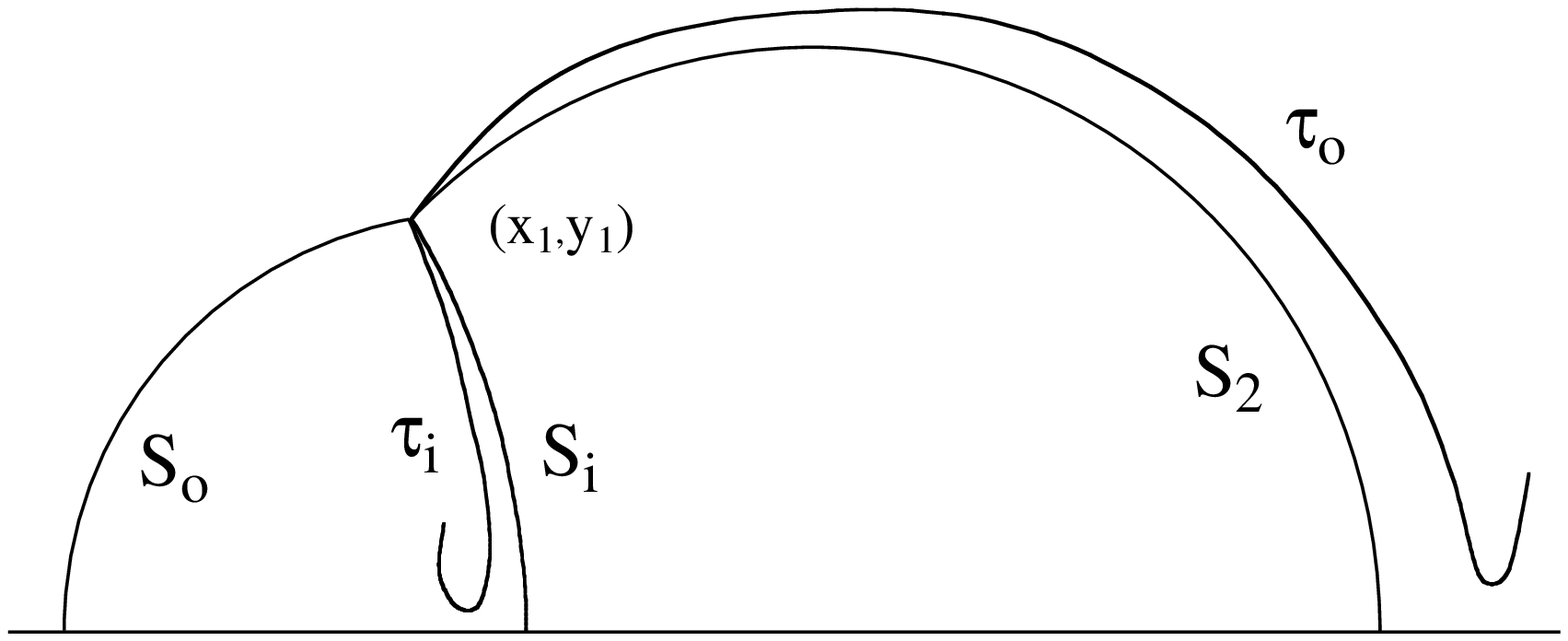}{400}
\smallbreak
{{Figure 14}. The mean curvature of $T_o$ cannot be too small if $\tau_o$ and $\tau_i$

 are
to meet again.}

\proclaim{{C}orollary} \label{T_o.is.nodoid}
$T_o$ is a nodoid{\rm .}
\endproclaim

{\it Proof}. If $ 60^\circ \le \theta_1 \le 120^\circ$ then Proposition~\ref{cotan} shows
that the mean curvature of $T_o$ is greater than that of the sphere which is tangent to it at
$(x_1,y_1)$, so Lemma~\ref{lem:trap} implies that $T_o$ is a nodoid. If $ \theta_1 <
60^\circ$ then $\tau_o$ starts out heading to the left. It must become
vertical, hence $T_o$ must be a nodoid. \phantom{poor puppy} \hfill\qed
\proclaim{Lemma} \label{circ.under.nodoid}
Let $N$ be a subcurve of a nodary generating a
nodoid of mean curvature $h>0${\rm ,} running
between a vertical tangency on the left at $(x_v,y_v)$
and a vertical tangency on the right at
$(x_v',y_v')${\rm ,} with a unique local maximum at $(x_M,y_M)${\rm . } Then a circle
$C$ of curvature $h$ which is tangent to $N$ at $(x_v,y_v)$ lies underneath $N${\rm .}
Moreover $x_v - x_m$ is greater than the radius $1/h$ of the circle{\rm .}
\endproclaim

{\it Proof}. In a deleted neighborhood of $(x_v,y_v)$, the curvature of $N$ is less than that of
$C$, so that $C$ locally lies beneath $N$. If $C$ crosses $N$ between $x_v$ and $x_v'$, let
$C'$ be the subcurve of $C$ starting at $(x_v,y_v)$ and running to the first point of intersection
of $C$ with $N$. Let $C_t' = C' + (0,t)$ be a vertical translate of $C'$ by a distance   $t$
along the $y$-axis.
Let $T = \sup \{ t : C_t' \cap N \ne \emptyset \}$.  Thus $C_T'$ is tangent to
$N$ at a point $P$ and lies above $N$ near $P$. The curvature of $N$ at
$P$ is given by $ k_m= h - k_p < h$. Since the curvature of $C_T'$ is greater than that of
$N$, it cannot lie above it, and it follows that $ x_M - x_v$ is larger than the radius of
$C$, i.e., $x_M - x_v > 1/h.$   \hfill\qed\medbreak

The next proposition rules out the existence of minimizing torus bubbles with $h_o$
and $\theta_1 $ both close to $0$. The numerical calculations that we will apply are badly
behaved in this region, so we use this geometric argument to exclude it.

\proclaim{Proposition} \label{0.2}
If a torus bubble has $ h_o \le 0.2$ and $\theta_1 \le 5.7^\circ${\rm ,} then it is not a minimizer{\rm .} 
\endproclaim

\demo{Proof} We will show that before $\tau_i$ and $\tau_o$
rejoin at $(x_2,y_2)$, either $\tau_o$ has two vertical tangencies, or
$\tau_i$ traverses more than a full period.
Note that our assumption implies that $h_i \ge 1.8$ 

Proposition~\ref{prop.Delaunay} implies that $\tau_i$ is a graph and its period is no
longer than $2\pi/h_i < 3.5$. The bubble is not a minimizer if $ x_2-x_1 > 3.5$ by
Proposition~\ref{period.size}. The initial angle of $\tau_i$ from the positive $x$-axis is $30^\circ -
\theta_1$, so it cannot turn (anti-clockwise) through an angle of more than
$60^\circ + \theta_1$. When they meet again $\tau_i$ and $\tau_o$ intersect at an angle of
$120^\circ$. It follows that $\tau_o$ does not intersect $\tau_i$ before it turns (clockwise) through
an angle greater than $180^\circ - \theta_1$. Therefore the two curves cannot intersect before the
nodary $\tau_o$ reaches its maximum, which happens when $\tau_o$ has
turned through an angle of $150^\circ - \theta_1$. Let $(x_M,y_M)$ be the point where $\tau_o$
reaches its maximum and $(x_v,y_v)$ the point where $\tau_o$ first goes vertical. Then 
$x_2 - x_1 > x_M - x_1$.

Now $ h_o = k_p + k_m$ where $k_m$ is the curvature of the
generating curve in the $xy$ plane and $k_p $ is the principal curvature due to rotation
around the axis. At $(x_v,y_v)$, $k_p = 0$ and $h_o = k_m \le .2$ by hypothesis.

Lemma~\ref{circ.under.nodoid} implies that $x_M - x_v > 1/h_o \ge 5 .$ 

The assumptions $h_o \le .2$ and $ \theta_1 \le 5.7^\circ$ give an upper bound for the force $f_o$
of $T_o$. Calculating at height $y_1 = \sin \theta_1  <\theta_1 \le  5.7 \pi / 180 < .1$ gives

\begin{eqnarray*}
| f_o | & = &  | {y_1}^2 h_o - y_1 \cos { (30^\circ + \theta_1 )}  | \\
& < & | {y_1}^2 h_o - y_1 \sqrt 3 /2 |  \\
& <  &  (.1)^2 h_o + .1 \sqrt 3 /2 \\
& < & .09 .
\end{eqnarray*}

We bound the value of $y_v$ from above by solving $f_o = y_v^2 h_o$.  
$$
 y_v^2 = f_o/h_o < .09/h_o,
$$
so
$$
y_v < \sqrt {.09/h_o}.
$$

We next bound $x_1 - x_v$ from above by
noting that $\tau_o$ is convex between $x_v$ and $x_1$ by Lemma~\ref{prop.Delaunay}.
This implies that $(x_v,y_v)$ lies to the right of the point $(x,y_v)$ at which the line
of slope $\dot {\tau_o}(x_1)$ has height $y_v$, 
so that 
\begin{eqnarray*}
|x_v - x_1| < |x-x_1| = |(y_v - y_1) / \dot {\tau_o}(x_1) | < \sqrt 3 (y_v - y_1) \\
 < \sqrt 3 y_v < \sqrt 3 \sqrt {.09/h_o} = \sqrt {.27/h_o}.
\end{eqnarray*}

This implies that
$$
x_M - x_1 =  (x_M - x_v) - (x_1 - x_v ) > 1/h_o -  \sqrt {.27/h_o}.
$$

For $ 0 < h_o \le .2 $, differentiation shows that the function 
$1/h_o - \sqrt {.27/h_o}$ is decreasing and its minimum value occurs at $h_o = .2$, so that
$x_M - x_1 > 1/(0.2) - \sqrt {.27/(0.2)} > 3.8$.  Thus $\tau_o$ travels a distance of at least $3.8$
to the right before hitting $\tau_i$. But the period of $\tau_i$ is less than $3.5$.
We have shown that either $\tau_o$ changes direction and has two vertical tangencies
or $\tau_i$ contains a full period.
In the first case it is not a minimizer by Lemma~\ref{lem.60},
and in the second case it is not a minimizer by Proposition~\ref{period.size}.
\enddemo

The next lemma will be used to give bounds on the size of a circle enclosing a loop of a
nodary.

\proclaim{Lemma} \label{curve.trap} If an embedded curve $\gamma$ in the
positive quadrant of the plane is tangent to the $y$\/{\rm -}\/axis at the origin and has
curvature $k > 1/r > 0${\rm ,} then it does not meet $C(r)${\rm ,}
the circle of radius $r$ centered at $(0,r)${\rm ,} at any point other than the origin{\rm .}
\endproclaim

{\it Proof}. If it does cross $C(r)$, denote by $\gamma'$ the subcurve of $\gamma$ starting at
the origin and running to the first point of intersection with $C(r)$, as
in Figure~15. The subcurve $\gamma'$ is a graph over the $x$-axis since it
is embedded. Let $\gamma'_t = \gamma' + (0,t)$ be a vertical translate of
$\gamma'$ by a vertical distance of $t$. Let $T = \sup \{ t : \gamma'_t \cap C(r) \ne
\emptyset \}$. Then $\gamma'_T$ meets $C(r)$ at a point $P$ in its interior,
but not the interior of $C(r)$, and thus is tangent to $C(r)$ at $P$.
Since the curvature of $\gamma'_T$ is greater than that of $C(r)$, this
is a contradiction. \hfill\qed

\figin{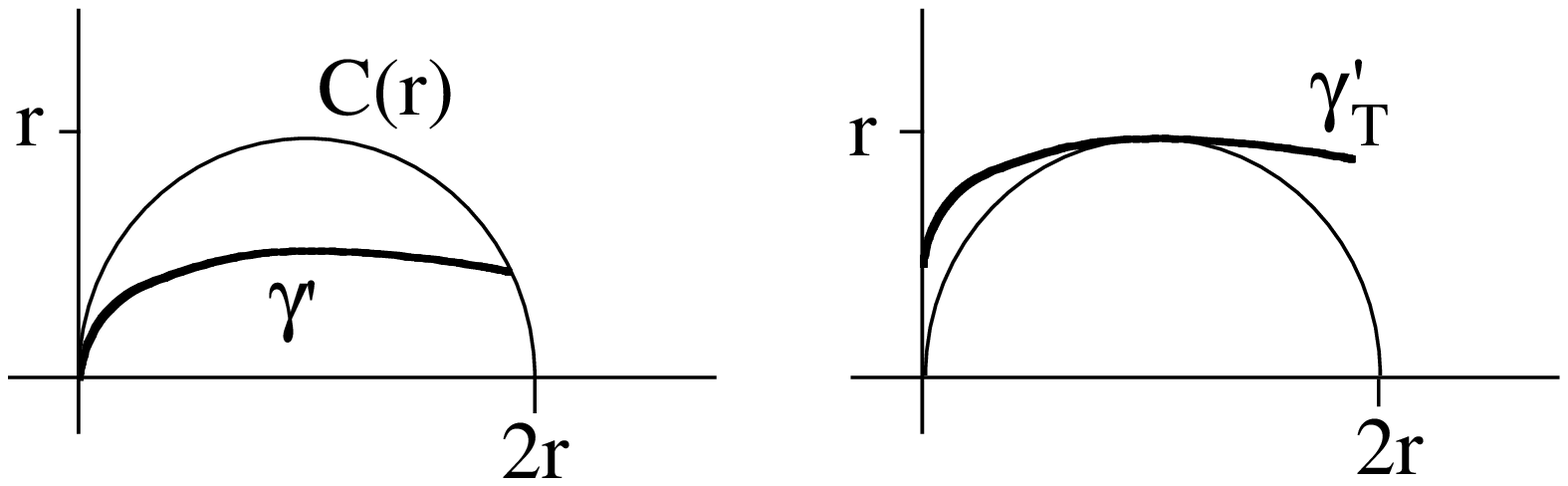}{400}

 {{Figure 15}. The curvature of the curve $\gamma$ cannot be greater than $1/r$

 if it crosses
$C(r)$ twice.}

\proclaim{Lemma} \label{circles.trap} If $\tau_o$ has curvature greater than $k${\rm ,}
and $\tau_i$ has curvature less than $-k${\rm ,}
then their union is contained in the intersection of two circles of radius $r = 1/k $ meeting at $120^\circ${\rm .}
\endproclaim

{\it Proof}. Construct circles $C_1$ and $C_2$ of radius $1/k$, meeting at $120^\circ$
at $(x_1,y_1)$, with $C_1$ tangent to $\tau_o$ and $C_2$ tangent to $\tau_i$
as in Figure~16.
Let $D$ denote the line segment joining the two intersection points of $C_1$ and $C_2$.
By Lemma~\ref{curve.trap} we know that each of $\tau_o$ and $\tau_i$ cannot leave the intersection of
these two circles without first crossing $D$.  Since each is convex in $C_1 \cap C_2$ it follows
that they remeet before either one leaves $C_1 \cap C_2$.
It follows that the intersection of the torus component
with the upper half-plane is contained in $C_1 \cap C_2$.
\hfill\qed

\begin{center}
\BoxedEPSF{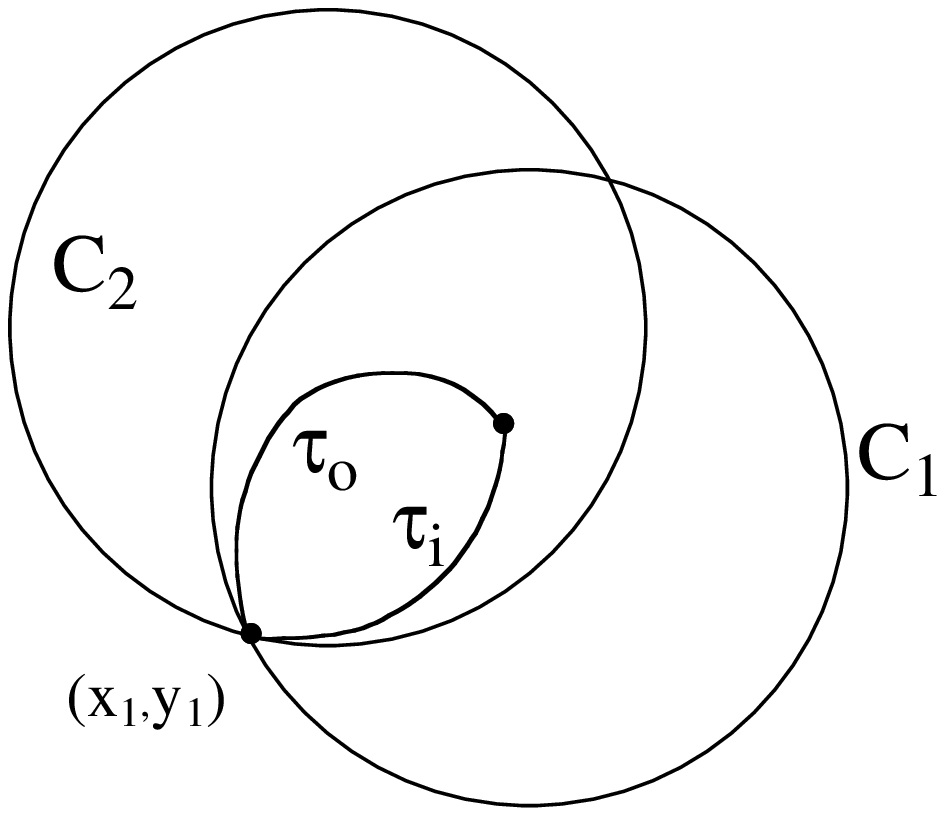 scaled 400}
 \end{center}

 {{Figure 16}. The torus component generating curves $\tau_i$ and $\tau_o$
are

 trapped inside two circles meeting at $120^\circ$.}
 
\proclaim{Lemma} \label{trap1} If a torus bubble has $h_i \le -k < 0$ and the intersection of
the torus component with the upper half\/{\rm -}\/plane has minimum $y$\/{\rm -}\/value equal to
$y_0$ with $y_0 > 1/k ${\rm ,} then the torus component contains volume
$$ 
v \le 2.5 \pi (y_0 + \frac{\sqrt 3}{2}\frac{y_0}{ky_0 - 1}) ( \frac{y_0}{ky_0 - 1})^2 \ .
$$ 
\endproclaim

\demo{Proof} 
In a torus bubble with $h_i \le -k < 0$ we have that $|h_o| = -h_i + 2  \ge~k$.
For a Delaunay curve passing through a point $(x,y)$ in the upper
half-plane,
\begin{eqnarray*}
|k_m| & \ge & |h| - |k_p| \\
& = & |h| - \frac{1}{ y \sqrt{1 + \dot {y}^2 }}\\
& \ge & |h| - \frac{1}{ y}.
\end{eqnarray*}
For the curves $\tau_i$ and $\tau_o$ we have $y  \ge y_0$, so that 
\begin{equation}
|k_m| \ge k - \frac{1}{ y_0}. \label{eq:k_m}
\end{equation}
So the curves $\tau_o$ and $\tau_i$ each have curvature with absolute value greater than $(ky_0 - 1)/{y_0}$. 

By Lemma~\ref{circles.trap} the intersection of the torus component with the upper half-plane is
contained in the intersection $w$ of two circles of radius
$r = y_0 / (k y_0 - 1)$ meeting at an angle of $120^\circ$.  The area of $w$ is 
$$ 
a = (2 \pi/ 3 - \sqrt 3 / 2) r^2 < 1.25 ( \frac{y_0}{ky_0 - 1})^2.
$$ 
The diameter of $w$ is 
$$
d =  \sqrt 3  \frac{y_0}{ky_0 - 1} .
$$
Since $w$ is centrally symmetric and has minimal $y$-value smaller than or equal to $y_0$, the center of mass
of $w$ has $y$-coordinate at most $y_0 + d/2$.
Pappus' theorem for the volume of a solid of revolution obtained by rotating this area
around the $x$-axis gives that 
\begin{eqnarray*}
v & \le & 2 \pi (y_0 + \frac{d}{2}) 1.25 ( \frac{y_0}{ky_0 - 1})^2,\\
& \le & 2.5 \pi (y_0 + \frac{\sqrt 3}{2}\frac{y_0}{ky_0 - 1}) ( \frac{y_0}{ky_0 - 1})^2.\\
\noalign{\vskip-36pt}
\end{eqnarray*}
\enddemo

\phantom{good lunch}

We apply these volume estimates to reduce the range of possible mean curvatures in a
minimizing equal-volume torus bubble.
 
\proclaim{Proposition} \label{h>-8}
In a minimizing equal\/{\rm -}\/volume torus bubble{\rm ,} $h_o \le 10${\rm .}
\endproclaim

{\it Proof}. The hypothesis $h_o \le 10$ is equivalent to $h_i \ge -8$.  
Suppose to the contrary that $h_i < -8$ in a minimizing equal-volume torus bubble. 

Proposition~\ref{torus.properties} states that $\tau_i$ is a graph, so
that the volume $v_i$ of the
ball component $B$ is strictly  larger than the volume
of the region under $S_2$. Lemmas~\ref{120} and
 \ref{lem.60} imply that $60^\circ < \theta_2 < 120^\circ$.
It follows that $v_i > 5 \pi / 24 $,
the volume under a unit radius spherical cap subtending an angle of $60^\circ$.
We now estimate the
volume $v_o$ of the torus component.
The curves $\tau_i$ and $\tau_o$ generating the torus
component contain the point $(x_2,y_2)$
where $\sqrt 3/2 \le y_2 \le 1 $. Their maximum and
minimum $y$-values differ by less than $2/|h_i| \le 0.25$
by Proposition~\ref{prop.Delaunay},
and therefore both curves are contained in
$\{ y \ge \sqrt 3 /2  - 0.25\}$. The expression
$\displaystyle \frac{y_0}{ky_0 - 1}$ is decreasing with respect to
$y_0$ when $y_0 > 1/k$, and decreasing with respect to
$k$, and also $0.6 < \sqrt 3 /2 - 0.25$. We apply Lemma~\ref{trap1} with $k \ge 8$
and $y_0 \ge 0.6$, giving that
\begin{eqnarray*}
v_o & < & 2.5 \pi (y_0 + \frac{\sqrt 3}{2}\frac{y_0}{ky_0 - 1}) ( \frac{y_0}{ky_0 - 1})^2 \\
 & < & 2.5 \pi (1 + \frac{\sqrt 3}{2}\frac{0.6}{4.8 - 1}) ( \frac{0.6}{4.8 - 1})^2 \\
 & < & 0.08  \pi  \ .
\end{eqnarray*}

Thus $v_o < 0.08  \pi <  5 \pi / 24 < v_i$, contradicting the assumption that the torus bubble
encloses equal volumes. \hfill\qed \medbreak

Our algorithm will need the following rather technical result about the
slope of a generating curve.

\proclaim{Lemma} \label{lem:convex.unduloid}
Suppose that $\tau_i$ satisfies $ h_i \le 2$ and $0 \le
\theta_1 \le 4^\circ${\rm .}  Let $\alpha$ be the angle made by $\tau_i$ with the positive $x$\/{\rm -}\/axis and
let $\alpha_1$ be the value of this angle at $(x_1,y_1)${\rm .} Then
$\alpha \ge \alpha_1$ for as long as $\tau_i$ remains below height $\{ y = 1/2 \}${\rm .} In
particular{\rm ,} $\alpha \ge 26^\circ$ below this height{\rm .}
\endproclaim

{\it Proof}. 
First note that if $h_i \le 0$ then $\tau_i$ is a nodary or catenary, and the angle 
 made by $\tau_i$ with the positive $x$-axis is increasing as one moves along $\tau_i$, which implies
the conclusion of the lemma. So we can assume that $0 < h_i \le 2 $.

Let $y$ be any point on $\tau_i$. Then 
$$
h_i y^2 - 2 y \cos \alpha - f_i = 0 \ .
$$
A given value of $\cos \alpha$ is realized by at most two $y$ values solving this quadratic equation. 
Suppose that as the $\tau_i$ goes from a minimum to a maximum, the angle 
$ 26^\circ \le \alpha_1 = 30^\circ - \theta_1 \le 30^\circ$ occurs twice. The two $y$-values which share
the same
$\alpha_1$ both solve the above equation. Thus they sum to $2   \cos \alpha_1 / h_i  \ge 2 (.8)
/ 2 \ge .8$. Since the smaller root has $y$-coordinate $\sin (4^\circ) < .1$ the larger root is greater
than $.7 > 1/2$. The
angle $\alpha$ is increasing near the smaller root. Thus
$\tau_i$ subtends an angle of at least $\alpha_1$ with the $x$-axis as long as $\tau_i$
remains below height $1/2$.  Since  $\alpha_1 = 30^\circ - \theta_1 \ge 26^\circ$, the final conclusion
holds. \hfill\qed

\section{Computation} \label{Computation} 

We have shown that to solve Conjecture~\ref{conj-double.bubble} in the equal volume case
it suffices to show that certain torus bubbles are less efficient at
enclosing two equal volumes then a double bubble.
In this section we will show how these torus bubbles are ruled out by a computation.
The algorithm for this computation is given below.

The basic idea is to consider a domain of torus bubbles corresponding to a product of small
intervals in each of $\theta_1$ and $h_o$, and to calculate as much as possible
about the geometry of the corresponding torus bubbles.
Various calculations are then applied which rule these torus
bubbles out as potential minimizers.
The accuracy of these calculations depends on the size of the domain
rectangle we start with. The computational scheme
will succeed if these can be chosen small enough to get sufficient
accuracy, yet large enough that a reasonable number of them
cover all the possibilities.

The computation was performed using double precision floating point
numbers. The fundamental data type used is the IEEE 754 32-bit real
number, see \cite{ANSI}.
This is a binary representation with a 23-bit mantissa (plus an
implied leading bit), an 8-bit binary exponent, and a sign bit.
The (finitely many) real numbers that can be represented in this standard
with no error are called {\it representable}.
The IEEE standard specifies that
the add, subtract, multiply, divide, and square root operations be performed as if done
exactly and rounded to the nearby representable number according to
the rounding mode in effect. There are three rounding modes that can
be chosen: up, down, or nearest. The values $ +\infty$ and $ -\infty$ are
representable and behave in a specified way. Floating point
exceptions are masked, but flags are sticky and available for clearing
and inspection.
Most computers in use today implement the IEEE standard.
We denote the result of rounding a real number $x$ down to a representable number
smaller or equal to $x$ by $\underline{x}$, and similarly denote by 
$\overline{x}$ the result of rounding $x$ up to a representable number
greater or equal to $x$.

Combined with the methods of interval arithmetic,
see Moore \cite{M} and Alefeld-Herzberger \cite{AH},
the IEEE standard allows numerical calculation with exact bounds on
accuracy. Interval arithmetic is a method by which a real-valued
function on the reals can be extended to an interval-valued function
of intervals.
An interval is formed from two representable reals. 
Mathematically, it represents the closed interval between the two reals. The add,
subtract, multiply, divide, and square root functions are extended to
intervals by the IEEE operations on reals along with directed rounding.
We only use the round-to-nearest mode in
calculating averages, (see the procedure avgwt below),
as these are only used to divide intervals and
do not need to have bounds calculated.

Generally, we say that an interval $X$ is an extension of a real number $x$ if $x \in X$ and
an interval-valued function of intervals $F$ is an extension of $f$ if $x \in X \Rightarrow
f(x) \in F(X)$.

In this section, real numbers are denoted by lower case, and intervals
with representable real endpoints by upper case.
The reals embed into the set of intervals by mapping a real
to the smallest interval containing it.
The lower and upper bounds to an interval are denoted
with lower and upper bars, so $X = [\underline{X},\overline{X}]$.
Arithmetic operations on representable reals and intervals
are interpreted according to IEEE and interval rules,
not by the usual mathematical definitions.

As an example, the sum of two intervals $A = [\underline{A},\overline{A}]$
and $B = [\underline{B},\overline{B}]$ is given by
$$
A+B = [\underline{ \underline{A}+\underline{B}}, \overline{ \overline{A}+ \overline{B}} ],
$$
where the left endpoints are rounded down when added,
and the right endpoints are rounded up when added.

Other operations on reals are extended to intervals
as in \cite{M} and \cite{AH}.
Relations are interpreted positively, so for example when $X$ and $Y$ are intervals, $X < Y$
means that for any $x \in X$ and any $y \in Y$, $x < y$, and 
$X \ne Y$ means that for any $x \in X$ and any $y \in Y$, $x \ne y$.

We have equivalent expressions
$X<Y \Leftrightarrow \overline{X} < \underline{Y} $,
and $X \ne Y \Leftrightarrow X \cap Y = \emptyset$.
The union of intervals $X$ and $Y$ is somewhat nonstandard; it is
the smallest interval containing both sets,
and is denoted $X \overline{\cup} Y$. The interval-valued functions {\bf
Absolute Value($X$), Max($X,Y$), Min($X,Y$),} and the
intersection of two intervals $X \cap Y$ are defined in the standard way,
without rounding.
There is no IEEE standard for transcendental functions,
so we designed our program to avoid all calls to trigonometric functions.

The operations of the IEEE floating point
standard can lead to undefined operations. For example
the quotient $0/0$ and the product $0 \times +\infty$ result
in an output of NAN (Not A Number), and are signaled by the
presence of an exception flag.
Operations of interval arithmetic can also lead to questionable operations, for example
when the quotient of two intervals $A/B$ is calculated and $B$ contains $0$.
Our implementation of the division operator on intervals returns
$[-\infty, +\infty]$ if $B$ contains 0.
The algorithm draws no conclusions when such a division
occurs, but rather calls for a subdivision of the input
into intervals of smaller size,
where the operations are repeated with greater accuracy.
Exceptions such as overflow and underflow,
sometimes an issue in computer assisted proofs,
are not an issue in our algorithm. In any case they do not occur.

Bounds obtained from repeated applications of
interval arithmetic are potentially far from sharp,
especially when using wide input intervals.
If just a single value needs to be computed rigorously,
one can work with intervals whose width is comparable to the roundoff error
of a single operation. For our purposes these types of intervals are far too
thin, and would make the number of calculations we need impractical.
Our intervals are ``fat'' in the
terminology of Fefferman \cite{Fefferman}, meaning that they are sized by the scope
of the problem we are solving rather than the size of the computational rounding.
Often the intervals are wide enough to make a
perfectly good formula look like nonsense. Sometimes it is possible to
narrow the intervals under consideration because we have
knowledge about what the possible legitimate values are that could
arise during a computation. For example, consider a
formula involving square roots. Assuming we're not using imaginary
numbers, the validity of the formula presupposes that the argument to
the square root is nonnegative. However, when we pass to an interval
extension of the formula, the interval argument to the square root
will often include negative values, even though these values cannot
arise from the problem we are interested in. The square root of a negative
number would normally trigger an exception,
but the nature of the calculation allows us to
define the interval square root function {\bf Sqrt} to
discard any negative portion of an interval argument. The
justification for this apparent sleight-of-hand rests on the validity of the
original formula over the reals. If our theorems tell us that the quantity whose square root
we are taking is nonnegative,
we are justified in truncating the interval to exclude negatives.

The interval function which always returns $[-\infty, +\infty]$ can be used to
represent any function, but not very usefully. We use this interval to return the value of a
division by an interval containing zero.
An interval $[a,b]$ is empty if $a > b$. All of our functions take 
nonempty interval arguments and return nonempty interval
results except for the intersection operation, which returns an
interval which can be empty.

The value of the nonrepresentable constant $\sqrt 3 $ in our
algorithm is expanded to a narrow interval.
Other constants used are representable.
Thus we use the expression $ 5 H_o  \le 1$ rather than
$H_o \le  0.2$, which involves the nonrepresentable constant $0.2$.
Representable constants such as 5 behave identically in interval arithmetic
to intervals where both endpoints have the same value, such as
$[5,5]$, and operations involving both intervals and representable constants
are treated as if the constant was first converted to a one-point interval.

The {\bf verify} statement is used to assure that a particular condition
holds. If the condition fails, then the entire program is stopped
(aborted). If no such failure occurs, the presence of such a {\bf verify}
statement constitutes a proof that the condition holds for the given
inputs.
One use of the {\bf verify} statement is in the procedure
{\bf avgwt} defined below, which gives a weighted average of two intervals.
To apply this we want to be sure that, as we expect in the cases
where it is applied, the first interval is strictly smaller than the second.
We check that this is so with a  {\bf verify} statement.
Running our program results in no violations of any statement which is tested
with a  {\bf verify} statement.

The program examines the set of all torus bubbles to see if any can
be a minimizer.
A range of hypothetical torus bubbles is specified by intervals $\Theta_1$ and $H_o$,
according to Proposition~\ref{determined}.
To avoid unnecessary use of trigonometric functions
our algorithm uses an equivalent parameterization of the space of torus bubbles by
intervals $Y_1$ and $H_o$,
where $Y_1 = \sin \Theta_1$. Rather than solve the mean curvature
differential equation directly to find intervals $(X_2,Y_2)$ containing the point $(x_2,y_2)$
where $\tau_i$ and $\tau_o$ intersect,
we assume existence of a torus bubble, derive $Y_2$, and then deduce
$X_2$ from numerical integrals.
Volume calculations for $V_i$ and $V_o$ are obtained from
additional numerical integrals.
At each stage in the calculation, we check whether the torus bubbles can be
rejected based on the instability results of Section~\ref{sec.torus}. If not, we
test whether the volumes are equal.
As it turns out all torus bubbles are rejected for these
reasons, so area calculations are not necessary and are not present in the
algorithm we present.
Most torus bubbles are rejected either because $\Theta_2$ is out of range, or
because the $x$-displacement of $\tau_i$ and $\tau_o$ differ,
or because the volumes enclosed in the two regions differ.

The computation of $\Theta_2$ and $Y_2$ is based on an analysis of the force function.
The integration process used to calculate the value
of $X_2$ is complicated by the fact that
the curves may not be graphs, so we have to choose an appropriate parameterization.
The most convenient parameterization is in terms of
$y$, because the ODEs satisfied by Delaunay curves
involve only $y$, and because the $y$-coordinates of
the endpoints are specified by $\Theta_1$ and $\Theta_2$.
This allows $x$-displacement and volume
to be expressed directly as integrals in terms of $y$, making
it unnecessary to actually generate an ODE solution which produces
$\tau_i$ and $\tau_o$.

The curves $\tau_i$ and $\tau_o$ are generally not graphs as functions of $y$,
so the integrals we use have singularities at points where $\dot{y} = 0$.
A change of variables resolves the singularity,
as worked out in Proposition~\ref{prop:dxcrit}.

Another difficulty is that the geometry sometimes degenerates at the boundary of the regions we
are examining. For example, we need to exclude torus bubbles with $\theta_1$ arbitrarily
close to $0$, but some of our formulas become singular when $\theta_1= 0$, and don't make
sense there. If the interval $\Theta_1$ contains $0$, then various intervals representing
$y$-coordinates along the Delaunay curves will also contain $0$,
and the integrals give results which will not allow torus bubbles to be
eliminated.
We get around this problem by using crude estimates which apply in a
small zone near the $x$-axis. Near the $x$-axis we have some information about
the initial angle and the concavity of the Delaunay curves, and this suffices
to give crude bounds to their behavior near the axis.
We compute integrals only when
at an appropriate distance above the $x$-axis. 
This approach suffices for all cases except for when $h_o$ is also near $0$.
In that case we rely on Proposition~\ref{0.2}.

We now describe in pseudo-code the main
interval procedures used in our algorithm.
As noted above, aside from the procedures described in this section
our algorithm makes use only of standard arithmetic operations on intervals,
namely addition, subtraction, multiplication, division, absolute value,
maximum, minimum and square.

The procedure {\bf avgwt} calculates a number in between two inputs,
used later to subdivide intervals. It returns a representable real number which is a weighted average. For convenience, we
set the rounding mode to round to nearest in this procedure. In some applications of this procedure, an interval
is assigned the output of the avgwt procedure, eg in Step 7 of the procedure
{\bf DivideAndCheckRectangle}. $\!$ The resulting interval
then consists of a single point. \medbreak

\setlength{\parindent}{0in}

{\tt real procedure avgwt \\
input: X, Y, w
  \begin{tabbing}
Col 1 \= Col 2 \= Col 3 \= \kill
\> verify $X < Y$\\
\> $z := (1-w) \overline{X} + w \underline{Y} $ \\
\> return $z$
\end{tabbing}
}

\setlength{\parindent}{.2in}

The interval procedure {\bf Sqrt(X)} returns an interval containing
$\{ \sqrt x : x \in X$ and $x \ge 0 \}.$
The real procedure {\bf width} calculates an upper bound for an interval's width. \smallbreak

The interval procedure {\bf Compare}  checks whether the interval $X$ is to the left of the
interval $Y$, returning $A$ if yes, $B$ if $X$ is to the right of
$Y$, and the union $A \overline{\cup} B$ if the intervals overlap. \vglue6pt

\setlength{\parindent}{0in}
{\tt
interval procedure Compare \\
input: $ X, Y, A, B$    
\begin{tabbing}
Col 1 \= Col 2 \= Col 3 \= \kill
\> if $X < Y$ then return $ A$     \\
\> else if $X > Y$ then return $B$     \\
\> else return $A \overline{\cup} B$    \\
\end{tabbing}
} 
\setlength{\parindent}{.2in}

\proclaim{Lemma} For any $x, y, a, b,$ in any $ X, Y, A, B,$ respectively{\rm ,} if $x \le y$
then $a \in$ {\bf Compare}$(X,Y,A,B)${\rm ,} and if $x \ge y$ then $b \in$ {\bf Compare}$(X,Y,A,B)${\rm .}
\endproclaim

{\it Proof}. The lemma follows from a straightforward case by case analysis.
\hfill\qed
\medbreak
\setlength{\parindent}{.2in}

{\bf Integrate} is the basic method for numerical integrals.
It gives relatively wide intervals, but
is good enough for our purposes. Its input is an interval-valued function
of one interval variable, and two intervals which give the limits of
integration.\smallbreak
 
\setlength{\parindent}{0in}
{\tt
interval procedure Integrate \\
input: $F, A, B$       
\begin{tabbing}
Col 1 \= Col 2 \= Col 3 \= \kill
\>$ R := 0$    \\
\> $H := (\underline{B} - \overline{A}) / 32$    \\
\> if $(\overline{H} > 0)$ then       \\
\> \>           begin    \\
\> \>           for $ i := 0$ to 31 do    \\
\> \> \> $R := R + F(\overline{A} +(0 \overline{\cup} H)+  i H )$  \\
\>      \> end \\
\> return $0\overline{\cup}F(A)$Width$(A) + 0\overline{\cup}F(B) $Width$(B) + R H$
\end{tabbing}
}
\setlength{\parindent}{.2in}

\proclaim{Proposition}
For any integrable $f$ and real numbers $a, b,$ with $a \le b${\rm ,} and interval
extensions $F, A, B${\rm ,} the integral of $f$ from $a$ to $b$ is contained in
${\rm Integrate}(F,A,B)${\rm .}
\endproclaim

\demo{Proof} If $A \le B$,
$$
\int_{a}^{b} f \ dx  =   \int_{a}^{\overline{A}}  f  \  dx  
+ \int_{\overline{A}}^{\underline{B}} f  \  dx   
+ \int_{\underline{B}}^{b}  f  \  dx   \  .
$$
The term
$ \displaystyle \int_{a}^{\overline{A}} f  \  dx $
is contained in  $0 \overline{\cup}F(A) {\rm Width}(A)$ and 
similarly $ \displaystyle  \int_{\underline{B}}^{b} f  \  dx $
is contained in $0\overline{\cup}F(B){\rm Width}(B)$.
Note that we need to replace $F(A)$ and $F(B)$
by  $0\overline{\cup}F(A)$ and $0\overline{\cup}F(B)$
to deal with the issue that $a$ can lie anywhere in $A$ and similarly $b$ in $B$.
Taking the union of $F$ with 0 ensures that
our sum contains the intervals $F([a,\overline{A}])$ and $F([\underline{B},b])$.

The value of 
$ \displaystyle \int_{\overline{A}}^{\underline{B}} f \ dx $ 
is bounded by lower and upper Riemann sums, using 32 equally sized intervals.
The term $ F(\overline{A} +(0 \overline{\cup} H)+  i H) $ evaluates
$F$ on an interval which contains $[\overline{A}+ih, \overline{A}+(i+1)h]$ for any $h \in H$,
and so contains both an upper and a lower bound for $F$ on this interval. 
The $RH$ term in the last line, when added to the
previous terms, then gives an interval containing
both lower and upper Riemann sums for 
$ \displaystyle \int_{a}^{b} f \ dx$.

Otherwise, $ A \le B$ is false, and since $a \le b$, $A$ and $B$ must overlap, so
we can decompose 
$$
\int_{a}^{b} f \ dx = \int_{a}^{c} f \ dx + \int_{c}^{b} f \ dx \ 
$$
where $a \le c \le b$ and $c \in A \cap B$.
Since $\underline{B} < \overline{A}$, $R = 0$ and the $RH$ term in the last line vanishes.
The result follows from
$[a,c] \subset [a,\overline{A}], [c,b] \subset [\underline{B},b]$ and inclusion monotonicity.
\enddemo 

\setlength{\parindent}{.2in}

The next several procedures define functions to be integrated.

The procedure {\bf Dx} is a straightforward interval extension of 
the expression for $dx/dy$ given in equation~\ref{eq:dx} in
Proposition~\ref{dx}.
It will be used to find the $x$-displacement from $(x_1,y_1)$ to $(x_2,y_2)$
by applying equation~\ref{eq:dx} to both  $\tau_i$ and~$\tau_o$.

Volumes under Delaunay surfaces are computed using the
volume integrand {\bf Dv}, which is obtained from {\bf Dx} as in
Proposition~\ref{dx}.
It is convenient to have {\bf Dv} calculate volume divided by $\pi$.
Since volumes are only compared to one another, we can
avoid the unnecessary step of multiplying by the nonrepresentable
constant $\pi$.
The volume element computation required that $dx/dy > 0$ for the volume to have
a positive sign.
The curve $\tau_o$ will sometimes start out in the negative $x$ direction,
when $\theta_1 < 60^\circ$, and after passing through a vertical tangency
switch to going in the positive $x$ direction.
Calculating the integral of {\bf Dv} in this case will give the volume under the part of
the nodoid generated by the arc going right minus the volume under the arc going left.
This is just what we want to compute the volume $V_o$
in the torus component in the case when such an overhang exists.

At local minima and maxima, {\bf Dx} and {\bf Dv}  are singular.
$ Y_{{\rm min}}$ and  $Y_{{\rm max}}$ are global variables which are accessed by
Dx$_{{\rm min}}$, Dx$_{{\rm max}}$, Dv$_{{\rm min}}$ and Dv$_{{\rm max}}$.
So near local extrema, we use a change of variables to
reparameterize the curve, as in Proposition~\ref{prop:dxcrit}.
Procedure {\bf Dx$_{{\rm min}}$} calculates
$dx/dz$ near a minimum and procedure {\bf  Dx$_{{\rm max}}$} near a maximum, where $z$ is given
in  Proposition~\ref{prop:dxcrit}.
We do an algebraic manipulation to make the denominator in the square root more computationally effective,
namely we replace the term 
$HY - 2 + H Y_{{\rm max}}$
in  procedure Dx$_{{\rm max}}$ by the equivalent expression
$  \displaystyle H (Y + {Y_{{\rm outer}}}_{{\rm min}} )  $.
Here $ \displaystyle {Y_{{\rm outer}}}_{{\rm min}}$ represents the minimum $y$-value of the nodoid which has a maximum at
$ Y_{{\rm max}}$.
The equivalence of the two expressions is derived using 
Proposition~\ref{prop.Delaunay} as follows:
$$
HY - 2 + HY_{{\rm max}} =  \frac{2}{Y_{{\rm max}} - {Y_{{\rm outer}}}_{{\rm min}}} (Y + Y_{{\rm max}}) - 2
=  H (Y + {Y_{{\rm outer}}}_{{\rm min}} ) .
$$
The volume calculation procedures {\bf Dv$_{{\rm min}}$, Dv$_{{\rm max}}$ }
give the formula for the volume integrand used near a critical point.
The volume calculations are applied only to graphs over the
$x$-axis.\medbreak

\setlength{\parindent}{0in}

{\tt
interval procedure Dx \\
input: $Y, H, F$    
\begin{tabbing}
Col 1 \= Col 2 \= Col 3 \kill
\>       $T := H Y^2 - F $  \\
\>       return $T / $Sqrt$( (2Y + T)(2Y - T))$  \\
\end{tabbing}
} 
 
{\tt
interval procedure Dx$_{{\rm min}}$   \\
input: $Z, H, F$    
\begin{tabbing}
Col 1 \= Col 2 \= Col 3 \= \kill
\> $Y := Y_{{\rm min}} + Z^2$  \\
\> $T := H Y^2 - F$  \\
\> return $2  T  / $Sqrt$( (2Y + T)(2 - H   Y_{{\rm min}} - Y H))$  \\
\end{tabbing}
}

{\tt
interval procedure Dx$_{{\rm max}}$   \\
input: $Z, H, F $   
\begin{tabbing}
Col 1 \= Col 2 \= Col 3 \= \kill
\> $Y := Y_{{\rm max}} - Z^2$ \\
\> $T := H Y^2 - F$ \\
\> return $2 T / $Sqrt$\displaystyle ((2Y + T) H (Y + {Y_{{\rm outer}}}_{{\rm min}} ) )$  \\
\end{tabbing}
}

{\tt
interval procedure Dv \\
input: $Y, H, F$    
\begin{tabbing}
Col 1 \= Col 2 \= Col 3 \= \kill
\>  return $ Y^2$ Dx($Y,H,F$)  \\
\end{tabbing}
}
\vglue-6pt
{\tt
interval procedure Dv$_{{\rm min}}$   \\
input: $Z, H, F$
\begin{tabbing}
Col 1 \= Col 2 \= Col 3 \= \kill
\> return $(Y_{{\rm min}} + Z^2)^2 $  Dx$_{{\rm min}}(Z,H,F)$  \\
\end{tabbing}
}
\vglue-6pt
{\tt
interval procedure Dv$_{{\rm max}}$   \\
input: $Z, H, F $
\begin{tabbing}
Col 1 \= Col 2 \= Col 3 \= \kill
\> return $ (Y_{{\rm max}} - Z^2)^2$ Dx$_{{\rm max}}(Z,H,F)$  \\
\end{tabbing}
}
\setlength{\parindent}{.2in}

\smallbreak
The next procedure is used for testing interval rectangles of torus bubbles to see if
they are potential minimizers and rejecting these intervals if they can be shown to
have a property which rules out any torus bubble contained in them.
The input is a pair of intervals $C_1$ and $H_o$, which
 give a range of possible values for $\cos \theta_1$ and $h_o$, respectively.
The interval of cosine values $C_1$ is used, rather than an interval
of $\theta_1$ values, to eliminate calls to trigonometric functions.
A variety of tests are performed based on the analysis of previous sections.
The various intervals do not
necessarily give sharp estimates,
so it is not {\it a priori} clear how many rectangles will be
rejected. As it turns out, the tests are sharp enough to reject all rectangles.\vglue6pt
\setlength{\parindent}{0in}

{\tt
Boolean procedure CheckRectangle \\
input: $C_1, H_o $
\begin{tabbing}
Col 1 \= Col 2 \= Col 3 \= \kill
1. \> if $1000 C_1 \ge 996$ and $ 5 H_o \le 1 $  \\
\> then return REJECT \\
\ \\

2.\> $H_i := 2 - H_o$ \\
\> $Y_1 := $ Sqrt$(1 - C_1^2)$ \\
\> $F_i := (H_i - 1)Y_1^2 - C_1 Y_1$Sqrt$(3)$ \\
\> $F_o := -F_i$ \\
\> $C_2 := [-C_1 \overline{\cup} C_1]  \cap  [-.5 \overline{\cup} .5] $  \\
\> $ T:= (2(H_i - 1)F_i + 3)/(3 + (H_i - 1)^2) - (1 - {C_1}^2)  $ \\
\> if $( {C_2}^2 +  T) \ne 1 $  return REJECT \\

\ \\

3. \> $ Y_2:= $ Sqrt$ ( T  \cap [0,1]  )$ \\
\> $\displaystyle C_2 := C_2 \cap \frac{ (H_i - 1)Y_2 - F_i/Y_2}{{\rm Sqrt}(3)}$ \\
\> if ${C_2}$ is empty then return REJECT \\
\ \\

4.\>  if $C_1 \le .5$ and $ H_o \le 1 - $Sqrt$(3) C_1/Y_1 $ \\
\> \> then return REJECT \\
\ \\

5.\>  if Width$(C_2) > .5$ then return NORESULT \\

\ \\
6.\> $W_{{\rm ends}} := (1-C_1)^2 (2+C_1)/3 + (1-C_2)^2 (2+C_2)/3$  \\
\> $Y_{{\rm min}} := - F_i / (1 +$ Sqrt$(1 + F_i H_i ))$  \\
\> $Y_{{\rm max}} := (1 + $Sqrt$(1 + F_o H_o)) / H_o$  \\
\> $Y := $ Compare$(C_1,$Sqrt$(3)/2,Y_{{\rm min}},Y_1)$ \\
\> if $Y H_i < -1$ then \\
\> \> begin \\
\> \> $ R :=  1/(- H_i - 1/Y)$  \\
\> \> $ W := 2.5 (Y + ($Sqrt$(3)/2)R) R^2$ \\
\> \> if $W < W_{{\rm ends}}$ then return REJECT  \\
\> \> end \\
 \ \\

7. \> $Y_{{\rm left}} := $ Sqrt$(F_o/H_o)$  \\
\> if ($Y_{{\rm min}} < Y_2$ and $Y_{{\rm left}} < Y_{{\rm max}}$) then $Y_{{\rm left}} := $ Max$ (Y_1,Y_{{\rm left}})$ \\
\> \> else return NORESULT  \\
\> $Y_4 := $ avgwt$(Y_{{\rm left}},Y_{{\rm max}},.5)$  \\
\> $Z_2 := $ Sqrt$(Y_{{\rm max}}-Y_2)$  \\
\> $Z_4 := - $Sqrt$(Y_{{\rm max}}-Y_4)$  \\
\> if $1000 C_1 \ge 998$ then \\
\> \> begin  \\
\> \> $ T := $ avgwt$(Y_1, Y_2 , 1/16)$   \\
\> \> $\Delta_i := (T - Y_1) (33/16) + $Integrate(Dx$(\cdot,H_i, F_i),T, Y_2)$ \\
\> \> $ \Delta_o := -(Y_{{\rm left}} - Y_1)$Sqrt$(3)$ \\
\> \> $ T := $  avgwt$(Y_{{\rm left}},Y_4,1/16)$   \\
\> \> $ \Delta_o := \Delta_o + $Integrate(Dx$(\cdot,H_o,F_o),T,Y_4)$\\
\> \> if $\Delta_o > \Delta_i$ then return REJECT   \\
\> \> $ \Delta_o := \Delta_o + $Integrate(Dx$_{{\rm max}}(\cdot,H_o,F_o),Z_4,Z_2)$  \\
\> \> if $\Delta_o > \Delta_i$ then return REJECT   \\
\> \> if $1 \in C_1$ then return NORESULT   \\
\> \> end  \\
 \ \\

8.\> verify ($C_1 \le $ Sqrt$(3)/2$ or $Y_1 \le Y_2$)  \\
\> $T := $ Sqrt$(Y_1-Y_{{\rm min}})$  \\
\> $Z_1 := $ Compare$(C_1,$Sqrt$(3)/2,-T,T)$  \\
\> $Z_3 := $ Sqrt$(Y_2 -Y_{{\rm min}})$  \\
\> $\Delta_i := $ Integrate(Dx$_{{\rm min}} (\cdot,H_i,F_i) ,Z_1,Z_3)$ \\
\> $\Delta_o := $Integrate(Dx$(\cdot,H_o,F_o),Y_1,Y_4)  $ \\
\> if $\Delta_i < \Delta_o$ then return REJECT \\
\> $\Delta_o :=  \Delta_o + $Integrate(Dx$_{{\rm max}} (\cdot,H_o,F_o) ,Z_4,Z_2)$ \\
\> if $\Delta_i \ne \Delta_o$ then return REJECT   \\
\ \\

9.  \> $W_{{\rm base}} := $ Integrate(Dv$_{{\rm min}} (\cdot,H_i,F_i),Z_1,Z_3) $  \\
\> $W_i := W_{{\rm ends}} + W_{{\rm base}}$   \\
\> $W_o := $Integrate(Dv$(\cdot,H_o,F_o),Y_1,Y_4) +$ \\
\> \> Integrate(Dv$_{{\rm max}}(\cdot,H_o,F_o),Z_4,Z_2) - W_{{\rm base}}$ \\
\>  if $W_i \ne W_o$ then return REJECT  \\
\> \> else return NORESULT \\
\end{tabbing}
}
\setlength{\parindent}{.2in}

\begin{center}
\BoxedEPSF{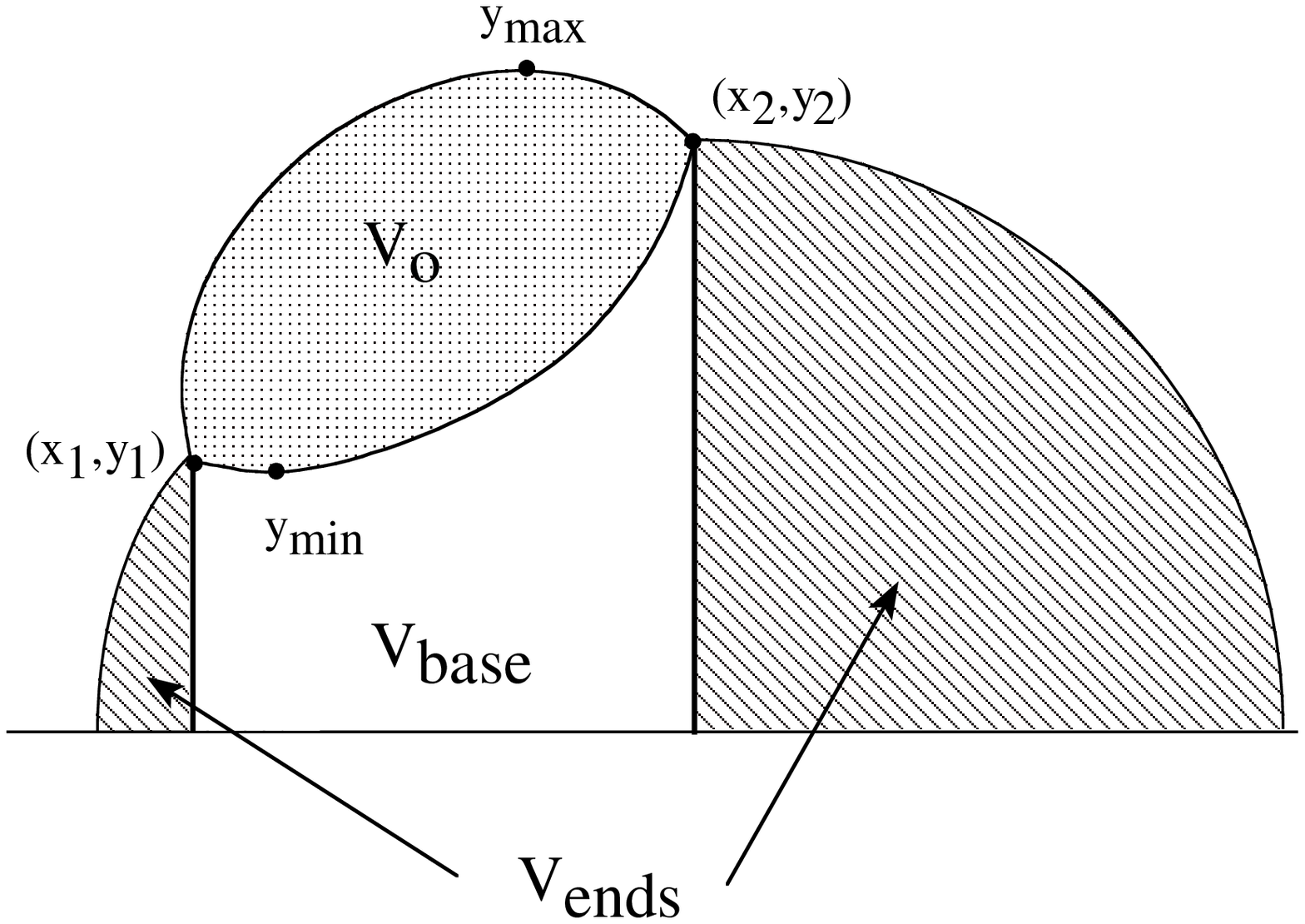 scaled 375}
\end{center}
\smallbreak
 {{Figure 17}.\enspace\ Compared volumes. If we can show that $V_o$ is not equal  to

$V_{{\rm base}} +V_{{\rm ends}}$
then we can discard the torus bubble as a potential equal

 volume
minimizer.}
 
\proclaim{Theorem} \label{CkRct} If {\bf CheckRectangle}($C_1,H_o$) returns {\rm REJECT} without
causing an {\rm IEEE} exception{\rm ,} then there is no area\/{\rm -}\/minimizing equal\/{\rm -}\/volume torus
bubble with $\cos \theta_1 \in C_1$ and $h_o \in H_o${\rm .}
\endproclaim

\demo{Proof} The procedure {\bf CheckRectangle} is fed a pair of intervals containing a range of values for 
$\cos \theta_1$ and $h_o$, these intervals being denoted by $C_1$ and $H_o$. We will call such
a range of input values an {\it input rectangle}. If every pair of values
$(\cos \theta_1, h_o) \in C_1 \times H_o$ in an input rectangle
can be ruled out as a possible minimizing torus bubble by some test,
then the procedure returns the value REJECT. Otherwise it returns the value NORESULT (which results in a subdivision
of the input rectangle into four smaller rectangles.) 
The tests applied to the rectangles of values for $\cos \theta_1$ and $h_o$ are based on the results developed in this
paper, as we will explain below.  At certain points we need to check whether some hypothesis
required in our tests are satisfied. We can do so by using a {\bf verify} statement.  This statement
will cause the program to terminate if it is not satisfied.
The proof proceeds by interpreting and justifying each step
of the algorithm.
\smallbreak
Step 1 rejects torus bubbles based on Proposition~\ref{0.2}. $C_1$ is a range of cosine values
for the angle $\theta_1$, so $1000 C_1 \ge 996$ corresponds to $\cos \theta_1 \ge .996$, which
implies that $\theta_1 \le 5.2^\circ$ and the proposition applies.  
\smallbreak
Step 2 first defines some values used in further steps.
The values of $F_o$ and $F_i$ are calculated using
the formulas in Proposition~\ref{force.values}. 
It then examines the possible values of $y_2 = \sin \theta_2$. 
Corollary~\ref{<120} implies $ -.5 \le \cos \theta_2 $, and
Lemma~\ref{lem.60} implies $ \cos \theta_2 \le .5 $.  Since we
assumed that $\theta_1 \le \theta_2$, we have $ \cos \theta_2 \le \cos \theta_1 $.  
Proposition~\ref{sum180} implies that $ -\cos \theta_1 \le \cos \theta_2 $.
Thus the values for $ \cos \theta_2$ arising from possible minimizing torus bubbles
lie in the intersection of the intervals $ [-.5 ,.5]$ and $[ -\cos \theta_1 ,\cos \theta_1]$.

The force equation for $T_i$,
$$
f_i(c) = (h_i-1) (1- c^2) - c \sqrt {3 (1- c^2) }  ,
$$
has at most two roots by
Lemma~\ref{2-1}, one of which is $c_1$.  The other one is $c_2$.   Letting
$y = \sqrt{1-c^2}$ and rearranging gives
\begin{eqnarray*}
f_i = (h_i-1) y^2 + cy \sqrt {3} \\
(f_i - (h_i-1) y^2)^2 = 3 y^2 (1 - y^2)   \ .
\end{eqnarray*}
This is a quadratic in $y^2$. The sum of the two roots in $y^2$ is 
$$
\frac{ 2(h_i-1)f_i + 3}{3 + (h_i-1)^2}
$$
and the known root is $1-{c_1}^2$, so subtraction gives the other root.
Thus the interval $T$ contains $y_2^2$, the square of the second root. 
We next add the interval ${C_2}^2$, containing ${\cos^2 \theta_2}$,
to the interval $T$ containing ${\sin^2 \theta_2}$ to see if they contain values
adding to one.
If not, there is no minimizing torus bubble in
the chosen input rectangle $C_1 \times H_o$, and this input rectangle is rejected. 
\smallbreak
Step 3 calculates the possible values of $x_2 = \cos \theta_2$ corresponding
to $y$-values whose square is in the interval $T$.
The only possible values for $y_2$ are in $[0,1]$, so we intersect
$T$ with $[0,1]$ before taking the square root of $T$ to get an
interval $ Y_2 $ which contains possible values of $y_2$.
Since 
$$
f_i = (h_i-1) (1- {c_2}^2) - {c_2} \sqrt {3 (1- {c_2}^2) }  
$$ 
by Proposition~\ref{force.values}, and $1- {c_2}^2 = y_2^2$,
we can solve to get an expression
for $c_2$ in terms of $h_i, f_i$ and $y_2$:
$$
C_2 = ( (H_i - 1)Y_2 - F_i/Y_2)/{\rm Sqrt}(3) \ .
$$
Note that this expression is chosen to give us $c_2$
with the correct sign. We intersect the resulting interval with the interval
$C_2$ computed in the previous step.
If the resulting interval is empty, there is no minimizing torus bubble in 
the chosen input rectangle and it is rejected. 
\smallbreak
Step 4 rejects torus bubbles based on Proposition~\ref{cotan}. It uses the fact that
$\cot(\theta_1) = c_1/y_1$. 
\smallbreak
Step 5 passes on a request to subdivide the rectangle if it is too wide.
It does not reject anything, so there is nothing to justify.
This step is included to speed the program's execution.
\smallbreak
Step 6 calculates rough bounds to check if the volumes
of the two regions can be equal.
We calculate volume divided by $\pi$ in the algorithm to avoid the needless
step of multiplying by the nonrepresentable constant $\pi$.
The interval $W_{{\rm ends}} = V_{{\rm ends}} / \pi $ is defined to be the volume of the solid of revolution
between the spherical ends and the $x$-axis, divided by $\pi$.
See Figure~17.
This solid is wholly contained in the ball component.
$Y_{{\rm min}}$ contains the $y$-coordinate of the local minimum of the
$\tau_i$, based on solving the force equation as in Part 7
of Proposition~\ref{torus.properties}.
This minimum is realized in the torus bubble if and only if
$\theta_1 \ge 30^\circ$.
Similarly $Y_{{\rm max}}$ contains the maximum height for $\tau_o$, as in
Part 8 of Proposition~\ref{torus.properties}. It is always achieved.
Step 6 goes on to reject potential torus bubbles if a calculation
shows that the torus component has less than half the total volume.
The volume of the ball component is bounded below by the volume
under the spherical caps, whose value is
given by $V_{{\rm ends}}$. The interval $Y$ is defined to ensure that
it contains the minimal $y$-coordinate of the torus component, which occurs either
in $Y_{{\rm min}}$ (if $\theta_1 \ge 30^\circ$) or in $Y_1$.
An upper bound $V$ for the torus component volume is calculated using
the formula derived in Lemma~\ref{trap1}.

\smallbreak
Step 7 tries to eliminate torus bubbles with small $\theta_1$, namely those
with $\cos {\theta_1} \ge .998$, a condition which forces
$\theta_1$ to be less than $3.7^{\circ}$ and allows
Lemma~\ref{lem:convex.unduloid} to be applied.
These small angles must be treated differently because the equations describing Delaunay
curves are not well behaved near $\theta_1 = 0$. 
Bounds on the $x$-displacement $x_2-x_1$ are calculated using
each of $\tau_i$ and $\tau_o$.
The interval $\Delta_i$ contains an upper bound for the width of 
$\tau_i$ and the interval $\Delta_o$ contains a lower bound for the width of
$\tau_o$.
If the resulting intervals satisfy $\Delta_i < \Delta_o$ then $\tau_i$ does not
travel as far to the right as $\tau_o$ and the
input rectangle can be rejected.

Step 7 begins by defining $Y_{{\rm left}}$ to be an interval containing
the $y$-value of the point where $\tau_o$ has a 
vertical tangency on the left as in Figures~18
and   19.
This value is computed by taking the square root of the nonnegative part of $F_o/H_o$.
Corollary~\ref{cor:signhf}
guarantees that this value is positive in a torus bubble,
so we are justified in truncating negative
values which may have appeared in the calculation before taking the square root.
The outer curve is always a nodoid by Corollary~\ref{T_o.is.nodoid},
and if extended to a vertical tangency, the
$y$-coordinate of that point is in the interval $Y_{{\rm left}}$. Replacing
$Y_{{\rm left}}$ with {\bf Max}($Y_{{\rm left}},Y_1$) assures that $Y_{{\rm left}}$
contains the $y$-coordinate of a
vertical tangency if there is one, and the initial point otherwise.
In particular it always contains the $y$-coordinate of the leftmost point.

The value $Y_4$ is used as an intermediate point for
chopping up the intervals of
integration for $\tau_o$ between its vertical and horizontal tangencies.
See Figures~18, ~19, ~20.
The limits  of integration $Z_2$ and $Z_4$ are defined for use when integrating near a maximum
of $\tau_o$, as in Proposition~\ref{prop:dxcrit}.
The square roots taken in these calculations
are justified in discarding the negative part of any interval, because for any
particular torus bubble the quantity whose root is being taken is nonnegative.

\figin{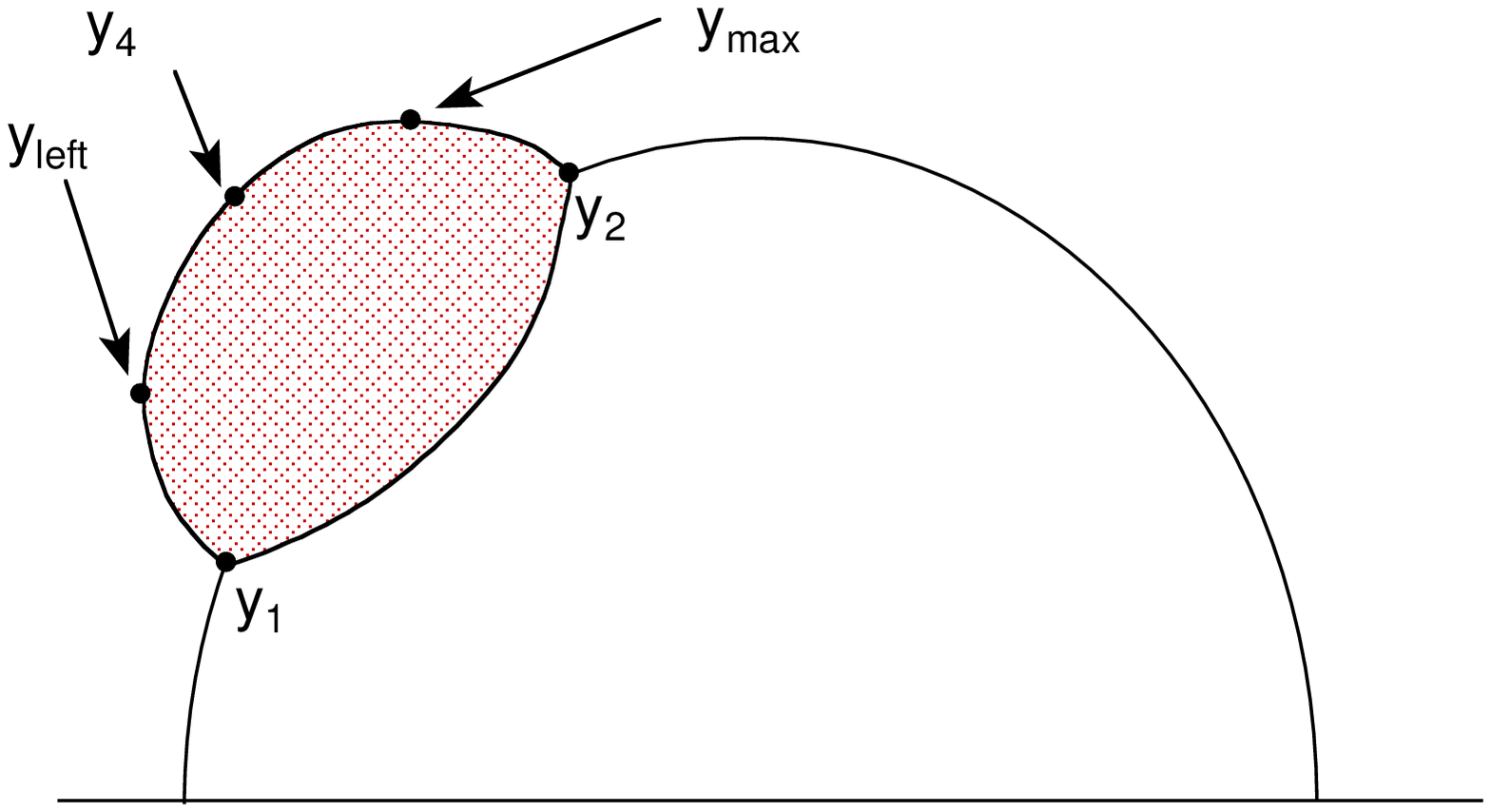}{340}
 \centerline{{Figure 18}. A double bubble configuration for $0 < \theta_1 \le 30$.} 
\figin{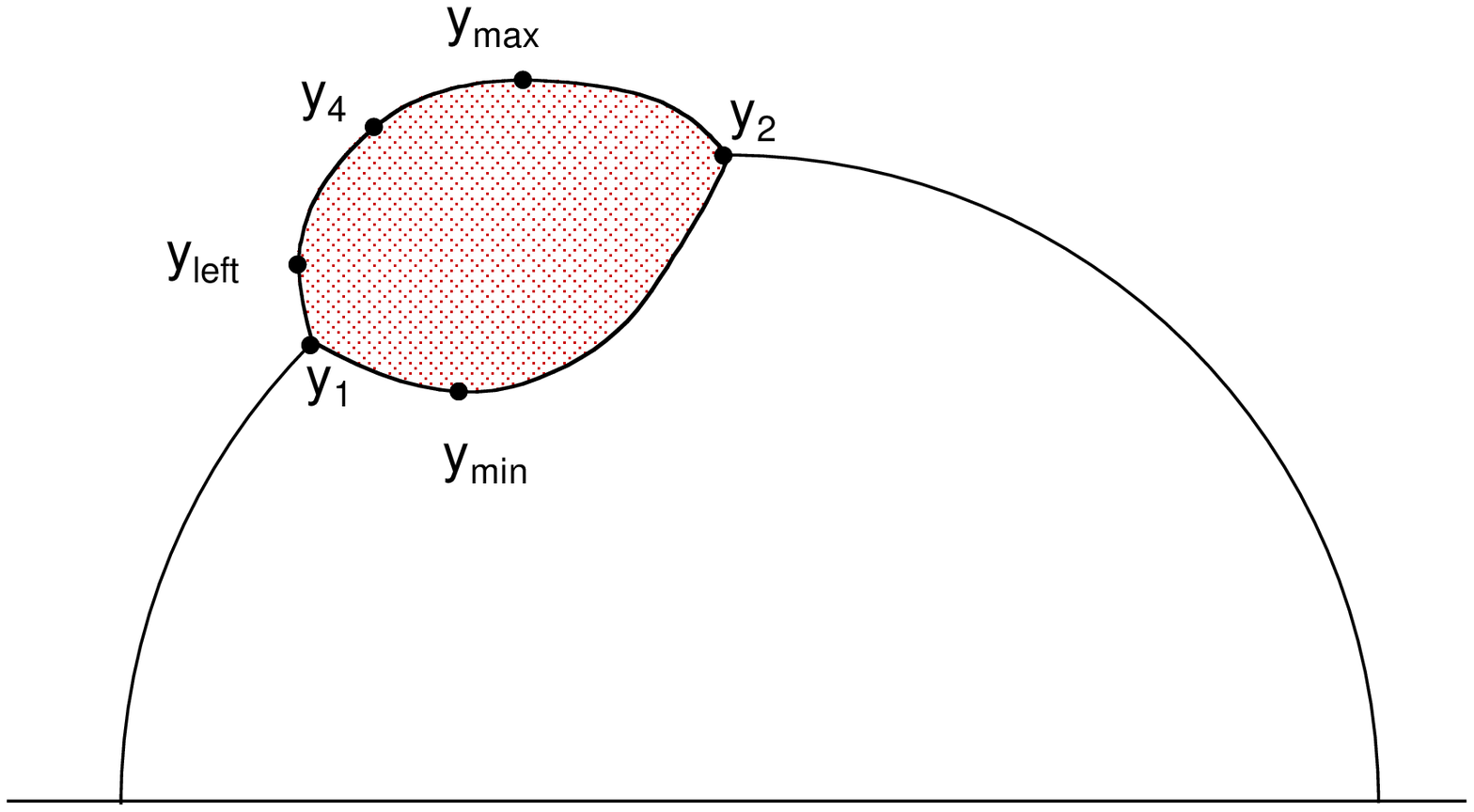}{340}
 \centerline{{Figure 19}. A double bubble configuration for $30 \le  \theta_1 \le 60$.} 
\figin{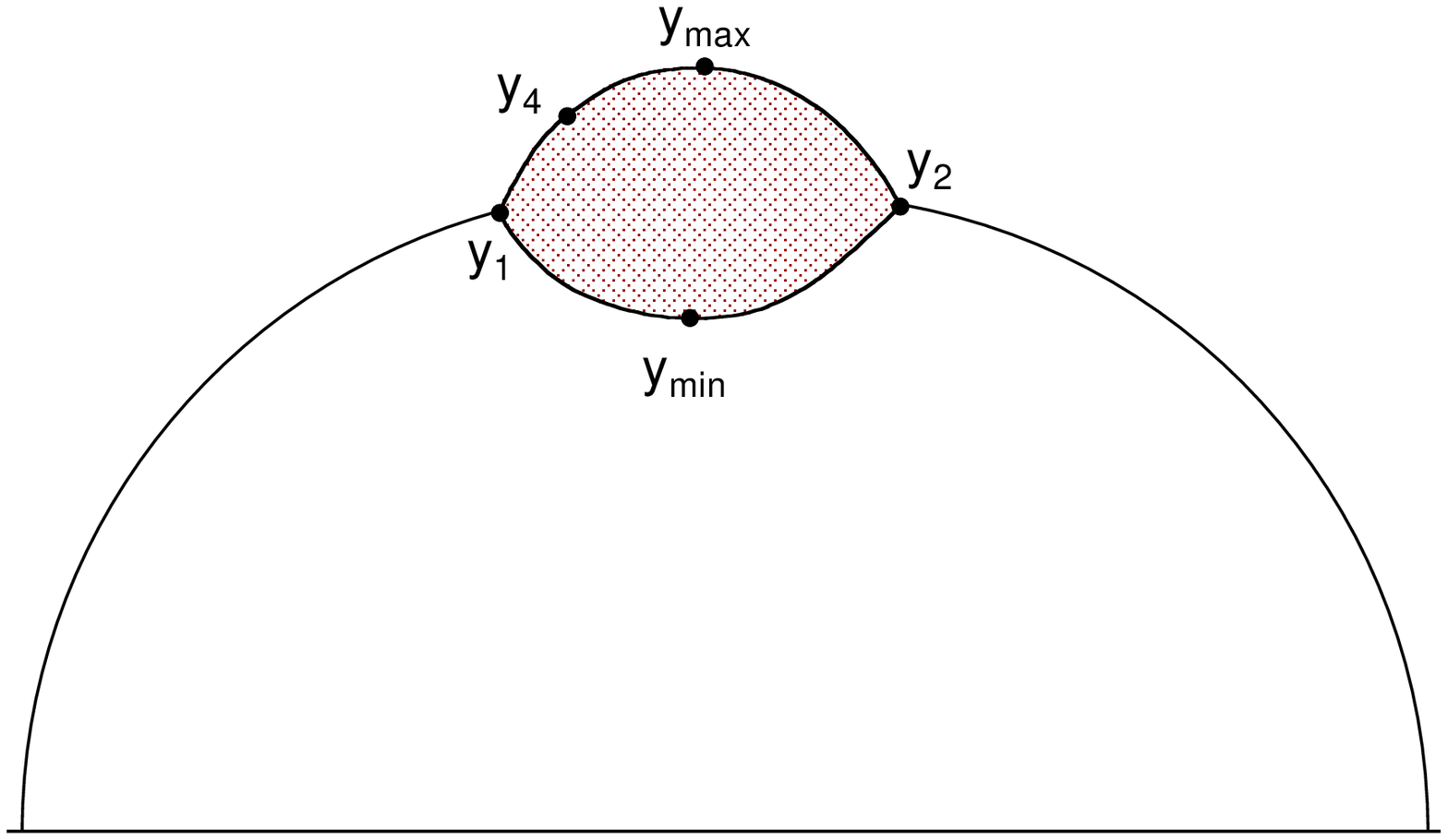}{340}
 \centerline{{Figure 20}. A double bubble configuration for $60 \le \theta_1  \le 90$.} 
\medbreak

An upper bound for the $x$-displacement of $\tau_i$ as it goes from
$y_1$ to $y_2$ is given by $\Delta_i$, and
a lower bound for the $x$-displacement of $\tau_o$ as it goes from 
$y_1$ to $y_2$  is given by $\Delta_o$.
If $\Delta_i < \Delta_o$, then a torus bubble cannot occur in the given 
input rectangle.

We use the Intermediate Value Theorem to estimate $\Delta_i$
up to height $ \displaystyle y_1 + \frac{y_2 - y_1}{16}$
and $\Delta_o$ up to height $y_{{\rm left}}$. 
The estimates for $\Delta_i$ and $\Delta_o$ near the $x$-axis are
justified by Proposition~\ref{prop.Delaunay} in the case of $\tau_o$ and by
Lemma~\ref{lem:convex.unduloid} in the case of $\tau_i$,
as shown in Figure~21.

\figin{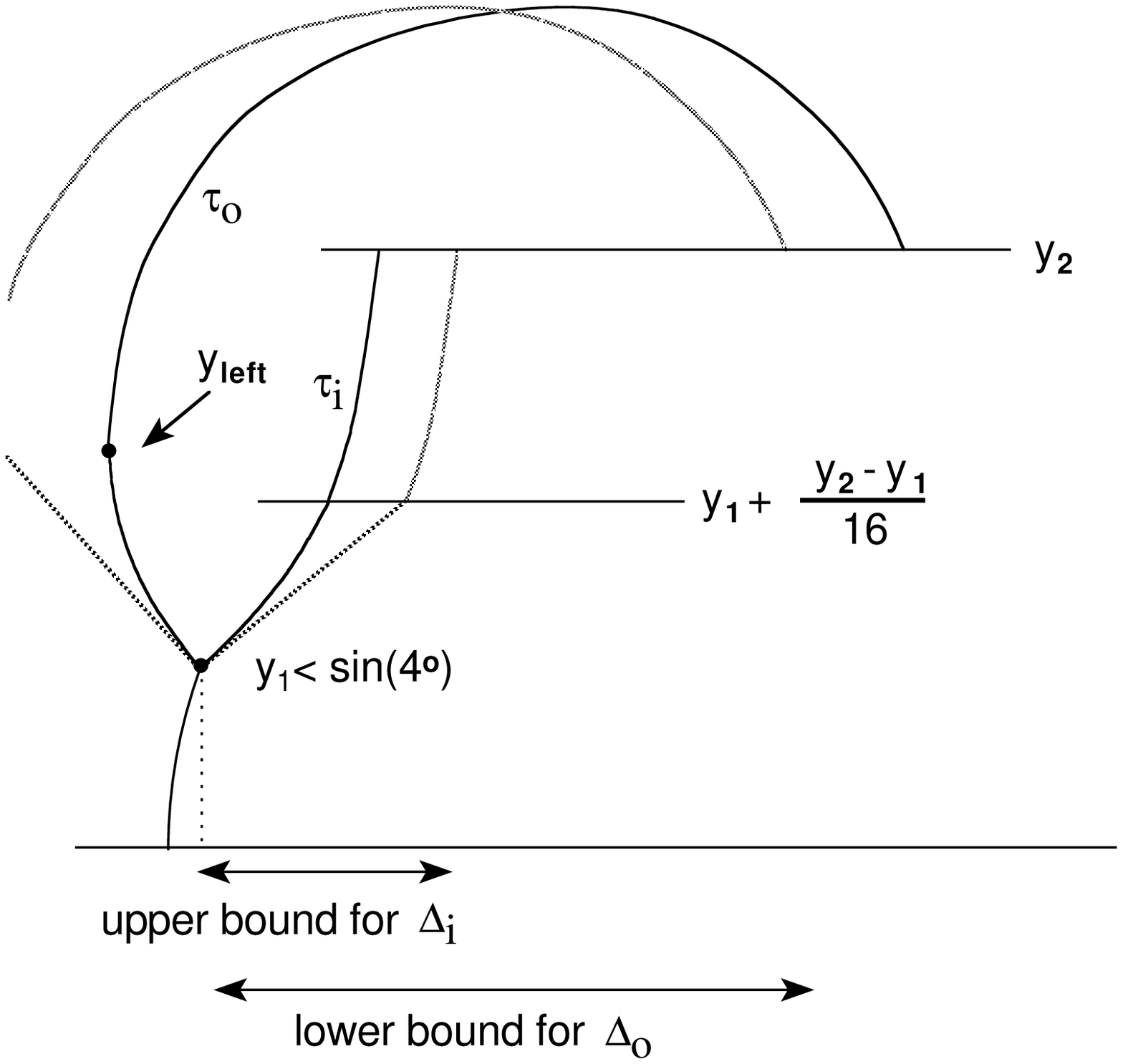}{400}%

 {{Figure 21}. \enspace Estimating the widths of $\tau_i$ from above and $\tau_o$ from below

with
 $\theta_1$ small.} 
\medbreak

For $\theta_1 < 4^\circ$, $\tau_i$ is
strictly increasing, so that the $x$-displacement is equal to the integral of {\bf Dx}
from $y_1$ to $y_2$. This can be broken up into an integral from $Y_1$ to $T$ and an
integral from $T$ to $Y_2$, where $T$ is set to be one sixteenth of the way from $Y_1$ to
$Y_2$. The values in $T$ are all less than 1/2, so that Lemma~\ref{lem:convex.unduloid} applies
and the slope of $\tau_i$ below height $T$ is
smallest at its initial point ($x_1,y_1$). It follows that 
$$
dx/dy \le \cot {(30 - \theta_1) } \le \cot  {26^\circ} < 33/16
$$
on this interval and the integral from $Y_1$ to $T$ is bounded
above by $(33/16)\break  (T-Y_1)$. Then
$\Delta_i = (T - Y_1) (33/16) + $Integrate(Dx$(\cdot,H_i, F_i),T, Y_2)$
gives an upper bound for the width
of $\tau_i$ as it increases from $Y_1$ to $Y_2$.

With $\tau_o$ we first find an upper bound on the leftward displacement between $y_1$ and
$y_{{\rm left}}$ using a similar method. In this case the slope has
absolute value larger than $\tan 30^\circ = 1/ \sqrt 3$
and $(Y_{{\rm left}} - Y_1)$Sqrt$(3)$ gives an
upper bound for the leftward displacement to height $y_{{\rm left}}$.
An integral is calculated
which gives a lower bound for the rightwards displacement of $\tau_o$ between $y_{{\rm left}}$ and 
$y_2$, by integrating $dx/dy$ from a point strictly above
$y_{{\rm left}}$ to $y_2$.
This integral is calculated in two pieces, the first to $Y_4$ and the second from
$Y_4$ through the local maximum to $Y_2$. 
A lower bound for $\Delta_o$ is gotten by taking the upper
bound for the leftwards displacement previously calculated and subtracting it from the sum of
these integrals. 
If $\Delta_o > \Delta_i$ the input rectangle is rejected.

This step can reject torus bubbles even if the interval 
$C_1$ contains the value~1, corresponding to $\theta_1 = 0$.
If this step test fails to reject a range of potential torus bubbles, we check whether
1 is in the interval $C_1$. If it is, then the next steps will not help, since
they are ineffective for $\theta_1 = 0 $. So if 1 is in $C_1$ we return NORESULT
and send the current input rectangle back for subdivision.

Step 8 calculates intervals for the $x$-displacement of $\tau_i$ and $\tau_o$ in the case
where the angle is larger than those treated in the previous step, namely $\cos \theta_1 \le .996$.
In this step $\Delta_o$ is an
interval containing the $x$-displacement of the outer curve and
 $\Delta_i$ an interval containing the $x$-displacement of $\tau_i$.
In contrast to Step 7, they do not contain bounds for these displacements,
but the actual widths of $\tau_i$ and $\tau_o$.
Thus we can reject the input rectangle in this step if 
$\Delta_i$ and $\Delta_o$ do not overlap.

The interval $\Delta_o$ is calculated by integrating {\bf Dx} from $Y_1$ to an intermediate
$Y_4$ and then using {\bf Dx$_{{\rm max}}$} to integrate from $Y_4$ across a maximum
to $Y_2$.
If the width of $\Delta_o$ is larger than $20$ then the
procedure returns NORESULT, calling for a subdivision of the input rectangle of intervals.
This is done because otherwise large amounts of computational time are wasted in inefficient calculations.

The calculation of $\Delta_i$ is complicated somewhat because we need
to deal with
both the case where $\theta_1 \ge 30^{\circ}$ and the case where  $\theta_1 \le 30^{\circ}$.
Only in the former case, where $C_1 \le $ Sqrt$(3)/2$, does $\tau_i$ achieve an interior minimum.
In that case, an integral along $\tau_i$ can be calculated by integrating 
{\bf  Dx$_{{\rm min}}$} from $Z_1$ to $Z_3$.
If an interior minimum is not
achieved, we can still integrate {\bf  Dx$_{{\rm min}}$} to calculate the
$x$-displacement from $Y_1$ to $Y_2$, but we need to reverse the sign of
$Z_1$ as in the last statement of  Proposition~\ref{prop:dxcrit}.
Our integral formulas for $x$-displacement need to have that $Y_1 \le Y_2$ in this case,
so we check  with a verify statement that either this holds or there is an interior minimum.
The {\bf Compare} statement ensures the correct sign is given to $Z_1$ unless
$C_1$ overlaps the interval {\bf Sqrt}$(3)/2$. 
$Z_1$ equals {\bf -Sqrt}$(Y_1-Y_{{\rm min}})$ if $\tau_i$ achieves an interior minimum and
{\bf Sqrt}$(Y_1-Y_{{\rm min}})$ otherwise.
The appropriate sign is given to $Z_1$ using a {\bf Compare} statement.
If $C_1$ overlaps the interval {\bf Sqrt}$(3)/2$ we can't tell which sign is correct, so
the {\bf Compare} statement returns an interval containing both $Z_1$ and $-Z_1$,
and the integral then contains the $x$-displacements for both cases.
If $\Delta_i$ and $\Delta_o$ are
disjoint, then the curves' final points cannot coincide, and the input rectangle
is discarded. 

Step 9 calculates volumes divided by $\pi$, as in Figure~17.
We again calculate volume divided by $\pi$ in the algorithm to avoid
an unnecessary multiplication of all volumes by $\pi$.
$V_{base}$ is the volume surrounded by $T_i$ and $W_{base} = V_{base}/ \pi$.
Adding $W_{{\rm base}}$ to $W_{{\rm ends}}$ gives $W_i = V_i/ \pi)$.
The calculation of $V_o$ is slightly more complicated because the
overhang on the left involves a subtraction, if there is an overhang.
Such an overhang occurs
when $\theta_1 < 60^\circ$, as shown in Figures~18 and
19.
The formulas are set up so that volume is counted with a
negative sign when $\tau_o$ is oriented to the left  and positive sign when $\tau_o$ is oriented
to the right, as discussed in  Proposition~\ref{dx}.
This is just what is required to calculate the volume of the overhang.
The volume $V_{{\rm base}}$ is subtracted from the volume inside $T_o$ to get $V_o$,
the volume of the torus component.
Step 9 then rejects an input rectangle if the values
$W_o = V_o/ \pi$ and $W_i = V_i/ \pi$ are unequal. \enddemo

\setlength{\parindent}{.2in}

Now that we have tests that are able to reject certain ranges of torus bubbles, we apply these
tests to see whether they in fact reject all the possibilities. Since we don't know
in advance how fat the ranges can be, we initially feed in a rectangle of
values with $0 \le \sin \theta_1 \le 1, 0 \le h_o \le 10$.
We recursively subdivide further as necessary.
We monitor IEEE exception flags at this level, so that no
computation is trusted if it raised an exception. \\

\setlength{\parindent}{0in}

{\tt
Boolean procedure DivideAndCheckRectangle \\
input: $Y_1, H_i$
\begin{tabbing}
Col 1 \= Col 2 \= Col 3 \= \kill
\> $ C_1 := $ Sqrt$(1 - Y_1^2)$ \\
\> call CheckRectangle($C_1,H_o$) \\
\> if result is REJECT then return SUCCESS \\
\> split $Y_1$ in half, into $Y_{1a}, Y_{1b}$ \\
\> split $ H_o$ in half, into $H_{oa}$, $H_{ob}$ \\
\> call DivideAndCheckRectangle($Y_{1a},H_{oa}$) \\
\> call DivideAndCheckRectangle($Y_{1a},H_{ob}$) \\
\> call DivideAndCheckRectangle($Y_{1b},H_{oa}$) \\
\> call DivideAndCheckRectangle($Y_{1b},H_{ob}$) \\
\> return SUCCESS \\
\end{tabbing}
}

{\tt
Main program
\begin{tabbing}
Col 1 \= Col 2 \= Col 3 \= \kill
\> begin   \\
\> clear exceptions \\
\> \> begin \\
\> \> $Y_1 := [0,1] $ \\
\> \> $H_o := [0,10]$ \\
\> \> Call DivideAndCheckRectangle($Y_1,H_o$) \\
\> \> end \\
\> verify no exceptions raised \\
\> Print ``All torus bubbles rejected.'' \\
\> end 
\end{tabbing}
}

\setlength{\parindent}{.2in}

\proclaim{Theorem} \label{no.minimizing.tori} The algorithm described in {\bf Main}{\rm ,} if run to
completion without causing an exception{\rm ,} shows that no area minimizing torus bubble
can enclose equal volumes{\rm .}
\endproclaim

\demo{Proof} By Proposition~\ref{determined}, any minimizing torus bubble is determined by
$\theta_1$ and $h_o$, where $\theta_1$ is the angle of first arc and $h_o$ is the
mean curvature of~$T_o$. It suffices to consider $0 \le \theta_1 \le 90^\circ$
by Proposition~\ref{sum180}, and $0 \le h_o \le 10$ by Proposition~\ref{param.bounds} and
Proposition~\ref{h>-8}. 

{\bf Main} calls {\bf DivideAndCheckRectangle} with a rectangle containing the parameters of
all possible torus bubbles.  {\bf DivideAndCheckRectangle} subdivides this rectangle,
depending on the success of its calculations. To avoid unnecessary use of trigonometric functions in the
computation, the space of angles is parameterized by a variable $y_1 = \sin \theta_1$. Thus the
intervals $Y_1$ cover the interval $[0,1]$, in 1-1 correspondence with
$ 0 \le \theta_1 \le 90^\circ$. {\bf DivideAndCheckRectangle} begins by using $Y_1$ to compute
$C_1$, an interval corresponding to values of the cosine of an interval of $\theta_1$ angles.  It
then passes the $C_1$ and $H_o$ intervals to {\bf CheckRectangle} where they are tested as
previously described.

{\bf Main} can only finish if {\bf CheckRectangle} rejects all rectangles. By
Theorem~\ref{CkRct}, torus bubbles rejected by {\bf CheckRectangle} are not minimizers.
\enddemo

We implemented the algorithm described above in C++, though
we could equally well have
used any language which allows for control of rounding
and follows the IEEE 754 standard, or any similar precise
prescription for floating point arithmetic.
We ran code implementing this algorithm
on a variety of machines, various Wintel, Sun, HP, SGI, UNIX and
Linux platforms, with identical results. The C++
code with embedded instructions is available at the web site \cite{hassweb}.
The code takes about 10 seconds to run on a fast (1999) PC.
The input rectangle was recursively subdivided up to ten times, resulting in a total
of 15,016 rectangles examined. The calculation involved a total
of 51,256 integrals.  

From examining the output of the program, we can see the reasons for
which various rectangle were rejected. Of the 15,016 rectangles
examined, Step 1 rejected  31, Step 2 rejected 1140, Steps 3, 4 and 5
rejected none, Step 6 rejected 3664, Step 7 rejected 541, Step 8 
rejected 6691 and Step 9 rejected 2949.
In Figure~22 we indicated which steps were used to
reject torus bubbles in which input rectangles.
This chart offers some insight into why torus bubbles fail to be minimizers.
However there is some arbitrariness to this, as doing the steps in 
a different order would have produced different data.

In running the program no exceptions were raised,
and the program returned the statement
``All torus bubbles rejected.''
As a consequence we obtain a proof of our main result.

\specialnumber{1} \proclaim{Theorem}  The unique surface of least area enclosing two equal volumes
in $\RR^3$ is a double bubble{\rm .}
\endproclaim

{\it Proof}. Theorem~\ref{no.minimizing.tori} eliminates the possibility of a minimizing torus
bubble. Theorem~\ref{two.types} implies that the only other possibility is a double bubble.

\begin{center}
\BoxedEPSF{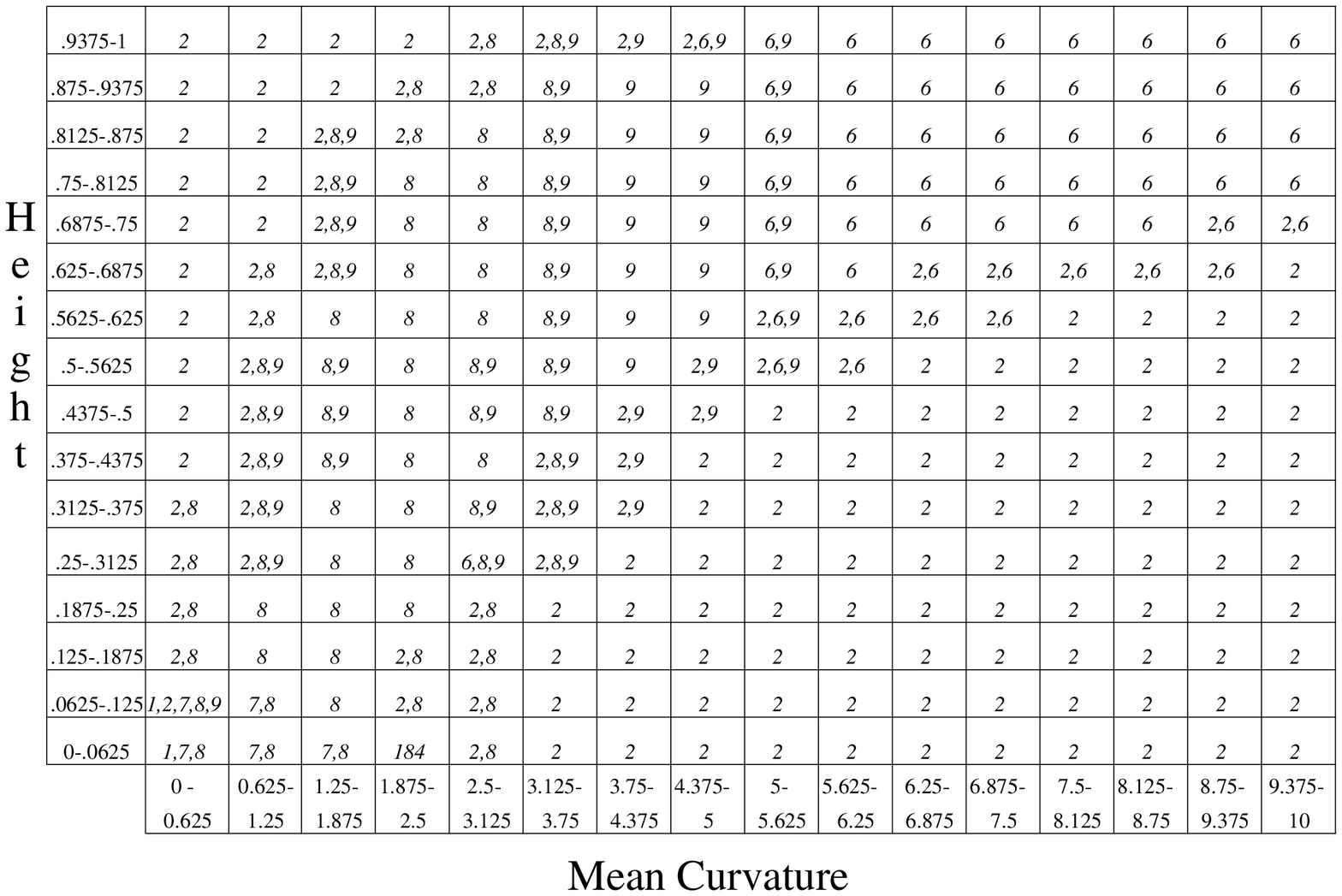 scaled 675}
\end{center}
\medbreak\noindent{{Figure 22}. Chart showing which tests were used in each input rectangle. In some 
rectangles there were several subdivisions, and
several tests were used. Steps 3, 4 and 5 were never used to reject input rectangles. }
\medbreak

\demo{Acknowledgements}
We are indebted to Frank Morgan for introducing us to this problem, to John Sullivan for
generously providing us with computer graphics of torus and double bubbles, and to F. Morgan, M.
Hutchings,\break W. Kahan and W. Rossman for helpful discussions.

\end{document}